\documentclass[onefignum,onetabnum]{siamonline190516}

\usepackage{lipsum}
\usepackage{amsfonts}
\usepackage{graphicx}
\usepackage{algorithmic}
\usepackage{ntheorem}
\usepackage{subfig}
\usepackage{enumitem}
\def \Q {\mathbb{Q}}
\usepackage{dsfont}	   

\newsiamremark{remark}{Remark}
\newsiamremark{hypothesis}{Hypothesis}
\newsiamremark{Def}{Definition}
\newsiamremark{lem}{Lemma}
\newsiamremark{pro}{Proposition}
\crefname{hypothesis}{Hypothesis}{Hypotheses}
\newsiamthm{claim}{Claim}

\title{On the conservation properties in multiple scale 
coupling and simulation for Darcy flow with 
hyperbolic-transport in complex flows}

\author{Eduardo Abreu\thanks{Department of Applied Mathematics, 
Institute of Mathematics, Statistics and Scientific Computing (IMECC)
University of Campinas (UNICAMP), Campinas 13083-970, SP, Brazil
  (\email{eabreu@ime.unicamp.br}).}
\and Ciro Diaz\thanks{Postdoctoral Fellow at Department of
Mathematics at Ryerson University, 350 Victoria Street
Toronto, Canada (\email{cirojdiaz@gmail.com}).}
\and Juan Galvis\thanks{Departamento de Matem\'aticas, Universidad Nacional de Colombia, Carrera 45 No. 26-85, Edificio Uriel Gutierr\'ez, Bogot\'a D.C., Colombia
  (\email{jcgalvisa@unal.edu.co}).}
\and Jonh P{\'e}rez\thanks{Metropolitan Institute of Technology 
(ITM) - University Institution, Calle 73 No 76A - 354 V\'{\i}a 
al Volador, Medell\'{\i}n, Colombia (\email{jhonperez@itm.edu.co}).}}

\begin{document}

\maketitle

\begin{abstract}
We present and discuss a novel approach to deal with 
conservation properties for the simulation of 
nonlinear complex porous media flows in the 
presence of: 1) multiscale heterogeneity 
structures appearing in the elliptic-pressure-velocity  
and in the rock geology model, and 2)  multiscale wave 
structures resulting from shock waves and rarefaction 
interactions from the nonlinear 
hyperbolic-transport model. For the pressure-velocity 
Darcy flow problem, we revisit a recent high-order 
and volumetric residual-based 
Lagrange multipliers saddle point problem to impose 
local mass conservation on convex polygons. We clarify 
and improve conservation properties on applications.
For the hyperbolic-transport problem we introduce a new
locally conservative Lagrangian-Eulerian finite volume 
method. For the purpose of this work, we recast our 
method within the Crandall and Majda treatment of the 
stability and convergence properties of conservation-form, 
monotone difference, in which the scheme converges to 
the physical weak solution satisfying the entropy 
condition. This multiscale coupling approach was 
applied to several nontrivial examples to show that 
we are computing qualitatively correct reference
solutions. 
We combine these procedures for the simulation 
of the fundamental two-phase flow problem with 
high-contrast multiscale porous medium, but 
recalling state-of-the-art paradigms on the of notion 
of solution in related multiscale applications.
This is a first step to deal with out-of-reach multiscale 
systems with traditional techniques. 
We provide robust numerical examples for verifying 
the theory and illustrating the capabilities of 
the approach being presented.
\end{abstract}



\section{Introduction}
In this paper, we are concerned with modeling, simulation 
and numerical analysis for approximate solutions in 
multiscale nonlinear Partial Differential Equations (PDE) 
related to highly complex systems. 
The large number of papers published 
in recent years is indicative of the relevance of the 
foundations of multiscale approach. It is a good measure
of the breadth and of the vitality of the 
area, and therefore, calling the multiscale modeling 
and simulation research community for new ideas and 
innovative approaches. In this work we give an overview 
of recent approaches
for the accurate and efficient simulation of complex
porous media flows and we also present some new results. 
We summarize below the main aspects of our work:

\begin{itemize}

\item We revisited a novel volumetric {\it locally 
conservative} and residual-based Lagrange multipliers 
saddle point reformulation of high-order 
 methods and present numerical results with 
realistic high-contrast multiscale coefficients. This is applied  to 
a second-order elliptic problem  
($\nabla \cdot \left[-K({\bf x}) \Lambda(S) \nabla p\right] = q(\textbf{x})$ )
instead of a traditional {\it first-order mixed 
formulations}. We clarifying and simplifying the 
presentation of its conservative properties.

\item We introduce a new robust and accurate forward 
tracking Lagrangian-Eulerian scheme for hyperbolic problems. 
Our method is carefully de\mbox{sign}ed to deal with multiscale wave structures 
resulting from shock wave interactions coming 
from abstract nonlinear {\it hyperbolic 
conservation laws}  such as
$w_t + \mbox{div} \mathcal{F}(t,{\bf x},w(t,{\bf x})) = 0$ 
 with fluxes of 
the form $\mathcal{F}(x,t,w) = {\bf v}(x, t)f(w(t,{\bf x}))$.
This formulation includes many problems of physical interest and,
in particular, we consider the case of complex multiscale 
flow in porous media (scalar and systems cases).

\item We improve the interpretation of the construction 
of numerically stable Lagrangian-Eulerian {\it no 
flow surface region} in two-space dimensions previously 
presented and analyzed in \cite{EAJP19} for 
one-dimensional balance and conservation laws.

\item We present a new approach with accurate 
multiscale resolution for two-phase flows. 
The method is able to handle multiscale rock 
geology in the elliptic-pressure-velocity system. 
Moreover, the multiscale wave structures from wave 
interactions present in  the hyperbolic-transport model 
seems to be properly simulated according to our numerical results.

\item We present numerical results with realistic 
high-contrast two-dimensional multiscale coefficients 
based on the 10th SPE Comparative Solution Project 
(SPE10). We address numerical issues of multiscale 
resolution and conservation properties.

\item We survey recent results on both: novel 
deterministic and probabilitc multiscale modeling 
paradigms and issues of mesh resolution inadequacy 
in multiscale complex systems.

\end{itemize}

The new elliptic solver \cite{AGP18a,AGP18b} is a 
general tool for imposing local conservation for 
high-order methods to deal with multiscale 
permeabilities. In particular is also applicable
to the Generalized Multiscale Finite Element Method (GMsFEM); see 
\cite{AGP19, GMSFEM, EGW2009jcp, PRESHO2016376, GMsFEMrecent} 
and references therein. On the other hand, the new 
Lagrangian-Eulerian method seems to be a promissing 
general approach to capture nonlinear wave 
interactions linked to multiscale behavior of 
interactions of waves in a wide range of models 
and, in particular, for systems (see 
\cite{EAJP19,APS17a,APS17b,APS18,AMPRB19,AAP19}). 
Numerical results show that our combined multiscale 
approach provides an accurate and robust procedure 
to study phenomena which couple distinct length 
or time scales. Our novel hyperbolic Lagrangian-Eulerian 
solver circumvents the use of adaptive and/or mesh 
generation. Indeed, the new method is also free of 
Riemann solvers. For Darcy flow, the method captures 
fine-scale effects using multiscale finite element techniques.
Finally, it is worth mentioning that for 
the purpose of this work we use Cartesian grids.
Thus convergence and error 
analysis reduces essentially to a one-dimensional 
problem and retains convergence result of approximated
solutions to the entropy weak solution by 
recalling  \cite{EGH95a,EGH95b,CCH99,CRMAJ80}. 
This is to say that our approach fits comfortably within
the classical theory. On the other hand, we provide 
 convergence and error analysis on triangular 
grids for hyperbolic conservation laws in another 
forthcoming work \cite{AAP19}. Moreover, the 
results in \cite{EAJP19}, strongly suggest that 
the {\it no flow surface region} encapsulates the domain 
of dependence even in the case of hyperbolic systems.

A reliable multiscale prediction of oil-water 
two-phase flow model through complex porous media 
requires the development and validation of coupling 
models for flow across length scales related to 
elliptic Darcy-pressure-velocity and hyperbolic 
saturation-conservation problems. Therefore, in 
this work we  introduce  a novel computational 
approach for multiscale computing of the oil-water 
flow model system. 

Let $\Omega$ be a domain in $\mathbf{R^2}$.
For simplicity 
of presentation, we consider two-phase flow in a 
viscous-dominated flow, and we neglect the effects 
of gravity\footnote{ In order to 
focus on ideas related to multiscale complexities  
we neglect  gravity terms in our formulation.
We mention also that in the case 
$\Omega \subset  \mathbf{R^3}$ gravity terms  
should be included.}, compressibility and capillarity and 
set the porosity equal to a constant, in which case it
has been scaled out by a change of the time variable. 
Thus, the fundamental differential multiscale system 
used to describe water-oil, incompressible, 
immmiscible displacement is given  by the {\it saturation} and 
{\it pressure} equations (see, e.g., 
\cite{CHM06}):
\begin{equation}
\frac{\partial S }{\partial t} + \nabla \cdot 
({\bf v} \, \mathcal{F}(S)) = 0, 
\quad \mbox{in } {\bf x} \in \Omega, 
\qquad 0 < t \leq T,
\label{hyper}
\end{equation}
\begin{equation}
\nabla \cdot {\bf v} = \nabla \cdot \left[-K({\bf x}) 
\Lambda(S) \nabla p\right] = q(\textbf{x}), \quad  
\mbox{in } {\bf x} \in \Omega, \qquad 0 < t \leq T.
\label{pres}
\end{equation}
Here, $S$ is the water saturation (water and oil 
saturation sum up to one), ${\bf x}$ is the 
absolute permeability, ${\bf v}$ is the total 
Darcy velocity,  $p$ is the pressure, $\Lambda$ 
is the total mobility and $\mathcal{F}$ is the 
fractional flow of water. We also have 
$\mathcal{F}=\Lambda_i / \Lambda$ with 
$\Lambda_i = k_i/\mu_i$ being  the phase mobility, 
$k_i$ the relative permeability and $\mu_i$ 
viscosity of phase $i$, respectively.

We end up with a multiscale system of two coupled 
nonlinear partial differential equations 
(\ref{hyper})-(\ref{pres}) exhibiting multiscale features
in both sides: the complex multiscale heterogeneity 
structures from rock geology appearing in the 
elliptic-pressure-velocity model (\ref{pres}) as well as 
multiscale wave structures resulting from shock 
wave interactions from the hyperbolic-transport 
model  (\ref{hyper}). For concreteness, we  complete 
equations (\ref{hyper}) and (\ref{pres}) with 
appropriate initial and boundary conditions 
aiming simulations in a slab geometry domain 
$\Omega$ and time interval $T =\left[t_0, t_f\right]$. 
 We consider the coupling system 
(\ref{hyper})-(\ref{pres}) in a two-dimensional 
rectangular slab domain $\Omega = (0, L_x) 
\times (0, L_y )$, with the boundary conditions
\begin{equation}
\begin{array}{ll}
{\bf v} \cdot {\bf n} = -q, & \text{on } x = 0, \\
p = 0, & \text{on } x = L_x, \\
{\bf v} \cdot {\bf n} = 0, & \text{on } y = 0, L_y ,
\end{array}
\label{bce}
\end{equation}
where ${\bf n}$ is the outward-pointing unit 
normal vector to $\partial \Omega$, and a 
uniform initial condition 
\begin{equation}
S(x, 0) = S_0.
\label{ice}
\end{equation}
The initial-boundary conditions (\ref{bce})-(\ref{ice}) 
simulate a left to-right waterflood. Water is 
injected uniformly (at a constant rate $q$) 
through the left vertical boundary $\Gamma_{inlet}$ 
($x = 0$) of $\Omega$ for all simulation time 
$t \in \left[t^0, t_f\right]$, no flow conditions are 
imposed along the horizontal boundaries 
$\Gamma_{no-flux}$ ($y = 0, L_y$), and fluid is 
produced from a well kept at constant (zero) 
pressure at the right vertical boundary 
$\Gamma_{outlet}$ ($x = L_x$).

Equations (\ref{hyper})-(\ref{ice}) can be viewed 
as multiscale water-oil Riemann-Goursat water 
injection problem with discontinuous flux function 
as closely related to  \cite{BADM15} (see 
\cite{SV03,EP10,EP09,EP07,EP05}), in which the 
hyperbolic scheme should be able to handle the 
main difficulty of the multiscale problem that 
consists in taking into account the jump discontinuities 
of the flux; see also \cite{EGH95a,EGH95b,ABFL19,EA14,ACP16,ACP17} 
and references cited therein. Another motivation 
is the occurrence of models such as (\ref{hyper})-(\ref{ice}) 
in numerous engineering problems, mainly the case 
of systems \cite{EA14,ACP16,ACP17,HFVS03,HFVS05,ASFM10,ASFM14,AMPZ02,DMBP01,ABFL19}. 
In addition, there are a number of prototype relevant 
models of hyperbolic conservation laws with discontinuous 
flux functions in {\it oil trapping phenomenon} 
\cite{BACC13,BACC14,CC10B,KEF99,SMJJ10,STIA12}, 
{\it a Whitham model of car traffic flow on a highway} 
\cite{HHNR95,LVFEB02} and a model of {\it continuous 
sedimentation in ideal clarifier-thickener units} 
\cite{BKRT04}, see also \cite{BBK04a,KT04}. Hyperbolic 
conservation laws with discontinuous flux functions 
also arise in sedimentation processes \cite{SD96}, 
in radar shape-from-shading problems \cite{DNO99} 
and also as building blocks in numerical methods 
for Hamilton-Jacobi equations \cite{KHKNHR02}. 
Hyperbolic conservation laws of the form (\ref{hyper}) 
with discontinuous (in $(t, {\bf x})$) flux function 
$\mathcal{H}({t,{\bf x},S})\equiv\{f(S)\left[-K({\bf x}) \Lambda(S) \nabla p\right]\}$ 
attracted much attention in the recent past, because 
of the difficulties of adaptation of the classical 
Kruzhkov approach developed for the smooth case, 
due in part to the presence of several different 
entropy solutions with same initial data. In the 
context of Buckley-Leverett equations as in 
(\ref{pres})-(\ref{ice}), each notion of solution 
is uniquely determined by the choice of a connection, 
which is made unique at the interface by a proper 
choice of an entropy solution from  many possible classes of 
entropy solutions, \cite{AMVG05, BACC13, BADM15}. Nonclassical solutions appear 
for Buckley-Leverett models with gravity and 
discontinuous flux functions \cite{CC10B,BACC13} 
as well as for three-phase flow problems with 
continuous flux functions, where recent and very 
relevant results can be found in 
\cite{ASFM10,ASFM14,AMPZ02,DMBP01} for a solution 
of Riemann problems and in \cite{HFVS03,HFVS05} 
concerning well-posedness. In \cite{AKR11}, a 
theory of $L^1$-dissipative solvers for scalar 
conservation laws with discontinuous flux was 
proposed, which is based on the corresponding 
$L^{1}$ contractive semigroups, some of which 
reflect different multiscale physical applications. 
In \cite{AKR11,BA15,CCP15} a number of the existing 
admissibility (or entropy) conditions are revisited
 and the so-called germs that underly these 
conditions are identified. The seminal survey 
article \cite{BA15} (see also 
\cite{RSS19,FB06,TTZT06,FBBP98,TT01,NNK76,HCCP19,KT04}) 
of recent developments helps to better understand 
the issue of admissibility of solutions in relation 
with specific modeling assumptions.

By looking at model problem (\ref{hyper})-(\ref{ice}) 
as a generalized problem linked to systems of 
conservation laws in several dimensions, we recall that some 
authors has advocated entropy measure-valued 
solutions, first proposed by DiPerna \cite{RJD85,RDAM87}, 
as the appropriate solution paradigm for systems of 
conservation laws; see  \cite{FKMT17,FLMW19} and references 
cited therein. In particular, these authors have 
presented some numerical evidence that state-of-the-art 
numerical schemes may not converge to an entropy 
solution of systems of conservation laws as the 
mesh is refined. Accordingly to \cite{FKMT17,FLMW19} 
 this has 
been attributed to the emergence of turbulence-like 
structures at smaller and smaller scales upon 
mesh refinement. See also \cite{GCJG19, GGJG12, LGS08}.

In \cite{FLMW19}, the authors point out that 
{\it intermittency is widely accepted to be a 
characteristic of turbulent flows \cite{UFTu95}}. 
It is believed that intermittency stems from the 
fact that turbulent solutions do not scale exactly 
as in the Kolmogorov hypothesis. On the other hand, 
in \cite{GCJG19} the authors offer a unified 
framework in which is possible to establish mathematical 
existence theories as well as a very innovative 
idea for the interpretation of numerical solutions 
through the identification of a function space in 
which convergence should take place. In a more 
general setting, the issue of multiscale modeling 
and simulation of chaotic mixing of distinct 
fluids has been addressed in  \cite{GCJG19,GGJG12,LGS08}. 
They mention that acceleration driven turbulent mixing 
is a classical hydrodynamic instability. See also 
\cite{FFFP03} to the case of two-phase flows.

We also mention also the very recent works 
\cite{GCJG19, GGJG12, LGS08,FKMT17,FLMW19} on 
nonlinear multiscale problems like (\ref{hyper})-(\ref{ice}) 
in which resolution of multiscale turbulent-like 
behavior is better resolved under mesh refinement. 
In fact, structures at smaller and smaller scales 
are formed as the mesh is refined to account for the 
complex heterogeneity structures 
from rock geology appearing in the elliptic model 
as well as multiscale wave structures resulting 
from shock wave interactions from the underlying 
hyperbolic problem. This is done such that solutions satisfy
a proper family of entropy inequalities 
\cite{BADM15}. Results in \cite{CLLS09,ECCL15} 
have demonstrated that entropy solutions may 
not be unique. In the direction from both rigorous 
mathematical analysis and numerical analysis, 
ingenious difficulties stem from the lack of 
regularity of solutions.

A better comprehension of multiscale fluid flow 
in subsurface is very hard,  challenging and 
undoubtedly still of current events. Multiscale 
sciences cuts across all of science from fluid 
dynamics to biology,  from meteorology to 
material science  and from physics to chemistry among 
many other directions. Multiscale issues are also central 
in subsurface flows ranging from complex geologic 
media to several time scales linked to the 
compositional and black-oil modeling fluid flow 
in oil reservoirs as well as several scale 
aspects of groundwater flow and related transport 
systems. Understanding the multiscale properties 
of subsurface flows is a major problem of modern 
approaches predicting groundwater level changes and 
predictive technologies  in petroleum reservoir. In this regard, 
many innovative techniques have been reported as 
such local-global upscaling approach 
\cite{LDHL16,GKD06,MDD06,CDGW03,WEH02,LJD91}, 
multiscale methods 
\cite{JEA04,APWY07,THXW97,JLT03,LWT08,OMKAL16a,OMKAL16b,GASP18,AHPV13}, 
model order reduction techniques 
\cite{HFEG19,JJLD17,MCLD10,SLW19,CLVW}, Two-scale 
homogenization theory 
\cite{FSYFY19,GA92,ABJ06,AB84,JKO12,ATMV16,ADH90}; 
see papers for a survey on recent development of 
multiscale computing and modeling approach 
\cite{SLW19,ADG19,THT19,GCLW19,GKNSVL19,VCEK19,HCCZC19}.

Despite the efforts of many researchers, no universal 
or unified multiscale modeling and methods for fluid 
flow through naturally complex geologic reservoirs 
has been achieved so far (see, e.g., 
\cite{ADG19,VCEK19,CNP18,RKNP18,EA14,SLW19,JJLD17,HFEG19,OMKAL16a,APWY07,LJD91,WEH02,LDHL16}). 
Under appropriate simplification assumptions 
compositional and black-oil models can be further 
simplified for the fundamental multiscale 
two-phase immiscible displacement with no mass 
transfer between phases (this is often appropriate 
in models describing displacements at the length 
scales associated with reservoir simulation grid 
blocks, see, e.g., \cite{LDHL16}). In addition, 
even for the two-phase case,  the multiscale modeling is 
very challenging in the presence  fractures and 
barriers for flow in porous medium  and their impact on the closure and 
constitutive relations as such multiscale relative 
permabilities and pressure difference (capillary 
pressure); see \cite{BACC13} and references cited 
therein for an interesting study on vanishing 
capillarity solutions of Buckley-Leverett equation 
with gravity in two-rock's medium (see also 
\cite{RKNP18,CNP18,BGJMP17} for multiscale modeling 
of Richards's equation and two-phase under 
non-equilibrium effects). However, in the case 
of scalar water-oil two-phase model a global-pressure 
formulation and Kirchhoff’s transform are not adequate 
when considering non-equilibrium effects as such 
hysteresis in the relative permeabilities \cite{ABFL19} 
and in the capillary pressure \cite{ABSH01} the 
issue of global pressure formulation for three-phase 
is not straightforward (see \cite{EA14,ABFL19,CAFM16} 
and reference cited therein for a detailed multiscale 
modeling for three-phase flow problem). Moreover, 
the case of discontinuous capillary pressure induced 
by mulsitcale modeling of fractures and barriers 
in the three-phase flow with gravity is hard and 
very intricate \cite{EA14,EAHYP12,CAFM16}. The 
degeneracy for three-phase and two-phase flows is 
also delicate \cite{EADC13,RPK08}. Altogether the 
fundamental multiscale water-oil model 
(\ref{hyper})-(\ref{ice}) is also useful both for 
describing some real cases (e.g., dead-oil systems) 
and for developing and studying numerical solution 
procedures. In this works we consider the multiscale 
approximation associated with reservoir multiscale 
simulation along with coupling 
techniques for elliptic (Darcy-pressure-velocity) 
and hyperbolic (conservation-saturation-transport) 
problems as pursued in this work.

The paper is organized as follows. In Section 
\ref{Ch2}, we construct an embedded high-order model 
for second-orde elliptic problems with local and 
global mass conservation, clarifying and simplifying 
the presentation of its conservative properties in 
lines as introduced in 
\cite{AGP18a,AGP18a,PRESHO2016376,AbGaDiSa2018,AGP19}. 
Next, we construct new a locally conservative 
Lagrangian-Eulerian method for hyperbolic-transport 
with focus on a novel approach for conservation 
properties of the no flow surface region for hyperbolic 
conservation laws in Section \ref{theLE2d}. In 
Section \ref{Coupling}, we present and discuss a 
the oupling conservative finite element method 
for Darcy flow problem with a locally conservative 
Lagrangian-Eulerian method for hyperbolic-transport, 
along with a set of representative computational 
results. In Section \ref{cpfinal}, a summary with 
concluding remarks and perspectives for future 
work are highlitghted. 

\section{Elliptic problem and mass conservation}\label{Ch2}

Many porous media related practical problems (scalar and systems) 
lead to the numerical approximation of the Buckley-Leverett type-models
given by the highly-nonlinear multiscale problem like (\ref{hyper})-(\ref{ice}). 
For the approximation of the pressure field given by the 
pressure-velocity Darcy-elliptic fundamental multiscale problem (\ref{pres}).
We use the method de\mbox{sign}ed and analyzed in \cite{AGP18a,PRESHO2016376};
see also \cite{AbGaDiSa2018,AGP19,GMSFEM}.  
We now recall this high-order conservative FEM formulation, 
clarifying and simplifying the presentation of its conservative
properties with focus on the conservation properties in multiple 
scale coupling and simulation for Darcy flow with 
hyperbolic-transport in complex flows. A general interpretation 
of this methodology goes as follows: given a Ritz approximation 
of the pressure and a computational procedure to obtain fluxes 
we can then, formulate an constrained (to local flux conservation) 
minimization problem to obtain approximated solution that satisfy 
local conservation properties in the form of average fluxes on 
local regions. See \cite{PRESHO2016376,AbGaDiSa2018,AGP19}.
\subsection{Imposing local mass conservation}
{\color{black} We follow the presentation in \cite{AGP18a,PRESHO2016376}.}
Denote $H^1_D (\Omega)$  the space of functions in 
$H^1(\Omega)$ which vanish on $\partial \Omega_D$. 
The variational formulation of problem (\ref{pres}) 
is to find $p\in H^1_D(\Omega)$ such that
\begin{equation}\label{eq:problem}
a(p,v)=F(v)   \quad \mbox{ for all } v\in H_D^1(\Omega),
\end{equation}
where $v\in H^1_D(\Omega)$ and the bilinear form $a$ is 
defined by
\begin{equation}\label{eq:def:a}
\vspace{-1.5mm}
a(p,v)=\int_\Omega 
 \Lambda (x)\nabla p(x)\nabla v(x) dx, 
\end{equation}
and the functional $F$ is defined by 
\vspace{-1.5mm}
$
{F(v)=\int_\Omega q({\bf x})v(x)dx.}
$\\

 Problem (\ref{pres}) is equivalent 
to the minimization problem: Find $p\in H^1_D(\Omega)$
such that 
\begin{equation}\label{eq:unconstraint}
\vspace{-1.5mm}
p  =\arg\min_{v\in H^1_D(\Omega)} \mathcal{J}(v),  
\quad \mbox{ where } \quad 
\mathcal{J}(v)=\frac{1}{2}a(v,v)-F(v).
\end{equation}
In the IMPES approach, for each time step the 
mobility can be thought as a function of position 
and reads simply as $\Lambda (x)$. Thus, in 
order to consider a general formulation for 
porous media applications we let $\Lambda$ be a 
$2\times 2$ matrix with entries in $L^\infty(\Omega)$ 
in Problem (\ref{pres}) to be almost everywhere 
symmetric positive definite matrix with eigenvalues 
bounded uniformly from below by a positive 
constant. 

\begin{figure}[ht!]
\centering
\includegraphics[scale=.5]{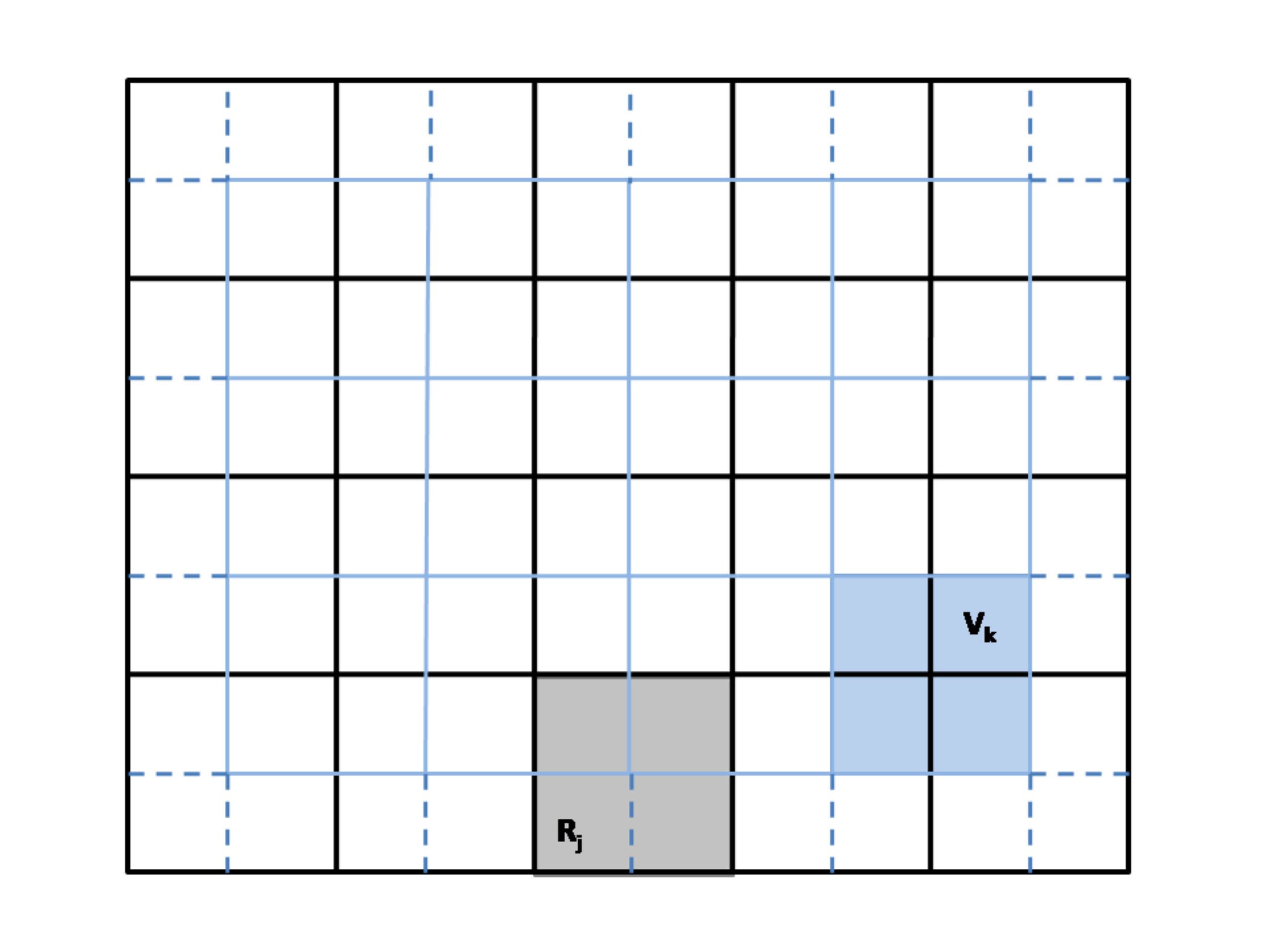}
\caption{Example of regular mesh made of squares and its dual 
triangulation.}
\label{Fig:Primal_Dual}
\end{figure}

In order to deal with mass conservation properties we 
follow the method introduced in \cite{PRESHO2016376,AbGaDiSa2018}. 
Let $ \mathcal{T}_h =\{  R_j \}_{j=1}^{N_h}$ be a 
primal mesh made of  elements that are triangles or 
squares. Here $N_h$ is 
the number of elements of the triangulation. We also 
have a dual mesh $\mathcal{T}_h^*=\{V_i \}_{i=1}^{N_h^*}$ 
where the elements are called control volumes and 
$N_h^*$ is the number of such volumes. In general 
it is selected one control volume $V_i$ per vertex 
of the primal not in $\partial\Omega_D$. In case 
$|\partial\Omega_D| = 0$, $N_h^*$ is the total 
number of vertices of the primal triangulation 
including the vertices on $\partial\Omega$. Figure  
\ref{Fig:Primal_Dual} illustrate primal and dual 
meshes made  of squares where $\partial\Omega_D = \partial\Omega$, 
in this case $N^*_h$ is equal to the number of 
interior vertices of the primal triangulation.

If $q\in L^2$ we have that solving (\ref{pres}) 
is equivalent to: Find {$p \in H_{\mbox{\small div},\Lambda}^1(\Omega)$} 
and such that
\begin{equation}\label{eq:problem-with-restriction2}
p=\arg\min_{v\in \mathcal{W}} \mathcal{J}(v),
\end{equation}
where the subset of functions that satisfies the mass 
conservation restrictions is defined by
\[
\mathcal{W}=\left\{ v\in H_0^1(\Omega),: \int_{\partial T}
- \Lambda  \nabla v\cdot\mathbf{n} =\int_{T}q \quad 
\text{ for all } T\in \mathcal{T}^*_h \right\}.
\]
and
\[
H_{\mbox{\small div},\Lambda}^1(\Omega) = \{  v \in H_D^1(\Omega) 
: \Lambda\nabla v \in H(\mbox{div},\Omega)\}
\]
with norm $\| v \|^2_{H_{\mbox{\small div},\Lambda}^1(\Omega)} = 
\|\Lambda\nabla v\cdot\nabla v \|^2_{L^2(\Omega)} 
+ \|\mbox{div}(\Lambda \nabla v) \|^2_{L^2(\Omega)}$.\\

Let  $M^h = \mathbb{Q}^0(\mathcal{T}^*_h)$ be the space 
of piecewise constant functions on the dual mesh $\mathcal{T}_h^*$. 
The Lagrange multiplier formulation of problem 
(\ref{eq:problem-with-restriction2}) can be written as:
Find $p  \in H_{\mbox{\small div},\Lambda}^1(\Omega)$  and 
$\lambda \in M^h$ that solves,
\begin{equation}
\label{eq:discrete-problem-with-restriction-lag}
\arg\max_{\mu\in \color{black} M^h}\min_{v \in\color{blue} H_{\mbox{\small div},\Lambda}^1(\Omega) } 
\mathcal{J}(v)-
(\overline{a}(p,\mu)-G(\mu)).
\end{equation}
Here, the total flux bilinear form 
$\overline{a}: H_{\mbox{\small div},\Lambda}^1(\Omega),\times {\color{black} M^h}\to \mathbb{R}$ is defined by
\begin{equation}\label{eq:def:overline-a}
\overline{a}(v,\mu)=\sum_{i=1}^{\color{black}N_h^*}
\mu_i\int_{\partial V_{i}} \Lambda  \nabla v\cdot 
\mathbf{n} \quad \mbox{ for all 
}  v\in H_{\mbox{\small div},\Lambda}^1(\Omega), \mbox{ and } \mu\in M^h .
\end{equation}
The functional 
$G:{\color{black} M^h} \to \mathbb{R}$ is defined 
by\[G(\mu)=\sum_{i=1}^{\color{black}N_h^*} 
\mu_i\int_{V_{i }}q\quad  \mbox{ for all } \mu \in  M^h .\]
The first order conditions of the min-max problem above  
give the following saddle point problem:
Find $p$ with $p \in H_{\mbox{\small div},\Lambda}^1(\Omega),$ 
and $\lambda\in M^h $ that solves,
\begin{equation}\label{eq:saddlepoint}
\begin{array}{llr}
a(p,v)+\overline{a}(v,\lambda)&=F(v) &\mbox{ for all } 
v\in H_{\mbox{\small div},\Lambda}^1(\Omega),, \\
\overline{a}(p,\mu) &=G(\mu)& \mbox{ for all } 
\mu \in M^h.\\
\end{array}
\end{equation}
{\color{black}For the analysis of this formulation see \cite{AGP18a}.}
Recall that we have introduced a primal mesh 
$\mathcal{T}_h =\{  R_j \}_{j=1}^{N_h}$ {made of 
elements that are triangles or squares}. Here $N_n$ is 
the number of elements of the triangulation. We also 
have given a dual mesh $ \mathcal{T}_h^*=\{  V_k \}_{k=1}^{N_h}$ 
where the elements are called control volumes. Figure 
\ref{Fig:Primal_Dual} illustrate a primal and dual mesh made 
of squares. See for instance \cite{MR2168342, PRESHO2016376,AbGaDiSa2018}. 

Let us consider $ P^h = \mathbb{Q}^r (\tau_h)$ the 
space of continuous polynomial functions   
of degree $r$ on each element of the primal mesh, and 
{let \color{black}$P^h_0 $ be the space the functions  in  $P^h$ that vanish in $\partial 
\Omega $}. Let 
$M^h = \mathbb{Q}^0(\tau^*_h)$ be the space of piece 
constant functions on the dual mesh $\tau_h^*$. 
{\color{black} For more details on the constructions see \cite{AGP18a,PRESHO2016376}.}
The discrete version of (\ref{eq:saddlepoint}) is to find 
$p^h \in {\color{black}P^h_0}$ { and} $\lambda \in M^h$ 
such that
\begin{align}\label{cela}
&a(p^h,v^h) + \overline{a}(v^h,\lambda^h) =  
     F(v^h) && \mbox{ for all } v^h \in {\color{black}P^h_0}\\
&\overline{a}(p^h,\mu^h) = 
      G(\mu^h) && \mbox{ for all } \mu^h \in M^h.
\end{align}

Let $\left\lbrace \varphi_i \right\rbrace$ be the 
standard basis of $P^h$. We define the matrix
\begin{equation}\label{eq:CoersMat}
A = \left[ a_{i,j} \right] \quad \mbox{ where } 
a_{ij} = \int_{\Omega} \Lambda \nabla\varphi_i 
\cdot \nabla\varphi_j.
\end{equation}

Note that $A$ is the finite element stiffness matrix 
corresponding to finite element space $P^h$. Introduce 
also the matrix
\begin{equation}\label{eq:restMat}
B = \left[ \overline{a}_{k,j} \right]  
\quad \mbox{ where } \overline{a}_{ij} = 
\int_{\partial V_k} \Lambda \nabla\varphi_j \cdot \textbf{n}.
\end{equation}

With this notation, the matrix form of the discrete 
saddle point problem is given by,
\begin{equation}\label{eq:saddle}
\mathcal{A}U^h = \left[
\begin{array}{cc}
A&B^T \\
B&O
\end{array} \right]
\left[ \begin{array}{c}
u^h\\
\lambda^h
\end{array}
\right] = \left[ \begin{array}{c}
f\\
g
\end{array}
\right]
\end{equation}
where vectors $f$ and $g$ are defined by, 

\begin{align}
f = \left[f_i\right] \quad \mbox{ with} \, f_{i} = 
\int_{\Omega} q \cdot \varphi_i && g = 
[g_i]_{i=1}^{M_f} = \int_{V_k} q.
\label{eq:fontVect}
\end{align}
{\color{black}For the analysis of the continuous and discrete problem  and corresponding error estimates see \cite{AGP18a}.}
As an important remark we mention  the optimality of the 
approximation error in both $H^1$ and $L^2$ norms.  From 
the analysis in \cite{AbGaDiSa2018,AGP19} we know that 
$u^h$ is an optimal approximation in the $H^1$ semi-norm 
(which becomes a norm when restricted to appropriate 
subspace). We also found in  \cite{AbGaDiSa2018} that, 
in each control volume, the $u^h+\lambda^h$  offers 
optimal approximation in the $L^2$ norm.  We real that 
in many approaches imposing conservation properties of 
second order formulations leads to a non-optimal $L^2$ 
approximation. For more details on this formulation in 
porous media applications; see
\cite{ADG19,PRESHO2016376,AbGaDiSa2018,AGP18a,AGP18b}.

\subsection{Numerical tests}

Let us illustrate the two main features of the 
HOCFEM method which are high-order approximation 
rate and conservation of mass.

\subsubsection{Homogeneous medium}

We consider the  Equation (\ref{pres})
with \textbf{$\Omega = \left[0,1\right]\times \left[0,1\right]$} in 
homogeneous medium ($\Lambda(x) \equiv 1$), to be discretized 
on a regular mesh made of $2^M \times 2^M$ squares. 
The dual mesh is  constructed by joining the centers 
of the elements of the  primal mesh as in Figure 
\ref{Fig:Primal_Dual}. We consider homogeneous 
Dirichlet's boundary conditions in $\partial\Omega$  
and construct the example (\ref{eq:Example:1}) by 
fixing the solution $p$ and computing a source term 
$q$ so we can compare the numerical solution with 
the exact solution
\begin{equation}
\begin{aligned}
q(x,y) &= 2\pi (\cos(\pi x) \sin(\pi y) 
- 3\sin(\pi x) \cos(\pi y) 
+ \pi \sin(\pi x) \sin(\pi y) (-x+3y)), \\
p(x,y) &= \sin(\pi x) \sin(\pi y) (-x+3y).
\end{aligned}
\label{eq:Example:1}
\end{equation}

\begin{figure}[ht!]
\centering 
\includegraphics[scale=.3, trim = 4cm 0 0 0]{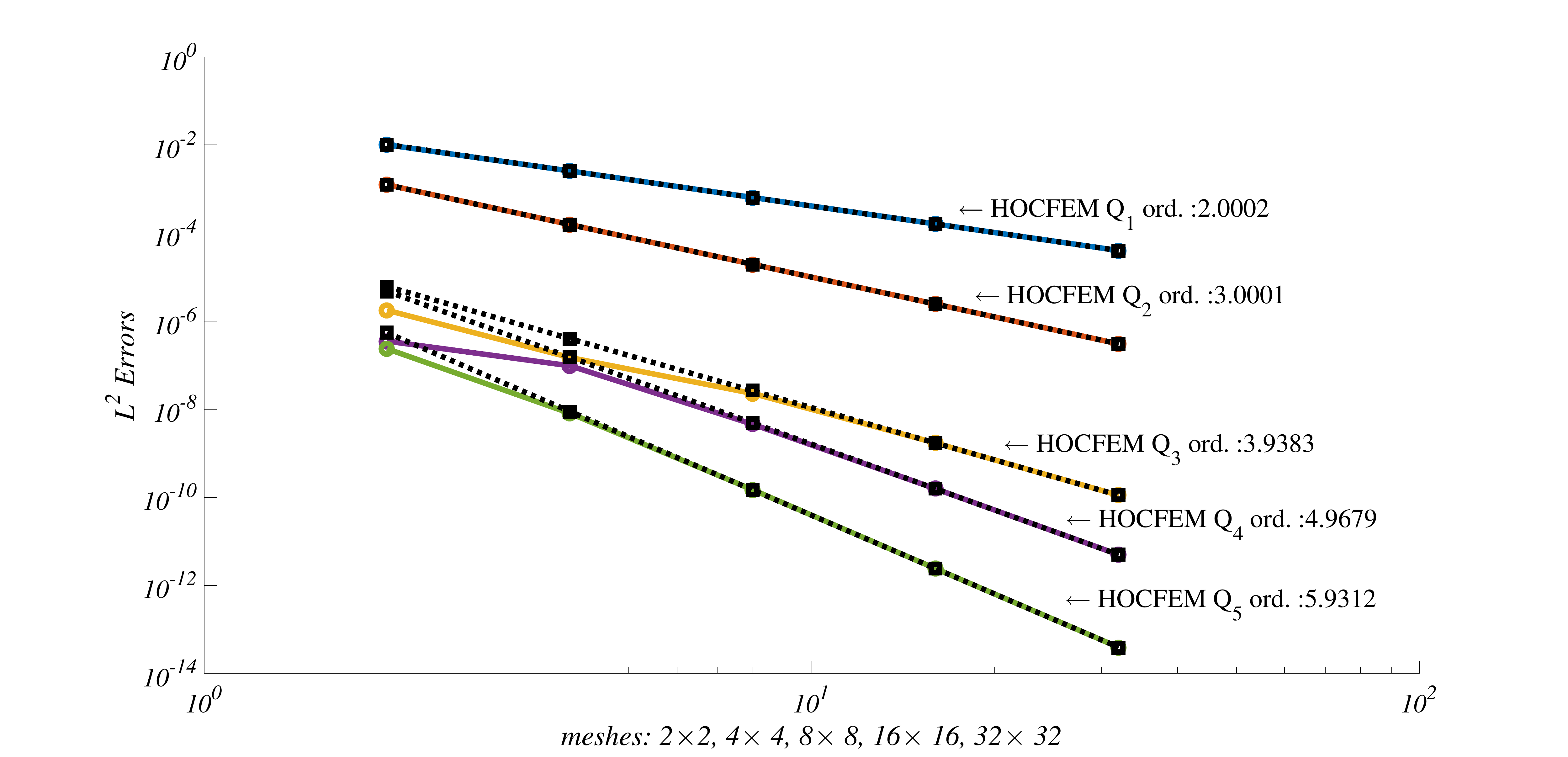}
\caption{Log-log graphic  of errors in $L^2$ 
norm for {HOCFEM} numerical solutions of Problem 
(\ref{pres}) with data (\ref{eq:Example:1}), using 
$\Q_1$, $\Q_2$, $\Q_3$, $\Q_4$ and $\Q_5$ basis 
through a mesh refinement.}
\label{logL22} 
\end{figure}

\begin{figure}[ht!]
\centering
\includegraphics[scale=.3, trim = 4cm 0 0 0]{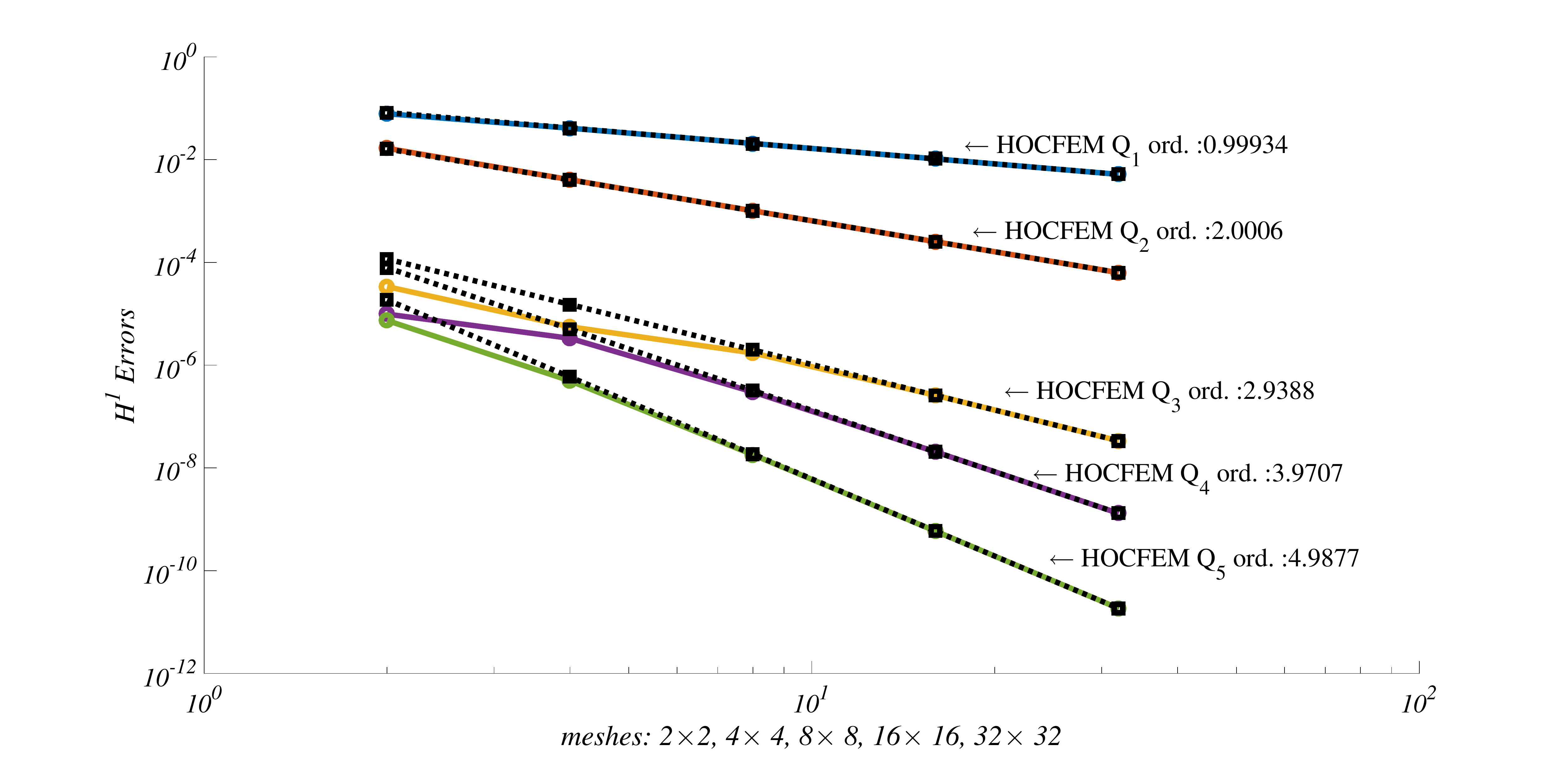} 
\caption{Log-log graphic  of errors in $H^1$ 
norm for {HOCFEM} numerical solutions of Problem 
(\ref{pres}) with data (\ref{eq:Example:1}), using 
$\Q_1$, $\Q_2$, $\Q_3$, $\Q_4$ and $\Q_5$ basis 
through a mesh refinement.}
\label{logH12}
\end{figure}

We apply HOCFEM to example (\ref{eq:Example:1}) 
with $\mathbb{Q}_1, \cdots ,\mathbb{Q}_5$ finite element 
spaces to the problem and compute HOCFEM and FEM solutions. 
We estimate the $L^2$ and $H^1$ errors and plot 
them in a log-log graphic shown in Figures 
\ref{logL22} and \ref{logH12}.


Numerical convergence is observed in $H^1$ norm 
with the same error variation  rate of FEM. The 
error $p-p^h$  in $L^2$ norm is not optimal but 
with the correction $p-(p^h+\lambda^h)$ (in each 
control volume) we recover convergence rate of 
FEM in this norm for HOCFEM. We compute conservation 
of energy indicator (\ref{eq:IndicatorEnergy}) and 
conservation of mass indicator (\ref{eq:IndicatorMass})
by the below formulas, respectively, 
\begin{equation}
E(p) = \frac{1}{2}\int_\Omega  
\Lambda\vert \nabla p\vert^2  - \int_\Omega qp, 
\label{eq:IndicatorEnergy}
\end{equation}
\begin{equation}
J(p) = \left(\sum_{R} \left(\int_{\partial R}  
\Lambda\nabla p\eta   - \int_{R} q \right)^2\right)^{1/2}.
\label{eq:IndicatorMass}
\end{equation}

\begin{table}[ht]
\centering
\renewcommand{\arraystretch}{1.1}
{\color{black}
\begin{tabular}{c|c|c||c|c} 
   $Q_r$  &  $ E(u_{FEM})$   &  $E(u_{HOCFEM})$ &  
   $ J(u_{FEM})$   &  $J(u_{HOCFEM})$    \\
    \hline  $1$  
    &  $-4.514912976$  &  $-4.514911724$ 
    &  $3.304137047\times 10^{-4}$  &  $5.893618438\times 10^{-15}$ \\ 
    \hline  $2$  
    &  $-4.523567134$  &  $-4.523565879$ 
    &  $3.308277779\times 10^{-4}$  &  $6.553391232\times 10^{-15}$ \\ 
     \hline  $3$  
    &  $-4.523568684$  &  $-4.523568684$ 
    &  $2.295180099\times 10^{-8}$  &  $1.902320206\times 10^{-14}$ \\ 
    \hline  $4$  
    &  $-4.523568684$  &  $-4.523568684$ 
    &  $2.295194635\times 10^{-8}$  &  $1.805166187\times 10^{-14}$ \\ 
    \hline  $5$  
    &  $-4.523568684$  &  $-4.523568684$ 
    &  $1.454599166\times 10^{-12}$  &  $3.207818336\times 10^{-14}$ \\ 
    \hline  $6$  
    &  $-4.523568684$  &  $-4.523568684$ 
    &  $3.489145981\times 10^{-12}$  &  $3.408881693\times 10^{-14}$ \\ 
\end{tabular}
}
\caption{{\color{black}Energy minimization and conservation 
indicator with $h = 2^{-9}$ for numerical solution 
of Problem (\ref{pres}) with data (\ref{pres}).}} 
\label{tab:conservationindicator1}
\end{table}

As we can see in Table \ref{tab:conservationindicator1} 
conservation  of Energy remains similar for both methods 
while HOCFEM exhibit much better conservation of mass than 
classical FEM. 

\subsubsection{Heterogeneous medium}

\begin{figure}[ht!]
\centering
\includegraphics[scale=0.85]{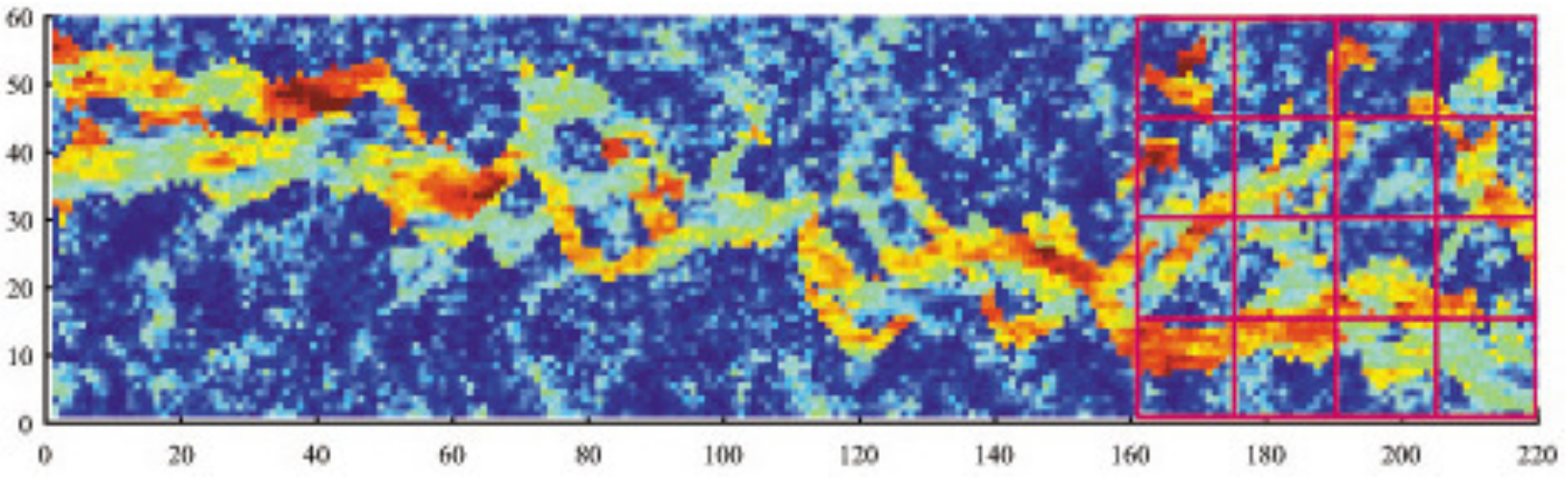} 
\caption{A 2D layer of the $SPE10$ 2D porous medium 
sample from \cite{TenhSPE}.}
\label{sp10}
\end{figure}

Let us move to a heterogeneous medium with high-contrast 
coefficients. The medium to be consider 
is the last $64 \times 64$ block of the 
geological $SPE10$ porous medium taken from 
\cite{TenhSPE} shown in Figure \ref{sp10}. This is 
a widely used heterogeneous porous medium for 
simulations (see for example \cite{TenhSPE}).

\begin{figure}[ht!]
\centering 
\includegraphics[scale=.3, trim = 4cm 0 0 0]{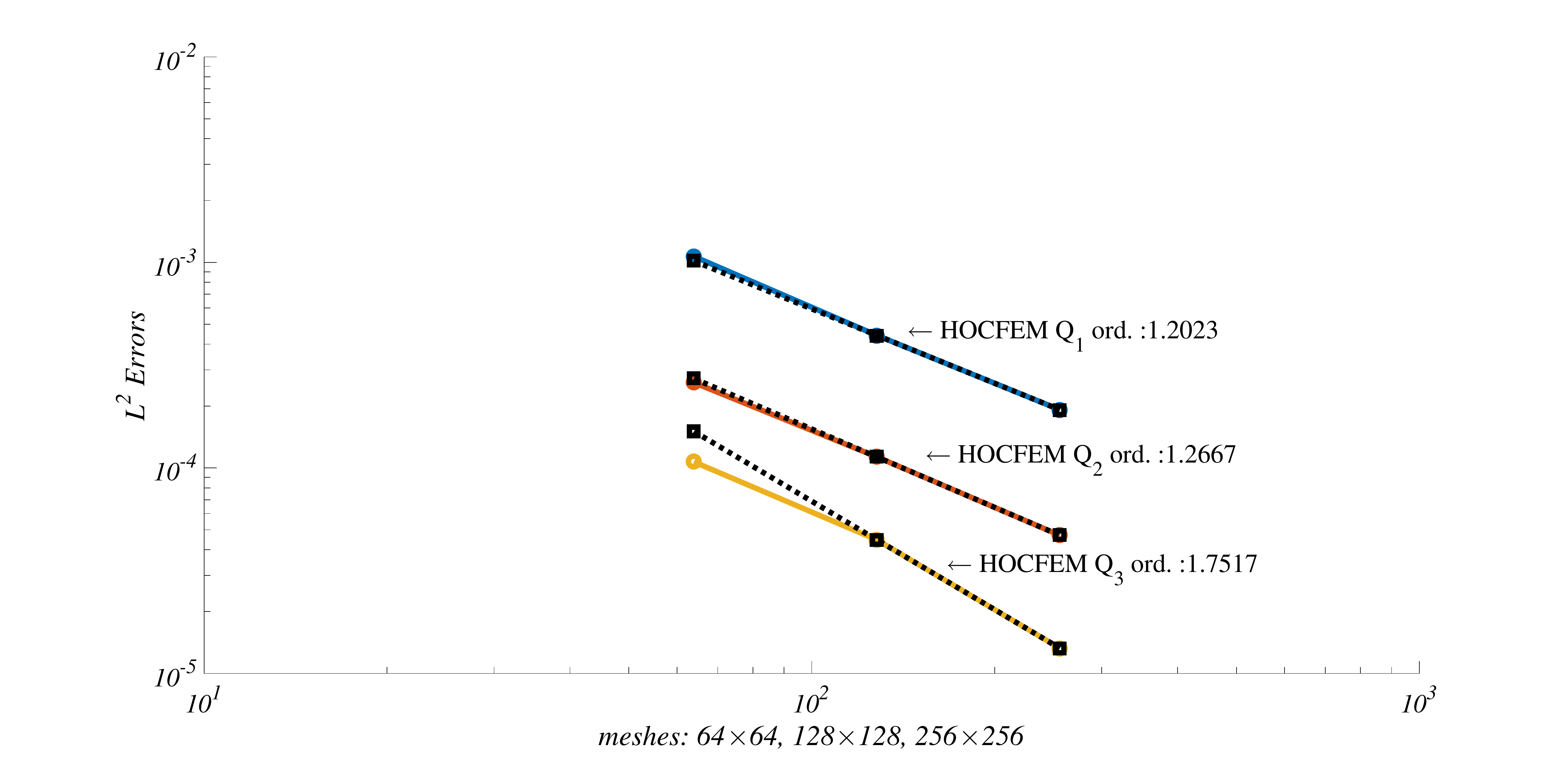}
\caption{Log-log graphic  of errors in $L^2$ norm 
for {HOCFEM} numerical solutions of Problem 
(\ref{pres}) in SPE10 medium, 
using $\Q_1$ and $\Q_2$ basis through a mesh refinement. 
The norm of the mesh  $h = 2^{-M}$ vary as $M=6,7,8$.}
\label{Fig:SPE10_L2} 
\end{figure}

\begin{figure}[ht!]
\centering
\includegraphics[scale=.3, trim = 4cm 0 0 0]{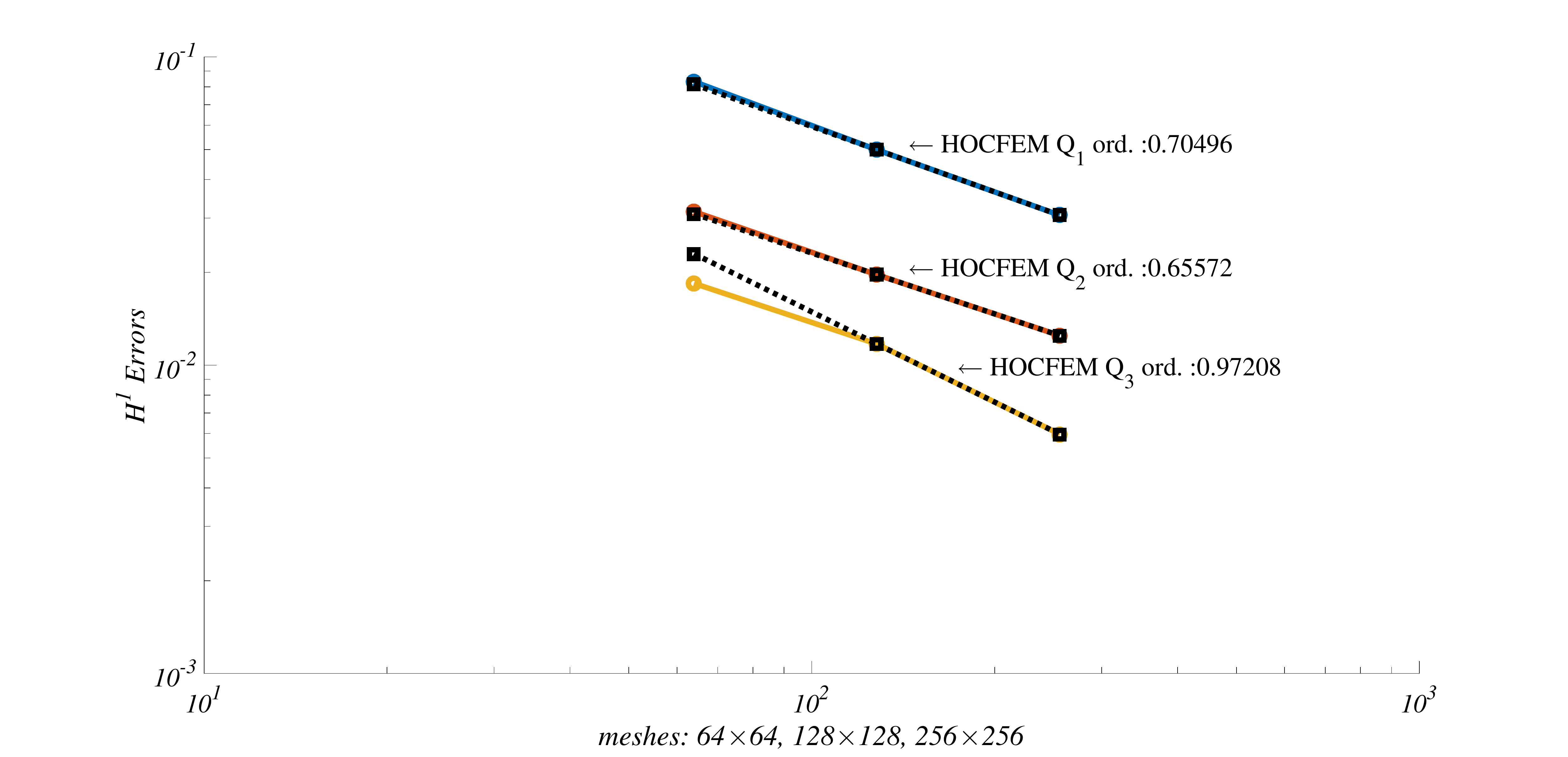} 
\caption{Log-log graphic  of errors in $H^1$ norm 
for {HOCFEM} numerical solutions of Problem 
(\ref{pres}) in SPE10 medium with data (\cite{TenhSPE}), 
using $\Q_1$ and $\Q_2$ basis through a mesh refinement. 
The norm of the mesh  $h = 2^{-M}$ vary as $M=6,7,8$.}
\label{Fig:SPE10_H1}
\end{figure}

We perform a numerical experiment to study  
convergence  and  conservation of  energy and mass  
of HOCFEM in such a realistic heterogeneous medium.
Let us consider the model problem (\ref{pres}) on 
$\Omega = \left[0,1\right]\times \left[0,1\right]$ 
with constant forcing 
term and homogeneous Dirichlet's boundary conditions 
over $\partial \Omega$. The mobility coefficient 
$\Lambda$ is taken from the SPE10 medium as 
described before.
We compute HOCFEM approximations using $\Q_1$ 
and $\Q_2$ basis  over $3$ square meshes of norm  
$h = 2^{-M}$ with $M = 6,7,8$ and compute errors 
in $L^2$ and $H^1$ norms against a reference solution 
calculated using $\Q_3$ basis in the finest mesh.  
Figures \ref{Fig:SPE10_L2} and \ref{Fig:SPE10_H1} 
show the $L^2$ and $H^1$ errors computed for both 
solutions using $\Q_1$ and $\Q_2$. {We observe 
that the  error variation rate is similar for both 
solutions. This behavior of the error in heterogeneous 
medium is different to homogeneous medium where 
rates where proportional to the grade of the 
polynomials used in basis.}

We also compute the approximated solution $u^h_{HOCFEM}$ 
solving system in (\ref{eq:saddle}) for a $256\times 256$ 
computational mesh and estimate the conservation of 
energy and mass indicators defined  in 
(\ref{eq:IndicatorEnergy}) and (\ref{eq:IndicatorMass}). 
Table \ref{tab:conservationindicator2} shows both,  
local mass and energy conservation indicators, for 
our high order HOCFEM formulation ($E(u_{HOCFEM})$ 
and $J(u_{HOCFEM})$) and for classical FEM ($E(u_{FEM})$ 
and $J(u_{FEM})$).


\begin{table}[ht!]
\centering
\renewcommand{\arraystretch}{1.1}
{\color{black}
\begin{tabular}{c|c|c||c|c} 
   $Q_r$  &  $ E(u_{FEM})$   &  $E(u_{HOCFEM})$ &  
   $ J(u_{FEM})$   &  $J(u_{HOCFEM})$    \\
    \hline  $1$  
    &  $-2.610659111$  &  $-2.606353663$ 
    &  $6.478664199\times 10^{-1}$  &  $2.734482311\times 10^{-12}$ \\ 
    \hline  $2$  
    &  $-2.755620550$  &  $-2.752274173$ 
    &  $5.472598738\times 10^{-1}$  &  $2.816845987\times 10^{-12}$ \\ 
     \hline  $3$  
    &  $-2.778407981$  &  $-2.778404839$ 
    &  $4.095003567\times 10^{-2}$  &  $9.330127650\times 10^{-12}$ \\ 
    \hline  $4$  
    &  $-2.786675138$  &  $-2.786672433$ 
    &  $3.959310237\times 10^{-2}$  &  $7.837454755\times 10^{-12}$ \\ 
    \hline  $5$  
    &  $-2.790778373$  &  $-2.790778322$ &  
    $8.962957140\times 10^{-3}$  &  $1.509701588\times 10^{-11}$ \\ 
    \hline  $6$  
    &  $-2.793138428$  &  $-2.793138381$ 
    &  $8.972699240\times 10^{-3}$  &  $1.289345855\times 10^{-11}$ \\ 
\end{tabular}
}
\caption{{\color{black} Energy minimization and conservation indicator  
for numerical solution of Problem (\ref{pres}) with 
data (\ref{eq:Example:1})  in a fixed mesh $64 \times 64$ 
using basis $\Q_1$, $\Q_2$, $\Q_3$, $\Q_4$, $\Q_5$ and $\Q_6$.}} 
\label{tab:conservationindicator2}
\end{table}

From Table \ref{tab:conservationindicator2} we see 
that conservation of global energy does not change 
from FEM to HOCFEM while conservation of mass is
superior in large with our new formulation.

\section{Conservation properties of the no flow surface 
region for hyperbolic conservation laws}\label{theLE2d}

The aim of this section is to present 
{\color{black} an} extension of the Lagragian-Eulerian scheme (see 
\cite{EAJP19,APS17a,APS17b,APS18,AMPRB19,AAP19}) 
for hyperbolic conservation laws in two-space dimensions 
with some initial condition coming from abstract 
nonlinear problems of hyperbolic conservation laws.
We can also consider problems of physical interest in 
fluid mechanics such  as multiscale flow in porous 
media scalar and systems treated in this work.

{\color{black}
We mention that our novel 
hyperbolic Lagrangian-Eulerian solver (it its simplest 
form) can be viewed as a monotone scheme (see Section 
\ref{monoLE}).
For coupling Darcy flow, the transport method captures 
fine-scale effects using a (conservative) fine-grid finite element 
technique combined with Lagragian-Eulerian scheme.
For the purpose of this work we use Cartesian grids
since for the case of {\it monotone scheme} convergence and error analysis 
reduces (essentially) to a one-dimensional problem and 
it retains  convergence results and approximation  
to the entropy weak solution by recalling 
\cite{EGH95a,EGH95b,CCH99,CRMAJ80}.


We improve the interpretation of the construction 
of numerically stable Lagrangian-Eulerian {\it no 
flow surface region} in two-space dimensions previously 
presented and analyzed in \cite{EAJP19} for 
one-dimensional balance and conservation laws. It turn out
that our {\it monotone Lagrangian-Eulerian} is a building 
block for construction of a novel class of Lagrangian-Eulerian 
shock-capturing schemes for first-order hyperbolic problems.} 
{\color{black}
The early monotone versions of the Lagrangian-Eulerian approach 
has been employed successfully in a number of very non-trivial 
problems and also developed theoretically 
\cite{EAJP19,APS18,APS17a,APS17b,AMPRB19,AAP19} 
linked to several transport models such as the Burgers' 
equation with Greenberg-LeRoux's and Riccati's 
source terms, the shallow-water system, Broadwell's rarefied 
gas dynamics, Baer-Nunziato's system linear, non-linear 
convex and non-linear non-convex 2D scalar conservation laws 
(see \cite{APS17a,EAJP19}). 
It is worth mentioning that the Lagrangian-Eulerian framework 
is able to compute qualitatively correct (entropy) solutions 
involving intricate non-linear wave interactions of rarefaction 
and shock waves. It is of significance to mention that the scheme 
is able to handle resonance effect associated to non classical 
transitional shock in a $2 \times 2$ three-phase flow 
water-oil-gas system \cite{APS17a,APS17b} and an intricate 
shock structure linked to a $5 \times 5$ isentropic Baer-Nunziato 
model (see \cite{EAJP19}). The Lagrangian-Eulerian scheme does 
handle { properly the sonic} rarefaction linked to Burgers' 
equation, namely, a typical small (and unphysical) discontinuity 
jump within the rarefaction structure; such discontinuation in 
the solution is unphysical, and thus with no mathematical 
relation with an entropy violating shock. Indded, our 
Lagrangian-Eulerian scheme does not produce the well-known 
spurious entropy glitch effect in the sonic rarefaction (as 
is the case of Rusanov and Godunov monotone schemes). In
addition, our first-order monotone Lagrangian-Eulerian scheme 
is less difusive than the classical Lax-Friedrichs scheme, 
but retains robustness and it is simple to implement and 
efficient for numerical computing \cite{AELP20}.

A key hallmark of the our Lagrangian-Eulerian (monotone) method is the 
dynamic tracking forward of the no-flow region (per time step). 
This is a considerable improvement compared to the classical 
backward tracking over time of the characteristic curves over 
each time step interval, which is based on the strong form of 
the problem. Indeed, in the case of systems and multi-D 
problems, we can say that backward tracking is not understood. 

Our new method can handle, with great simplicity, nontrivial 
scalar and systems problems in 1D and multi-D \cite{EAJP19,APS17a}.
Another key hallmark of the our Lagrangian-Eulerian (monotone) 
method is the {\it a flux separation strategy} and its impact 
on the balancing (multiple scale) discretization between the first-order 
approximation of the hyperbolic flux and the source term to take into 
account nonlinear wave interactions preserving conservation properties. 
For instance, in \cite{EAJP19} the numerical tests show that the 
discretizations resulting from the {flux separation strategy} 
when applied to the 2 by 2 shallow-water system and 5 by 
Baer-Nunziato's system seem to be of good quality. Moreover, 
such strategy seem to be very appropriate to deal with convex 
and nonlinear non-convex 2D scalar conservation laws.

With respect to the theory of monotone scheme, the {\it no flow 
region} (see Figure \ref{fII}) is the control volume where the 
(local) wave interaction (always in the fine mesh of any multiscale 
method approach) takes place. On the other hand, in light of 
modern reasearch (see \cite{ECCL15, BLIS15} and \cite{EAJP19}), the 
{\it no flow region} is a space-time cutoff to account the 
complex and intricate nonlinear wave group interaction 
within control volume {\it per time step} in the overall
simulation (\cite{CAFM16,ASFM10,DMBP01,EAJP19,APS17a}). In 
computing practice, the {\it no flow region} parallel with 
the CFL stability criterion associated with the space-time 
discretization of many numerical methods. 

Therefore, the monotone Lagrangian-Eulerian approach is a 
interesting novel framework for hyperbolic conservation 
laws and multiscale transport flow models.
}

\subsection{Lagragian-Eulerian technique with conservation properties}
We discus our new Lagragian-Eulerian technique with 
conservation properties for the approximation of the 2D 
initial value problem for hyperbolic of conservation laws,
\begin{equation}
 \begin{cases}
  \displaystyle\frac{\partial u}{\partial t} + \frac{\partial f(u)}{\partial x}
   + \frac{\partial g(u)}{\partial y} = 0, \qquad (x,y,t) \in \Omega \times (t^0,T], \\
   u(x,y,t^0) = \eta(x,y), \quad (x,y) \in \Omega, 
  \end{cases}
 \label{2D1H}
\end{equation}
where $\Omega$ is a interior square domain in $\mathds{R}^2$, 
whit boundary $\partial \Omega$ and $T=t_f >0$. \\

For the finite dimensional function spaces we introduce 
the following standard notation. The space region
$\left(\mathds{R} \times \mathds{R}\right) \times \overline{\mathds{R}} 
= \{(x,y,t) : \,\, - \infty <  x,y < \infty, \,\, t > 0 \}$ 
is replaced by the lattice 
$\left(\mathds{Z} \times \mathds{Z}\right) \times \mathds{N} 
= \{(i,j,n) : \,\, i,j = 0, \pm 1, \pm 2,\dots; \,\, n = 0,1,2,\cdots\ \}$.  We
consider the sequence $U^n = (U^n)_{i,j}$, $i,j \in \mathds{Z}$ for $n = 0,1,2,...$, 
for a given grid size $\Delta x, \Delta y > 0$ and time level
\begin{equation}
 t^n = \sum_{i = 0}^n \Delta t^i, \quad \text{with $t^0 = 0$},
\quad \text{for non-constant time steps $\Delta t^i$'s}, 
 \label{time2D}
\end{equation}
subject to 
the CFL constraint  (which determines the maximum 
allowable time-step to garantee the desire 
conservation properties of the no flow surface region).\\

In the time level $t^n$, we have, 
$$(x_i^n,y_j^n) = {\color{black}(i\Delta x,j\Delta y)}\quad \mbox{ and } \quad \left(x_{i + \frac{1}{2}}^n,x_{i + \frac{1}{2}}^n\right) = 
\left(i \Delta x  + \frac{\Delta x}{2},j \Delta y  + \frac{\Delta y}{2}\right)$$ on the 
uniform local grid (or primal grid). Here
$$h_{x,i}^n = \Delta x^n = x_{i + \frac{1}{2},j}^n - 
x_{i - \frac{1}{2},j}^n, 
\quad \mbox{ and } \quad  h_{y,j}^n = \Delta y^n = y_{i,j + \frac{1}{2}}^n - 
y_{i,j - \frac{1}{2}}^n, \quad \mbox{ for } i,j \in \mathds{Z},$$  where 
$(x_{i \pm \frac{1}{2}}^n,y_{i \pm \frac{1}{2}}^n)$ 
are the corners of the $(i,j)$-cell. 
For the non-uniform grid we have 
$h_{x,i}^{n+1} = \overline{\Delta x}^{n+{1}} = 
\overline{x}_{i + \frac{1}{2}}^{n+{1}} - 
\overline{x}_{i - \frac{1}{2}}^{n+{1}}$
and 
$h_{y,j}^{n+{1}} = 
\overline{\Delta y}^{n+{1}} = 
\overline{y}_{j + \frac{1}{2}}^{n+{1}} - 
\overline{y}_{j - \frac{1}{2}}^{n+{1}}$, in the time 
level $t^{n+1}$. 
\begin{figure}[ht!]
\centering
\includegraphics[width=0.4\textwidth]{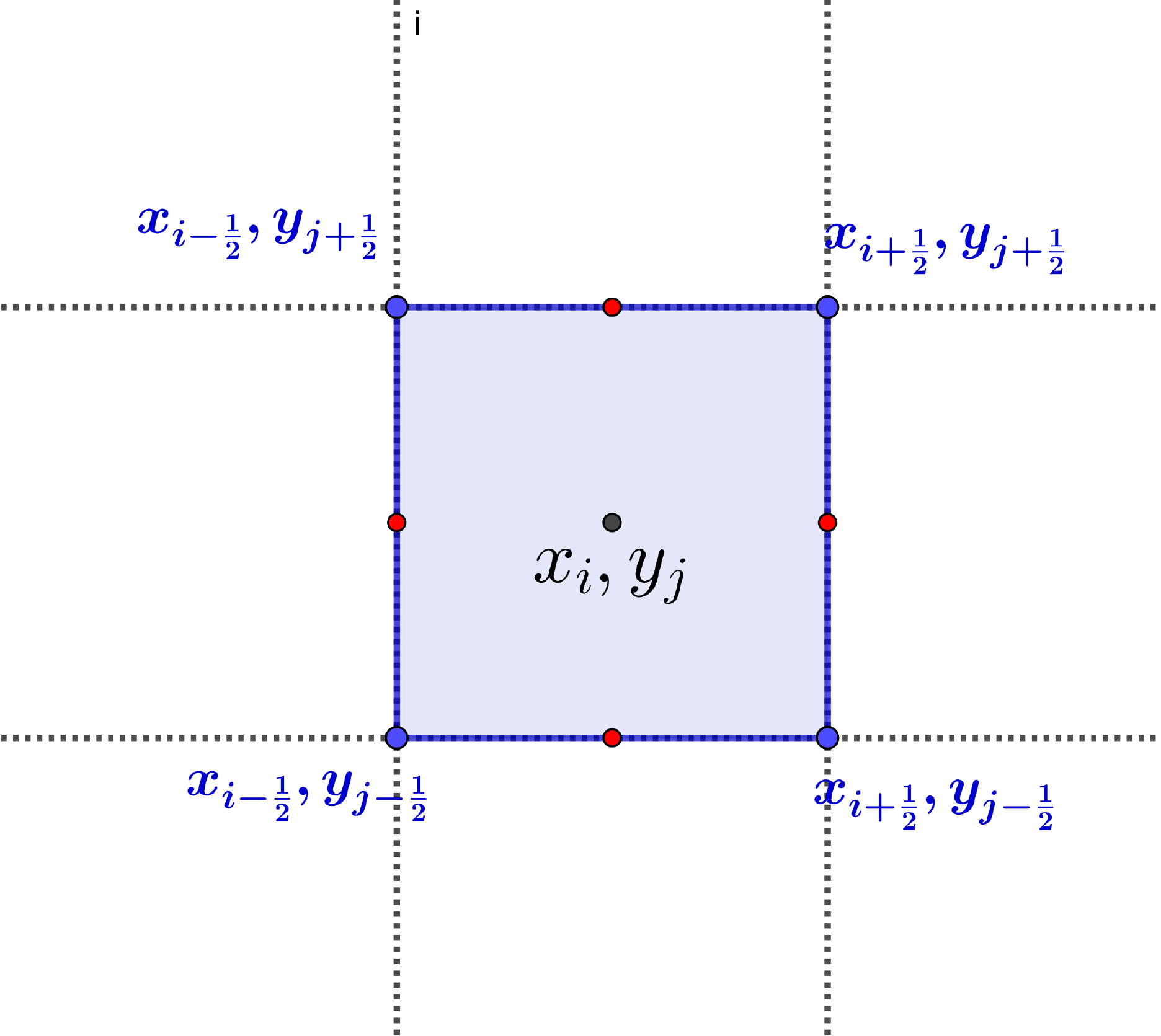}
\caption{Illustration of the notation related to the ,$(i,j)$-cell.} 
\label{fI}
\end{figure}

\begin{figure}[ht!]
	\centering
  \includegraphics[width=0.37\textwidth]{3dint-eps-converted-to.pdf}
	\includegraphics[width=0.37\textwidth]{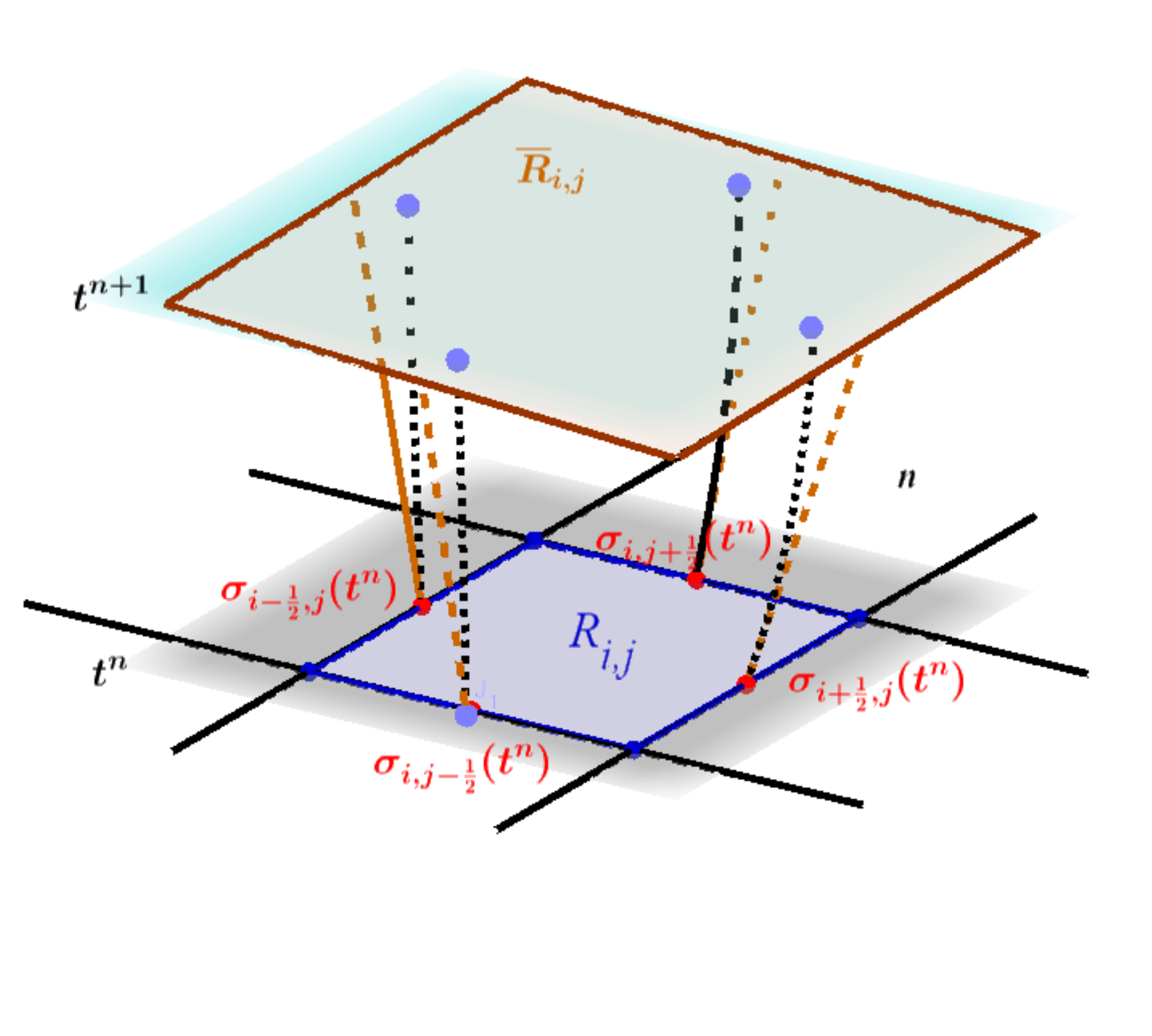}
	\caption{Analitic (left) and approximated (right) no flow surface region.} 
	\label{fIV}
\end{figure}

\begin{figure}[ht!]
	\centering
  \includegraphics[width=0.4\textwidth]{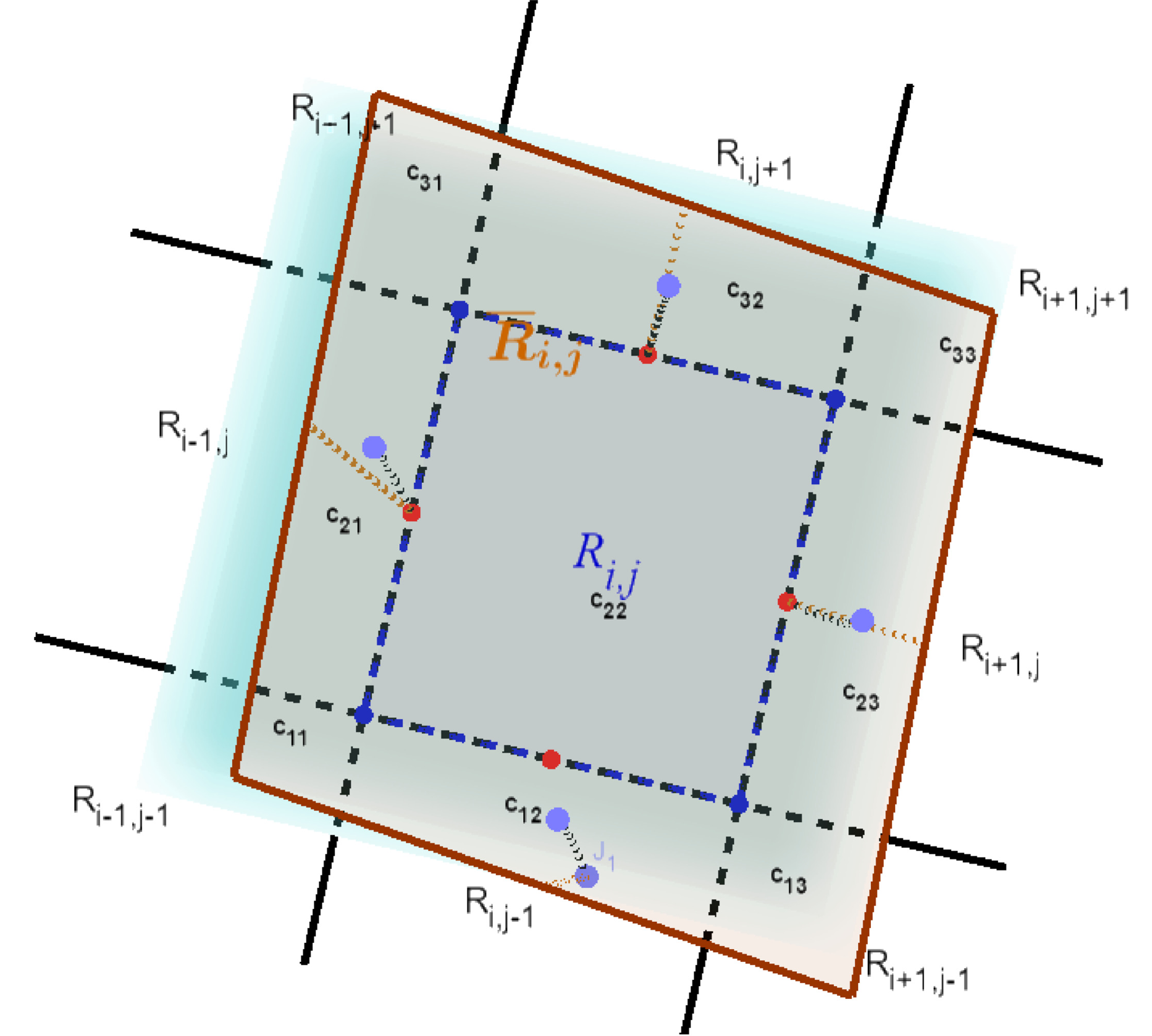}
	\caption{Projection over original grid.} 
	\label{fV}
\end{figure}

The pair $(x_i^n,y_j^n)$ is the centers of the $(i,j)$-cell, 
$i,j \in \mathds{Z}$. From now on, for short,  when there 
is no chance of misunderstanding, the limits of integration 
will indicate the time level where integration calculation 
takes place with respect to the pair $(x_i,y_j)$ of the 
$(i,j)$-cell, $i,j \in \mathds{Z}$. In each cell
$\left[x_{i - \frac{1}{2}}^n,x_{i + \frac{1}{2}}^n \right] \times 
\left[y_{j - \frac{1}{2}}^n,y_{j + \frac{1}{2}}^n\right]$
(see Figure \ref{fI}), the approximate solution for (\ref{2D1H}),  
is defined by
\vspace{-2mm}
\begin{equation}
 U(x_i,y_j,t^n) = U_{i,j}^n \equiv \frac{1}{\Delta x \Delta y} 
                \int_{x_{i - \frac{1}{2}}^n}^{x_{i 
               + \frac{1}{2}}^n} \int_{y_{j - \frac{1}{2}}^n}^{y_{j 
               + \frac{1}{2}}^n} u(x,y,t^n)\,dx\,dy, 
 \label{2DAp1}
\end{equation}
\vspace{-5mm}
\begin{equation}
 \overline{U}(x_i,y_j,t^{n+1}) = \overline{U}_{i,j}^{n+1} 
\equiv \frac{1}{h_{x,i}^{n+1} h_{y,j}^{n+1}} 
 \int_{\overline{x}_{i - \frac{1}{2}}^{n+1}}^{\overline{x}_{i 
+ \frac{1}{2}}^{n+1}} 
 \int_{\overline{y}_{j - \frac{1}{2}}^{n+1}}^{\overline{y}_{j 
+ \frac{1}{2}}^{n+1}} u(x,y,t^{n+1})\,dx\,dy, 
 \label{2DAp2}
\end{equation}
along with the initial condition
$U(x_i^0,y_j^0,t^0) = U_{i,j}^0$ in the cells 
$\left[x_{i - \frac{1}{2}}^0,x_{i + \frac{1}{2}}^0\right] \times 
\left[y_{j - \frac{1}{2}}^0,y_{j + \frac{1}{2}}^0\right]$,
$i,j \in \mathds{Z}$. 
It is worthy to mention that the approximation value 
$\overline{U}(x_i,y_j,t^{n+1})$ is performed over the region 
$\overline{R}_{i,j}^{n+1}$; see the right picture in Figure 
\ref{fIV} as well as Figure \ref{fV}, for an illustration of 
the projection procedure over original grid in control volumes.
Note that in  (\ref{2DAp1}) and (\ref{2DAp2}), the quantity 
$u(x,y,t)$ is a solution of (\ref{2D1H}). The discrete 
counterpart of the space $L^p(\mathds{R}^2)$ is 
$l_{\Delta x, \Delta y}^p$, the space of sequences 
$U = (U_{i,j})$, with $i,j \in \mathds{Z}$, with norm given by 
$$\|U\|_{l_{\Delta x ,\Delta y}^p} = \left(\Delta x \Delta y\sum_{i \in
\mathds{Z}}\sum_{j \in \mathds{Z}} |U_{i,j}|^p \right)^{\frac{1}{p}},
\quad \text{ where } 1 \leq p < \infty. $$ 

To build the new two dimensional scheme we extend the concept 
of {\it no flow surface region} {\color{black}$D_{ij}^n$} (see \cite{EAJP19,APS17a}) 
to three dimensional variables ($x$,$y$ and $t$) as 
$D_{i,j}^n \subset \mathds{R}^3$, where $i$ and $j$ refer to
$(x_i,y_j)$ and $n$ refers to time state $t^n$. The border 
of the control volume $D_{i,j}^n$ is 
represented by $\partial D_{i,j}^n = R_{i,j}^n \cup S_{i,j}^n
\cup \overline{R}_{i.j}^{n+1}$ where (see Figure \ref{fI}),
\begin{itemize}
\item $R_{i,j}^n = \left[x_{i-\frac{1}{2}}^n, x_{i+\frac{1}{2}}^n\right] \times 
\left[y_{j-\frac{1}{2}}^n, y_{i+\frac{1}{2}}^n\right]$ in 
$\mathds{R}^2$ is the entry of the no flow surface region
\item $\overline{R}_{i,j}^{n+1} = 
\left[\overline{x}_{i-\frac{1}{2}}^{n+1}, 
\overline{x}_{i+\frac{1}{2}}^{n+1} \right] \times 
\left[\overline{y}_{j-\frac{1}{2}}^{n+1}, 
\overline{y}_{i+\frac{1}{2}}^{n+1} \right]$ is the exit of 
the no flow sufarce region, and
\item  $S_{i,j}^n$, in $\mathds{R}^3$, is the 
lateral surface of the no flow surface region. 
\end{itemize}
We consider now (\ref{2D1H}) in the generalized space-time 
divergence form,
\begin{equation}
   \nabla_{t,x,y} \left[\begin{array}{c}
                  u \\f(u) \\ g(u)
                    \end{array}\right] = 0.
 \label{LG2DExpl1}
\end{equation}
Integration over the {\it control volume} and the use 
of the divergence theorem gives,
\begin{equation}
                   \int_{{\color{black} \partial D_{i,j}^n}} \left[\begin{array}{c}
                  u \\ f(u) \\ g(u)
                    \end{array}\right] \cdot 
\vec{\textsl{n}}\,\,d(\partial D_{i,j}^n) = 0. 
  \label{LG2DExpl2}
\end{equation}
The normal vector in the entry of the no flow surface region,  
$R_{i,j}^n$, is $\left[-1 \,\, 0 \,\, 0\right]^T$ and the 
vector normal in the exit of the no flow surface region, 
$\overline{R}_{i,j}^{n+1}$,  is $\left[1 \,\, 0 \,\, 0\right]^T$. 
Then, the right 
side of (\ref{LG2DExpl2}) can be written as
\begin{equation}
 \int_{R_{i,j}^n} \left[\begin{array}{c}
                  u \\ f(u) \\ g(u)
                    \end{array}\right] \cdot \left[-1 \,\, 0 \,\, 0\right]^ T \,\, dA
                    + \int_{S_{i,j}^n} \left[\begin{array}{c}
                  u \\ f(u) \\ g(u)
                    \end{array}\right] \cdot \vec{\textsl{n}}\,\,dS
                    + \int_{\overline{R}_{i,j}^{n+1}} \left[\begin{array}{c}
                  u \\ f(u) \\ g(u)
                    \end{array}\right] \cdot \left[1 \,\, 0 \,\, 0\right]^T \,\,dA = 0. 
  \label{LG2DExpl3}
\end{equation}
We assume there is not flow through the surface $S_{i,j}^n$ (that is, $S_{i,j}^n$ 
is impervious; this is natural in many applications 
\cite{EAJP19,APS17a,APS17b,APS18,AAP19}). Therefore 
{\it surface integral} of $S_{i,j}^n$ is zero, i.e.,
\begin{equation}
  \int_{\overline{R}_{i,j}^{n+1}} u(x,y,t^{n+1}) \,\, 
dA = \int_{R_{i,j}^n} \,\, u(x,y,t^n) \,\, dA,
  \label{LG2DExpl5}
\end{equation}
which we call {\it conservation identity}.
The numerical approximations $U_{i,j}^n$ and $\overline{U}_{i,j}^{n+1}$ 
appearing in (\ref{2DAp1}) and (\ref{2DAp2}), respectively, can be 
defined from equation (\ref{LG2DExpl5}) with the desired 
conservation properties and reads,
\begin{equation}
  \overline{U}_{i,j}^{n+1} = \frac{1}{\overline{R}_{i,j}^{n+1}}
\int_{\overline{R}_{i,j}^{n+1}} u(x,y,t^{n+1}) \,\, dA = 
	\frac{R_{i,j}^{n}}{\overline{R}_{i,j}^{n+1}}   
\frac{1}{R_{i,j}^{n}}\int_{R_{i,j}^n} \,\, u(x,y,t^n) \,\, dA = 
	{\color{black}\frac{R_{i,j}^{n}}{\overline{R}_{i,j}^{n+1}}} U_{i,j}^n.
  \label{ReltnTnp1}
\end{equation}
\begin{figure}[ht!]
	\centering
	\includegraphics[width=0.4\textwidth]{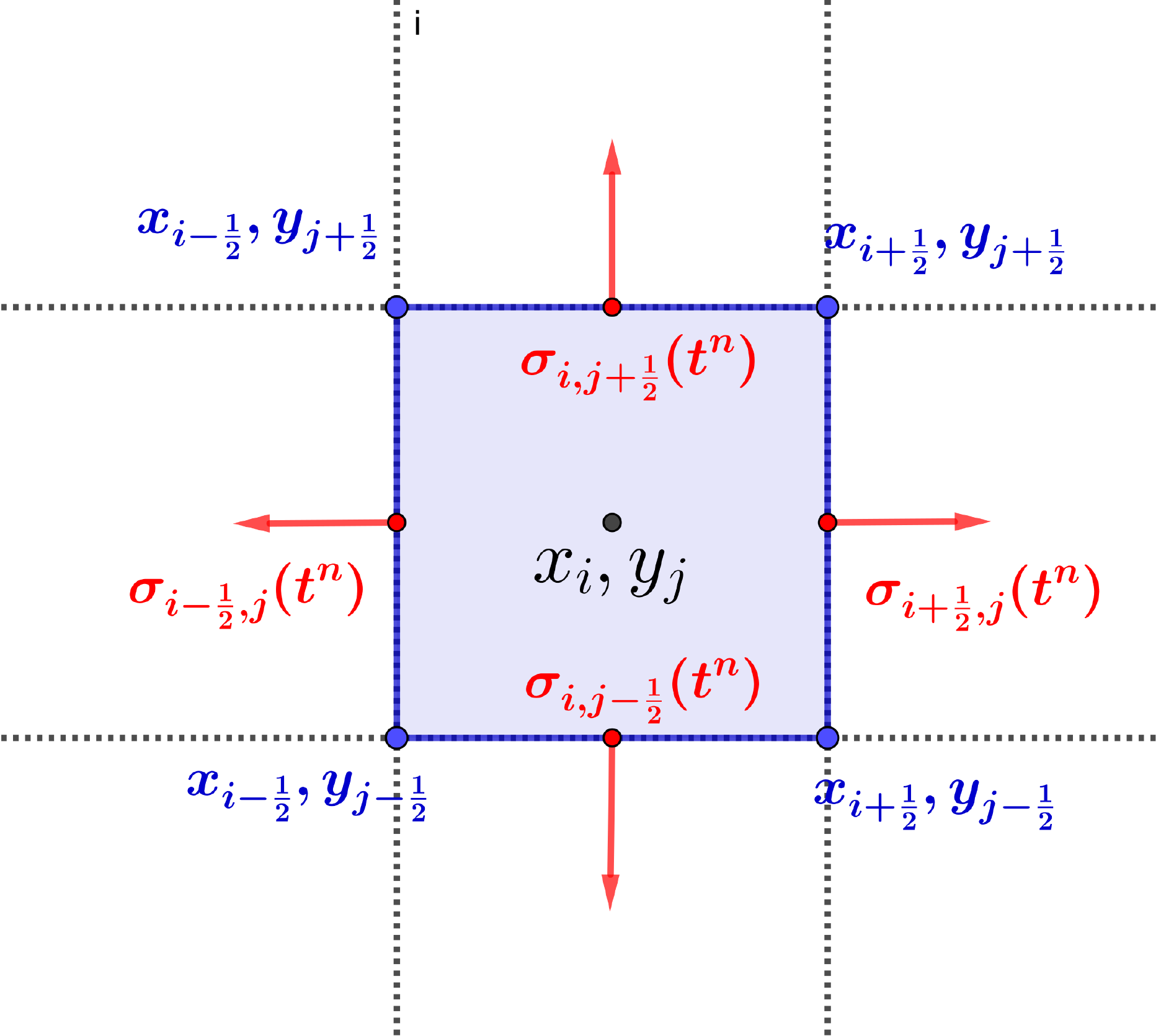}
	\includegraphics[width=0.4\textwidth]{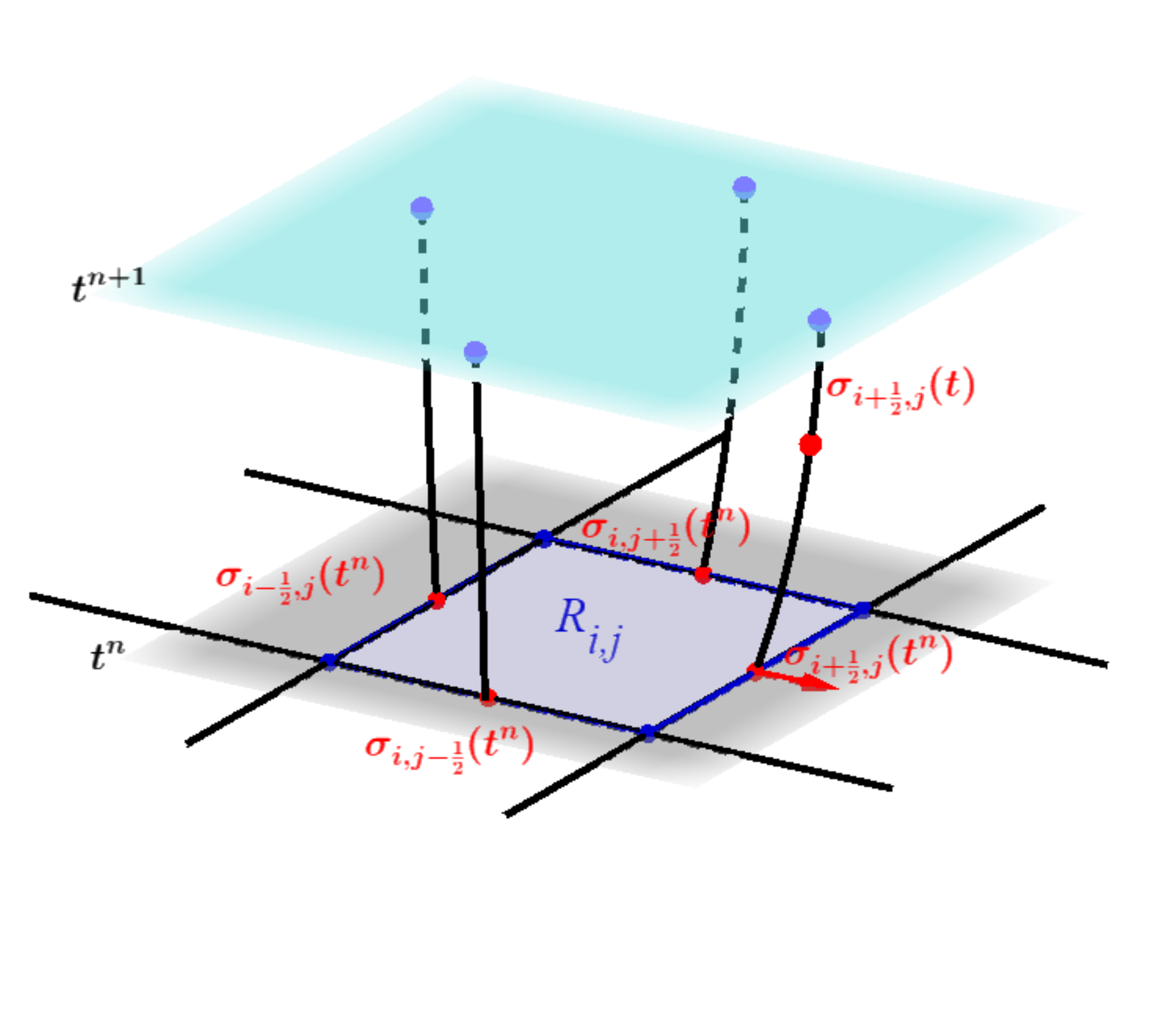}
	\caption{Normal vectors on each side of $R_{i,j}^{n}$-cell 
(left). No flow surface region and parametrized curves from
        time level $t^n$ to $t^{n+1}$ (right).} 
	\label{fII}
\end{figure}

On other hand, from (\ref{LG2DExpl1}) and by the natural conservation 
properties of the no flow surface region for hyperbolic conservation 
law it follows
\begin{equation}
        \int_{S_{i,j}^n}
				\left[\begin{array}{c}
           u \\ f(u) \\ g(u)
        \end{array}\right] \cdot \vec{\textsl{n}}\,\,dS = 0
\qquad \text{ and } \qquad
\left[\begin{array}{c}
           u \\ f(u) \\ g(u)
        \end{array}\right] \cdot \vec{\textsl{n}} = 0. 
\end{equation}
Now, let ${\color{black}\sigma_{i-\frac{1}{2},j}^{n}(t)} = (t,\sigma_1(t),\sigma_2(t))$, 
$t^n \leq t < t^{n}$,  be a parametrized curve such that 
$${\color{black}\sigma_{i-\frac{1}{2},j}^{n}(t^n)} = (t^n,\sigma_1(t^n),\sigma_2(t^n)) 
= (t^n,x_{i-\frac{1}{2}},y_{j})$$ and ${\color{black}\sigma_{i-\frac{1}{2},j}^{n}(t)}
\in S_{i,j}$.  Analogously, we can define parameterized 
curves correspoding to other sides of $R_{i,j}^{n}$ such that 
${\color{black}\sigma_{i+\frac{1}{2},j}^{n}(t), \sigma_{i,j-\frac{1}{2}}^{n}(t), 
\sigma_{i,j+\frac{1}{2}}^{n}(t)}  \in S_{i,j}$ with initial 
point ($t = t^n$) in the respective center of the side 
of $R_{i,j}^{n}$ (see Figure \ref{fII}).
				
Construction of a lateral curve of the {\it no flow surface 
region}. First, such construction is not unique. Actually, 
this might lead to a family of methods; this interesting 
issue will not be addressed in this work. Make fixed the 
point ${\color{black}\sigma_{i-\frac{1}{2},j}^{n}(t^n)}$. The normal vector 
to the corresponding side of $R_{i,j}^{n}$ is the vector 
$\left[t^n,-1,0\right]^T$. Moreover, the normal vector on 
the curve ${\color{black}\sigma_{i-\frac{1}{2},j}^{n}(t)}$ at $t^n \leq t < t_{n 1}$ 
is a orthogonal vector to vector 
${\color{black}(\sigma_{i-\frac{1}{2},j}^{n})'(t)} = \left[1,\sigma_1'(t),0\right]$; 
see right frame in Figure \ref{fII}. Indeed, the vector 
at point ${\color{black}\sigma_{i-\frac{1}{2},j}^n(t)}$ may be calculated 
as $n = \left[-1,\frac{1}{\sigma_{1}'(t)},0\right]$ and 
follows:
\begin{equation}
0 = \left[u \,\,\, f(u) \,\,\, g(u)\right]^T 
  \cdot \vec{\textsl{n}} = 
    \left[u \,\,\, f(u) \,\,\, g(u)\right]^T 
  \cdot \left[-1, \frac{1}{\sigma_{1}'(t)},0\right]= 
    -u + \frac{f(u)}{\sigma_{1}'(t)} {\color{black} \mbox{ and then }}
\label{odecp1}
\end{equation}
\begin{equation}
\sigma_{1}'(t) = \frac{f(u)}{u}
\label{odecp2}
\end{equation}
\noindent with 
${\color{black}\sigma_{i-\frac{1}{2},j}^n(t^n)} = (t^n,x_{i-\frac{1}{2}},y_{j})$.
Finally, since ${\color{black}\sigma_{i-\frac{1}{2},j}^n(t)}$ is in the plane 
$y = y_j$, then ${\color{black}\sigma_{i-\frac{1}{2},j}^n(t)} = (t,\sigma_1(t),y_{j})$. 
We point out that an analogous reasoning as in 
(\ref{odecp1})-(\ref{odecp2}) might lead to the parametrized curves 
${\color{black}\sigma_{i+\frac{1}{2},j}^n(t)} = \left[t,\gamma_1(t),y_j\right]$, 
${\color{black}\sigma_{i,j-\frac{1}{2}}^n(t)}=\left[t,x_i,\alpha_2(t)\right]$, 
${\color{black}\sigma_{i,j+\frac{1}{2}}^n(t)} = \left[t,x_i,\theta_2(t)\right]$, 
$\gamma_1(t)$, such that $\alpha_2(t)$ and $\theta_2(t)$ must 
satisfy {\it the exact} conditions,
\begin{equation}
\begin{array}{l l l}
  \begin{cases}
   \gamma_1(t)'(t) \!=\! \frac{f(u)}{u}, \\
   \gamma_1(t)'(t^n) \!=\! (t^n,x_{i+\frac{1}{2}},y_j),
  \end{cases}  
      \!\!\!\!\!\! & \!\!\!\!\!\!
  \begin{cases}
   \alpha_2(t)'(t) \!=\! \frac{g(u)}{u}, \\
   \alpha_2(t)'(t^n) \!=\! (t^n,x_i,y_{j-\frac{1}{2}}),
  \end{cases}
      \!\!\!\!\!\! & \!\!\!\!\!\!
  \begin{cases}
   \theta_2(t)'(t) \!=\! \frac{g(u)}{u}, \\
   \theta_2(t)'(t^n) \!=\! (t^n,x_i,y_{j+\frac{1}{2}}).
  \end{cases}
\end{array}
\label{syst}
\end{equation}

\begin{figure}[ht!]
	\centering
	\includegraphics[width=0.4\textwidth]{3dint-eps-converted-to.pdf}
	\includegraphics[width=0.4\textwidth]{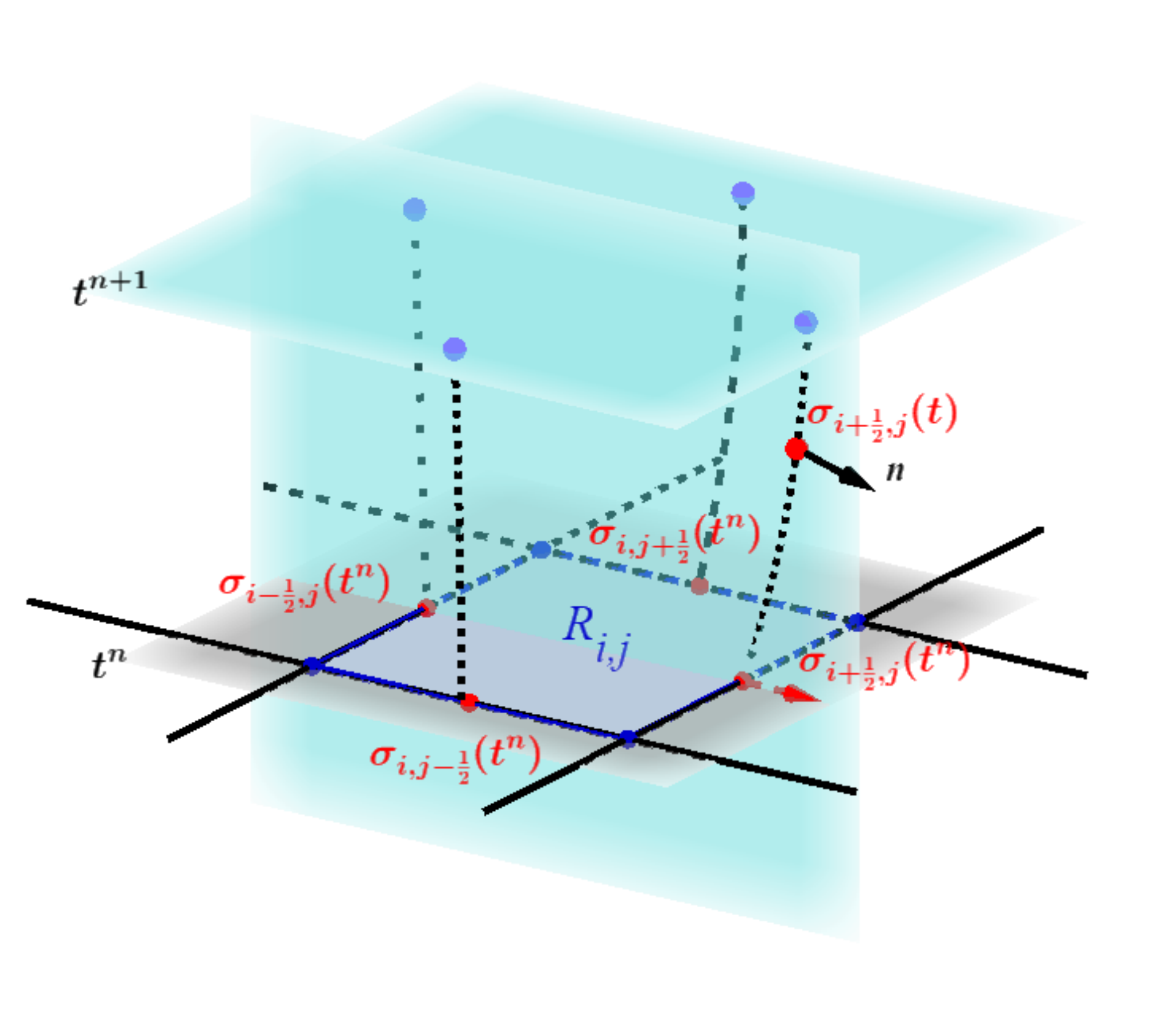}
	\caption{No flow surface region (left). Normal 
         vector $\vec{\textsl{n}}$ in $y = y_j$-plane (right).} 
	\label{fIII}
\end{figure}
			
{\color{black}
\noindent{\bf Remark}: We point out that solutions for 
the generalized ODE system (\ref{syst2}) to compute 
$\sigma_{i-1/2,j}^n(t)$ in Eq.(\ref{auxSigma}) by the 
differential equation 
$\frac{d \sigma_{i-1/2,j}^n(t)}{dt}=\frac{f(U_{i-\frac{1}{2},j}^n)}{U_{i-\frac{1}{2},j}^n}$ 
on the edge of the {\it no flow surface region} (see Figure \ref{fIII} and 
Figure \ref{fIV}) can be facilitated by suitable linear resconstrution 
$L(x,t)$ such as $U_{i-\frac{1}{2},j}  = 
\displaystyle\frac{1}{h}\int_{x_{i-1,j}^n}^{x_{i,j}^n} L(x,t) dx$.
The situation is similar for quantities $U_{i+\frac{1}{2},j}^n$,
$U_{i,j-\frac{1}{2}}^n$ and $U_{i,j+\frac{1}{2}}^n$;
see Section \ref{ImODE} for details.}

We choose the simplest approximation of 
system \eqref{syst}, by setting $U_{i,j}^n$ 
at $t=t^n$ we get:

\begin{equation}
 \begin{array}{l l l}
  \begin{cases}
    \gamma_1(t)'(t) = 
\frac{f\left(U_{i+\frac{1}{2},j}^n\right)}{U_{i+\frac{1}{2},j}^n}, \\
    \gamma_1(t)'(t^n) = x_{i+\frac{1}{2}},
  \end{cases}  
      \!\!\! & \!\!\!
  \begin{cases}
    \alpha_2(t)'(t) = 
\frac{g\left(U_{i,j-\frac{1}{2}}^n\right)}{U_{i,j-\frac{1}{2}}^n}, \\
    \alpha_2(t)'(t^n) = y_{j-\frac{1}{2}},
  \end{cases}
      \!\!\! & \!\!\!
  \begin{cases}
    \theta_2(t)'(t) = 
\frac{g\left(U_{i,j+\frac{1}{2}}^n\right)}{U_{i,j+\frac{1}{2}}^n}, \\
    \theta_2(t)'(t^n) = y_{j+\frac{1}{2}},
\end{cases}
\end{array}
\label{syst2}
\end{equation}
where
{\color{black}
\begin{equation}
\frac{f(U_{i+\frac{1}{2},j}^n)}{U_{i+\frac{1}{2},j}^n} 
\equiv f_{i+\frac{1}{2},j}, \,\,
\frac{f(U_{i-\frac{1}{2},j}^n)}{U_{i-\frac{1}{2},j}^n} 
\equiv f_{i-\frac{1}{2},j}, \,\,
\frac{g(U_{i,j-\frac{1}{2}}^n)}{U_{i,j-\frac{1}{2}}^n} 
\equiv g_{i,j-\frac{1}{2}}\mbox{ and } 
\frac{g(U_{i,j+\frac{1}{2}}^n)}{U_{i,j+\frac{1}{2}}^n} 
\equiv g_{i,j+\frac{1}{2}}.
\label{auxSigma}
\end{equation}
}
Thus, we can approximate curves of the {\it no flow surface region} at 
$t^n < t < t^{n+1}$ as: 
$$\sigma_1(t) \approx {\color{black}x_{i-\frac{1}{2}}} + (t - t^n)f_{i-\frac{1}{2},j},$$ 
$$\gamma_1(t) \approx x_{i+\frac{1}{2}} + (t - t^n)f_{i+\frac{1}{2},j},$$
$$\alpha_2(t) \approx y_{j-\frac{1}{2}} + (t - t^n)g_{i,j-\frac{1}{2}},$$
 and $$\theta_2(t) \approx y_{j+\frac{1}{2}} + (t - t^n)g_{i,j+\frac{1}{2}}.$$
The approximation of the volume ${\color{black}D_{i,j}^n}$ gives (see right 
frame in Figure \ref{fIV}):
$${\color{black}D_{i,j}^n} = \{ (t,x,y)  /  t^n \leq t < t^{n+1}, \,\,\, 
\sigma_1(t) \leq x < \gamma_1(t), \,\,\, \alpha_2(t) 
\leq y < \theta_2(t) \}.$$

The new conservative Lagrangian-Eulerian scheme 
is given by very simply formulas:

\vspace{2mm}
{\bf STEP I} (Lagrangian Evolution, see Figure \ref{fIV}, 
{\color{black} and below $h\equiv\Delta y = \Delta x$})
\begin{equation}
  \overline{U}_{i,j}^{n+1}  = 
	\frac{A(R_{i,j}^{n})}{A(\overline{R}_{i,j}^{n+1})} U_{i,j}^n,
        \quad
        \text{ with }
    A(R_{i,j}^{n}) = {\color{black}h^2}  
        \text{ and }
        A(\overline{R}_{i,j}^{n+1}) = h_i^n * h_j^n
  \label{ReltnTnp2}
\end{equation}
where $h_i^n * h_j^n = ({\color{black}\Delta x} - (f_{i-1/2,j} + f_{i+1/2,j})
\Delta t)*({\color{black}\Delta y} - (g_{i,j-1/2} + g_{i,j+1/2})\Delta t)$.

\vspace{2mm}
{\bf STEP II} (Eulerian Projection, see Figure \ref{fV}) 
\begin{equation}
{\color{black}U_{i,j}^{n+1}} = \frac{1}{A(\overline{R}_{i,j}^{n+1})} (C1+C2+C3),
\label{proy}	
\end{equation}
\noindent where						
$C1\equiv c_{11}\overline{U}_{i-1,j-1}^{n+1} + 
c_{12}\overline{U}_{i,j-1}^{n+1} + c_{13}\overline{U}_{i+1,j-1}^{n+1}$,
$C2\equiv c_{21}\overline{U}_{i-1,j}^{n+1} + 
c_{22}\overline{U}_{i,j}^{n+1}  + c_{23}\overline{U}_{i+1,j}^{n+1}$ and
$C3\equiv c_{31}\overline{U}_{i-1,j+1}^{n+1} + 
c_{32}\overline{U}_{i,j+1}^{n+1} + c_{33}\overline{U}_{i+1,j+1}^{n+1}$.

\subsection{Improving numerically the solution of 
the generalized ODE system (\ref{syst2}) with 
conservation and robustness}\label{ImODE} Solutions 
$\sigma_{i-1/2,j}^n(t)$ of the differential system 
can be obtained using the approximations
\begin{equation}
\begin{array}{ll}
U_{i-\frac{1}{2},j} & = 
\displaystyle\frac{1}{h}\int_{x_{i-1,j}^n}^{x_{i,j}^n} L(x,t) dx 
=\frac{1}{h}\left( \int_{x_{i-1,j}^n}^{x_{i-\frac{1}{2},j}^n} 
L_{i-1,j}(x,t) dx
+ \int_{x_{i-\frac{1}{2},j}^n}^{x_{i,j}^n} 
L_{i,j}(x,t) dx \right) \\ \\
& = \displaystyle\frac{1}{2} 
(U_{i-1,j} + U_{i,j}) + \frac{1}{8} (U'_{i,j} - U'_{i-1,j}).
\end{array}   
\end{equation}
\noindent Additional  and even high-order approximations 
are also acceptable for $\frac{d \sigma_{j-1/2}^n(t)}{dt}=
\frac{f(u)}{u}$ (see  (\ref{syst2})). As in 
\cite{CNMAC2016}, the piecewise constant numerical 
data is reconstructed into a piecewise linear 
approximation through the use of MUSCL-type interpolants:
$L_{i,j}(x,t) = u_{i,j}(t)+(x-x_j)\frac{1}{\Delta x} u'_{i,j}$.
For the numerical derivative $\frac{1}{\Delta x}u'_{i,j}$, 
there are several choices of slope limiters for scalar 
case; in the book \cite{LVFEB02} there is a good compilation 
of many possible options. Finally, in order to show the flexibility 
of the reconstruction we use the nonlinear Lagrange 
polynomial in $U_{i-1,j}^n$, $U_{i,j-1}^n$, $U_{i,j}^n$, 
$U_{i,j+1}^n$ and $U_{i+1,j}^n$. Therefore, equation 
(\ref{ReltnTnp2}) reads
\begin{equation}
\overline{U}_{i,j}^{n+1}  =  
\displaystyle\frac{1}{h_j^{n}}
\displaystyle\int_{s_{j-\frac{1}{2}}^{n}}^{s_{j+\frac{1}{2}}^{n}}P_2(x,y)ds,
\label{Reco}
\end{equation}
where $s = x,y$ and $P_2(x,y) = U_{i-1,j}^n \, 
L_{-1}(x - x_{i}) + U_{i,j-1}^n \, L_{-1}(y - y_{j}) 
+ U_{i,j}^n \,
L_{0}(x - x_{i})+ U_{i,j}^n \,
L_{0}(y - y_{j}) + U_{i+1,j}^n\, L_1(x - x_{j}) + 
U_{i,j+1}^n\, L_1(y - y_{j})$ and 
\begin{equation}
L_{\pm \,1}(x) = \frac{1}{2}\left[\left(\frac{x}{h} 
\pm \frac{1}{2}\right)^2
- \frac{1}{4}\right], \quad \quad \quad  
L_0(x) = 1- \left(\frac{x}{\Delta x}\right)^2.
\label{Lagran}
\end{equation}

\subsection{A Lagrangian-Eulerian CFL stability 
constraint} Next, by definition of $f_{i-1/2,j}$, 
$f_{i+1/2,j}$, $f_{i,j-1/2}$ and $f_{i,j+1/2}$ 
in (\ref{syst2}),  we obtain the resulting 
coefficients of the Eulerian projection formula 
(\ref{proy}) as follows. Let the vectors 
$$C_x = \left[C_{xl}, {\color{black}\Delta x} - C_{xl} - C_{xr},C_{xr}\right],$$ 
where $C_{xl} = 0.5(1 + \mbox{sign}(f_{i-1/2,j}))f_{i-1/2,j}\Delta t$, 
$C_{xr} = 0.5(1 - \mbox{sign}(f_{i+1/2,j}))f_{i+1/2,j}\Delta t$ and
$$C_y = \left[C_{yr}, {\color{black}\Delta y} - C_{yr} - C_{yr},C_{yr}\right],$$
where $C_{yl} = 0.5(1 + \mbox{sign}(g_{i,j-1/2}))g_{i,j-1/2}\Delta t$, 
$C_{yr} = 0.5(1 - \mbox{sign}(g_{i,j+1/2}))g_{i,j+1/2}\Delta t$.
			
We define the coefficients of projection formula 
(\ref{proy}) as the coefficients of the matrix
\begin{equation}
C = (c_{i,j}) = C_x^T  C_y,
\label{matrixcfl}
\end{equation}
under the CFL-condition ({\color{black}along with $h \equiv \Delta x = \Delta y$})
\begin{equation}
\max_{i,j} \{{\color{black}\left|f'(U_{i,j}^n)\right|},{\color{black}\left|g'(U_{i,j}^n)\right|}, \left|f_{i-1/2,j}\right|, \left|f_{i+1/2,j}\right|, 
\left|g_{i,j-1/2}\right|, \left|g_{i,j+1/2}\right| \} \Delta t 
< \frac{{\color{black} h}}{2}.
\label{matrixcfl2}
\end{equation}


\subsection{A connection with monotone convergent 
entropy stable numerical scheme}\label{monoLE} 
For the purpose of this work, we  invoke the 
solid theory of general monotone difference schemes 
(see, e.g., \cite{EGH95a,EGH95b,CCH99,CRMAJ80}) to 
illustrate the generality or our Lagrangian-Eulerian 
approach (\ref{ReltnTnp2})-(\ref{proy}), under 
the CFL stability constraint (\ref{matrixcfl2}). 
By definition of matrix, 
\begin{equation}
A_{i,j} = \left[\begin{array}{lll}
\overline{U}_{i-1,j-1}^{n+1} & \overline{U}_{i,j-1}^{n+1} 
& \overline{U}_{i+1,j-1}^{n+1} \\
\overline{U}_{i-1,j}^{n+1}   & \overline{U}_{i,j}^{n+1}   
& \overline{U}_{i+1,j}^{n+1} \\
\overline{U}_{i-1,j+1}^{n+1} & \overline{U}_{i,j+1}^{n+1} 
& \overline{U}_{i+1,j+1}^{n+1}
\end{array}\right],
\label{MatrixProj}
\end{equation}
the Eulerian Projection step (\ref{proy}) over 
original grid may be recast as
\begin{equation}	
U_{i,j}^{n+1} = C_x A_{i,j}^{T} C_y^T,
\label{proyF}
\end{equation}
or in the form of conservative monotone scheme as,
\begin{equation}	
U_{i,j}^{n+1} = U_{i,j}^n - \lambda ^x \Delta _{+}^{x} 
     F(U_{i-p,i-r}^n,...,U_{i+q+1,j+s+1}^n) 
     - \lambda ^y \Delta _{+}^{y} 
      G(U_{i-p,j-r}^n,...,U_{i+q+1,j+s+1}^n)
\label{proyF1}
\end{equation}
where, along with (\ref{ReltnTnp2}) and taking 
$h = \Delta x = \Delta y$, we have, 
\begin{itemize}
\item[] $F(U_{i-1,j-1}^n,...,U_{i+1,j+1}^n) =$ 
\item[] \quad $F_R(U_{i,j-1}^n,U_{i-1,j-1}^n,
U_{i-1,j}^n,U_{i,j}^n,U_{i,j-1}^n) - 
F_L(U_{i-1,j+1}^n,,U_{i-1,j}^n,U_{i,j-1}^n,U_{i,j}^n,U_{i,j+1}^n),$
\end{itemize}
where
\begin{itemize}
\item[] \qquad $F_R = h_y C_{xl}  
\left( \overline{U}_{i+1,j}^{n+1} - 
\overline{U}_{i,j}^{n+1} \right)-$ 
\item[] \qquad $C_{xl} C_{yr} 
\left(\overline{U}_{i-1,j+1}^{n+1} - 
\overline{U}_{i-1,j}^{n+1} - ( \overline{U}_{i,j+1}^{n+1} 
- \overline{U}_{i,j}^n ) \right) + 
\frac{1}{2}(f(U_{i-1,j}^n)+f(U_{i+1,j}^n)),$
\end{itemize}
\begin{itemize}
\item[] \qquad $F_L = h_y C_{xl} C_{yl}
\left( \overline{U}_{i-1,j}^{n+1} - 
\overline{U}_{i,j}^{n+1} \right)-$ 
\item[] \qquad $C_{xl} C_{yl} 
\left(\overline{U}_{i-1,j-1}^{n+1} - 
\overline{U}_{i,j-1}^{n+1} - ( \overline{U}_{i-1,j}^{n+1} - 
\overline{U}_{i,j}^{n+1} )\right)  + 
\frac{1}{2}(f(U_{i-1,j}^n)+f(U_{i+1,j}^n)).$
\end{itemize}

\noindent and

\begin{itemize}
\item[] $G(U_{i-1,j-1}^n,...,U_{i+1,j+1}^n) =$
\item[] \quad $G_R(U_{i-1,j}^n,U_{i,j}^n,
U_{i+1,j+1}^n,U_{i,j+1}^n,U_{i+1,j}^n) - 
G_L(U_{i+1,j}^n,,U_{i,j}^n,U_{i,j-1}^n,U_{i+1,j}^n,U_{i,j+1}^n),$
\end{itemize}

\begin{itemize}
\item[] \qquad $G_R = h_x C_{yl}  
\left( \overline{U}_{i,j+1}^{n+1} - 
\overline{U}_{i,j}^{n+1} \right)-$
\item[] \qquad $C_{xr}*C_{yl} 
\left(\overline{U}_{i+1,j-1}^{n+1} - 
\overline{U}_{i+1,j}^{n+1} - ( 
\overline{U}_{i,j-1}^{n+1} - 
\overline{U}_{i,j}^{n+1} ) \right) + 
\frac{1}{2}(g(U_{i,j-1}^n)+g(U_{i,j+1}^n)),$
\end{itemize}

\begin{itemize}
\item[] \qquad $G_L = h_x C_{yl}  
\left( \overline{U}_{i,j-1}^{n+1} - 
\overline{U}_{i,j}^{n+1} \right)-$
\item[] \qquad $C_{xr}*C_{yr}  
\left(\overline{U}_{i+1,j+1}^{n+1} - 
\overline{U}_{i+1,j}^{n+1} - ( 
\overline{U}_{i,j+1}^{n+1} - 
\overline{U}_{i,j}^{n+1} )\right)  +  
\frac{1}{2}(g(U_{i,j-1}^n)+g(U_{i,j+1}^n)).$
\end{itemize}

We can note that, $F$ and $G$ satisfy condition 
(\ref{sclpcons2d}), this implies 
consistentcy with (\ref{2D1H}) and thus the 
numerical method to 2D-hyperbolic equations is monotone. 

In order for the above scheme be consistent with 
(\ref{2D1H}), we must have:
%
\begin{equation}
 F_1(u,u,...,u) = f(u) \quad \text{and} \quad
 F_2(u,u,...,u) = g(u), \quad u \in \mathds{R}.
 \label{sclpcons2d}
\end{equation}
Here, the functions $F_1$ and $F_2$, are the corresponding numerical 
fluxes of the perninent approximation. The difference approximation 
is monotone on the interval $\left[a, b\right]$ if $\mathcal{G}$ a nondecreasing 
function of each argument $U_{i,j}^n$ so long as all arguments lie 
in $\left[a, b\right]$. Write $u(x, y, t) = (S(t)u0)(x, y)$, where 
$S(t) :L^1(\mathds{R}^2) \cap L^{\infty}(\mathds{R}^2)
\rightarrow L^1(\mathds{R}^2) \cap L^{\infty}(\mathds{R}^2)$
for each $t \geq 0$ and $t \rightarrow S(t)u_0$ is continuous into 
$L^1(\mathds{R}^2)$. To compute this solution numerically we set 
\begin{equation}
 u^{\Delta \, t} = \sum_{n=0}^{\infty} \sum_{k = -\infty}^{\infty} 
 U_{j,k}^n \, \mathcal{X}_{j,k}^n,
 \label{concep8defu2d}
\end{equation}
where $\mathcal{X}_{j,k}^n$ is the characteristic 
function in the respective cell. Indeed, it turns 
out that conservative monotone schemes converge to 
entropy solutions. Therefore, convergence toward 
the entropy solution to out our Lagrangian-Eulerian 
approach (\ref{ReltnTnp2})-(\ref{proy}) is proven. 
In \cite{AAP19}, we were able to establish entropy 
convergence and error estimates for conservative 
Lagrangian-Eulerian method on triangular grids.
	
\subsection{Numerical experiments with the 
Lagrangian-Eulerian with conservation properties}

We present a benchmark comprehensive set of numerical tests which 
explore the role of accuracy of our new 2D Lagrangian-Eulerian 
scheme with conservation properties. 

\subsubsection{A numerical convergence study 
for a linear 2D advection flow model}

\begin{figure}[ht!]
\centering
\includegraphics[width=0.45\linewidth]{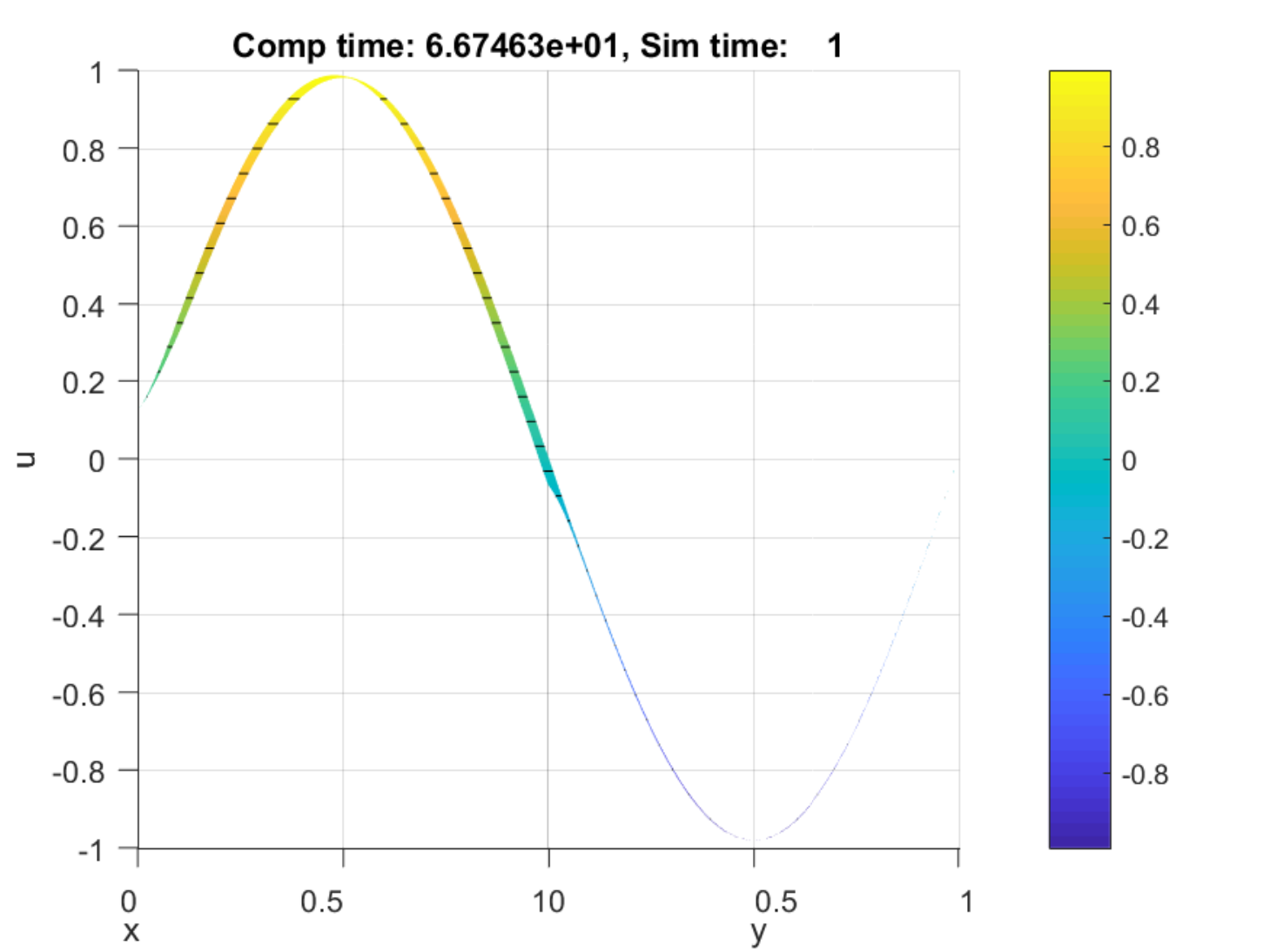}
\includegraphics[width=0.45\linewidth]{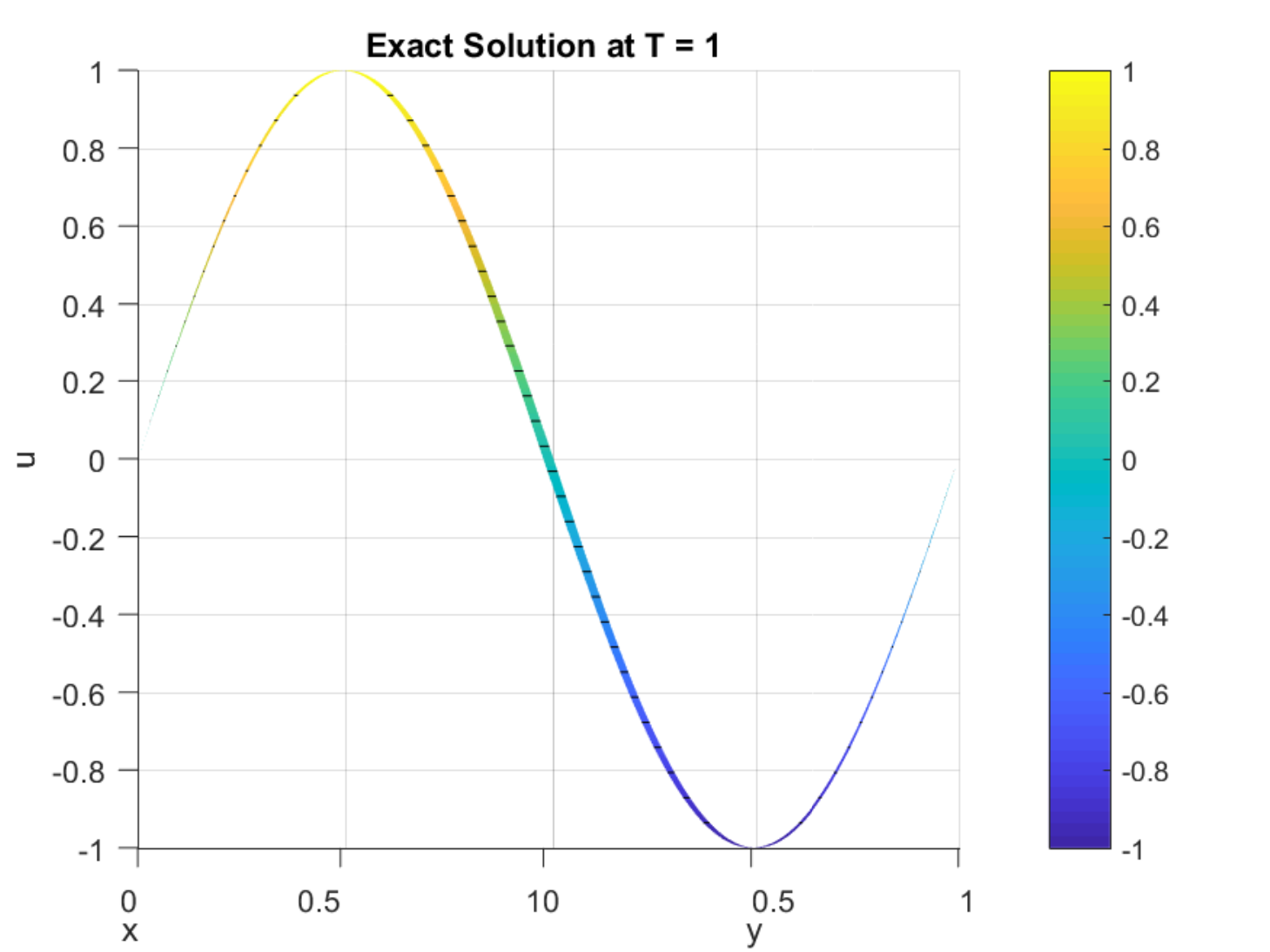}
\caption{Computed solution by our 2D Lagrangian-Eulerian scheme
(\ref{ReltnTnp2})-(\ref{proy}), under the CFL stability constraint 
(\ref{matrixcfl2}) to problem (\ref{linear1})-(\ref{linear2}) 
at simulation time $t=1$ (left) and exact solution (right) both 
projected over oblique plane $x=y$.}
\label{gridref128}
\end{figure}

\begin{figure}[ht!]
\includegraphics[scale=0.33]{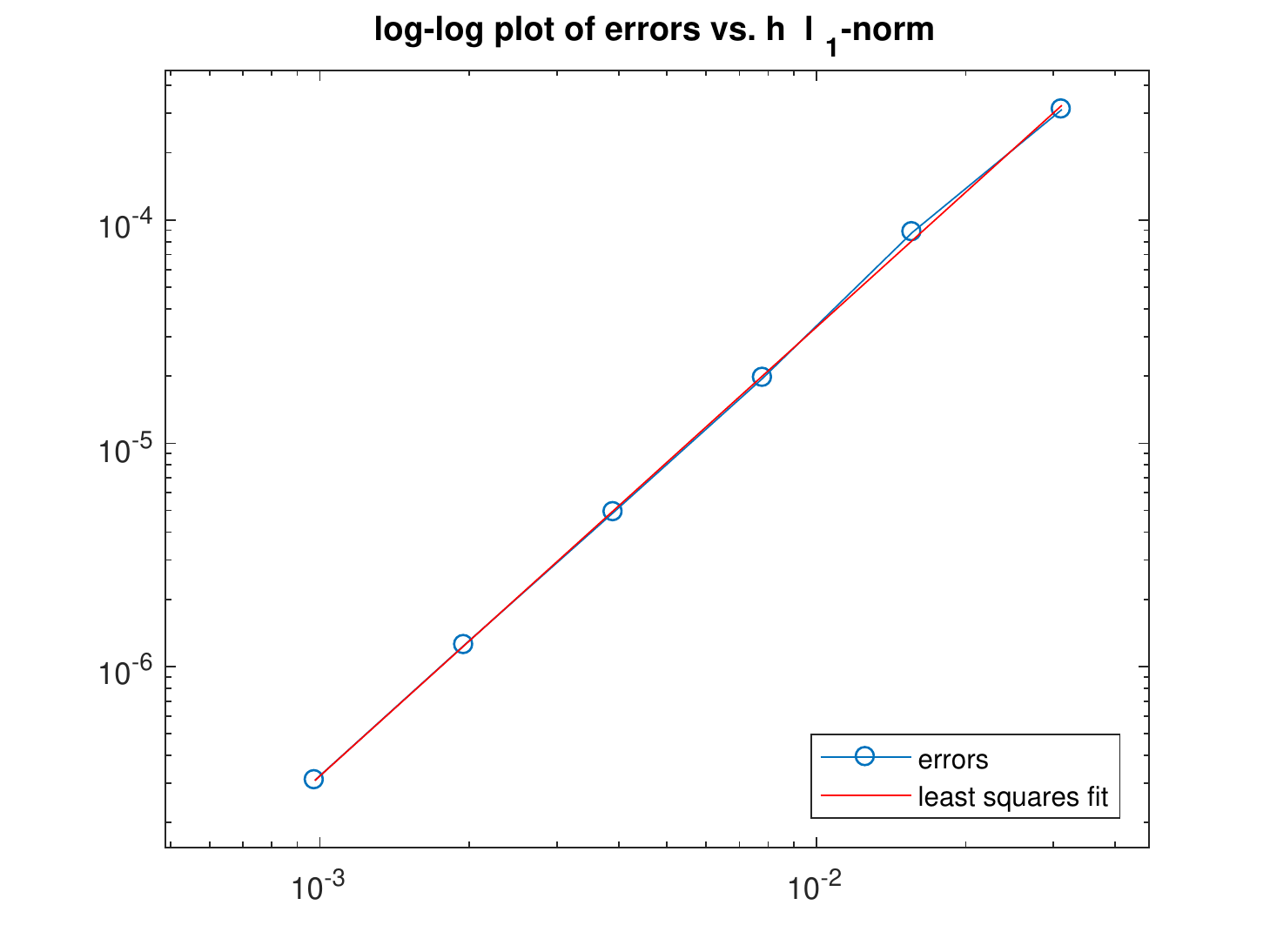}
\includegraphics[scale=0.33]{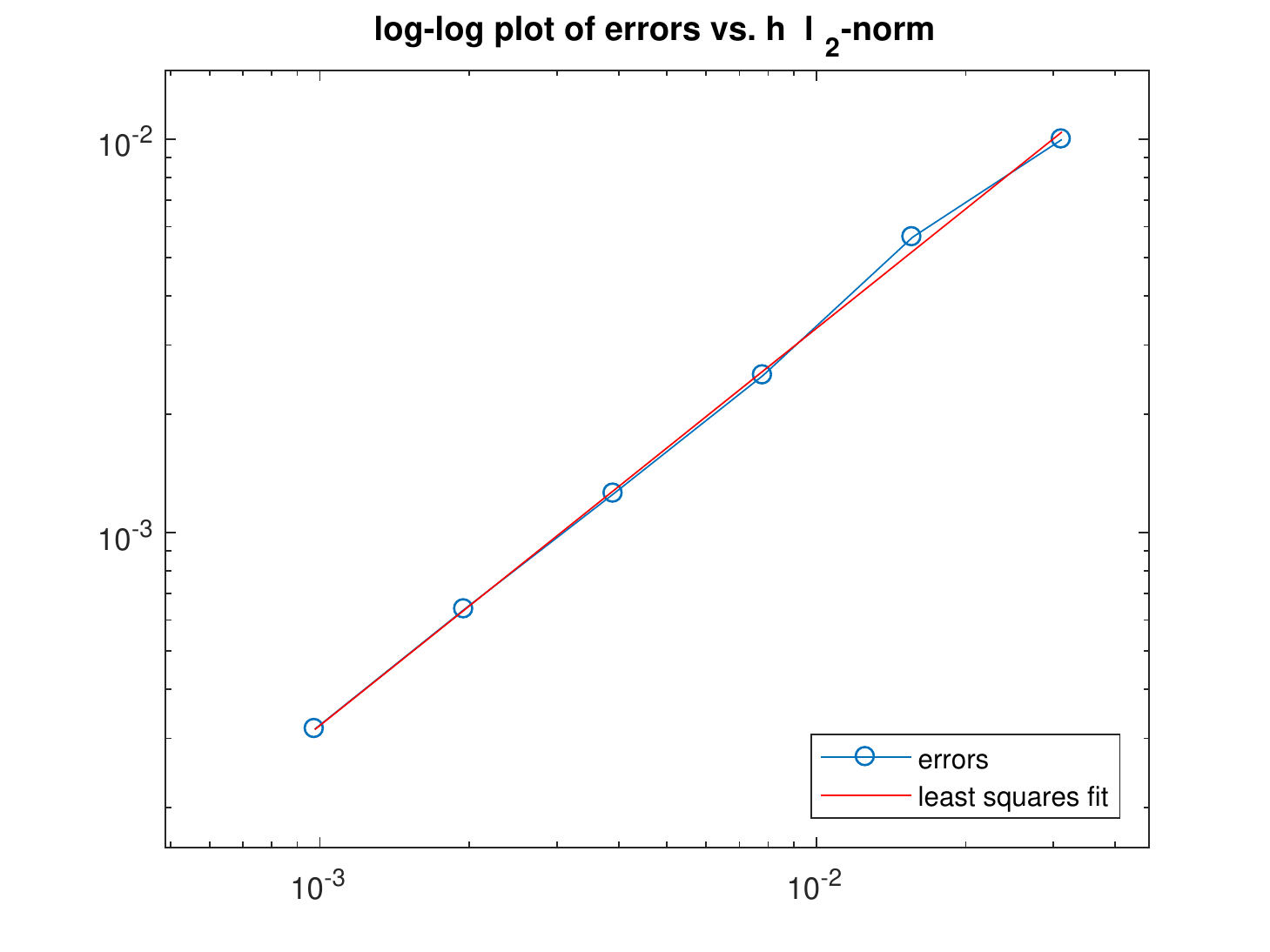}        
\includegraphics[scale=0.33]{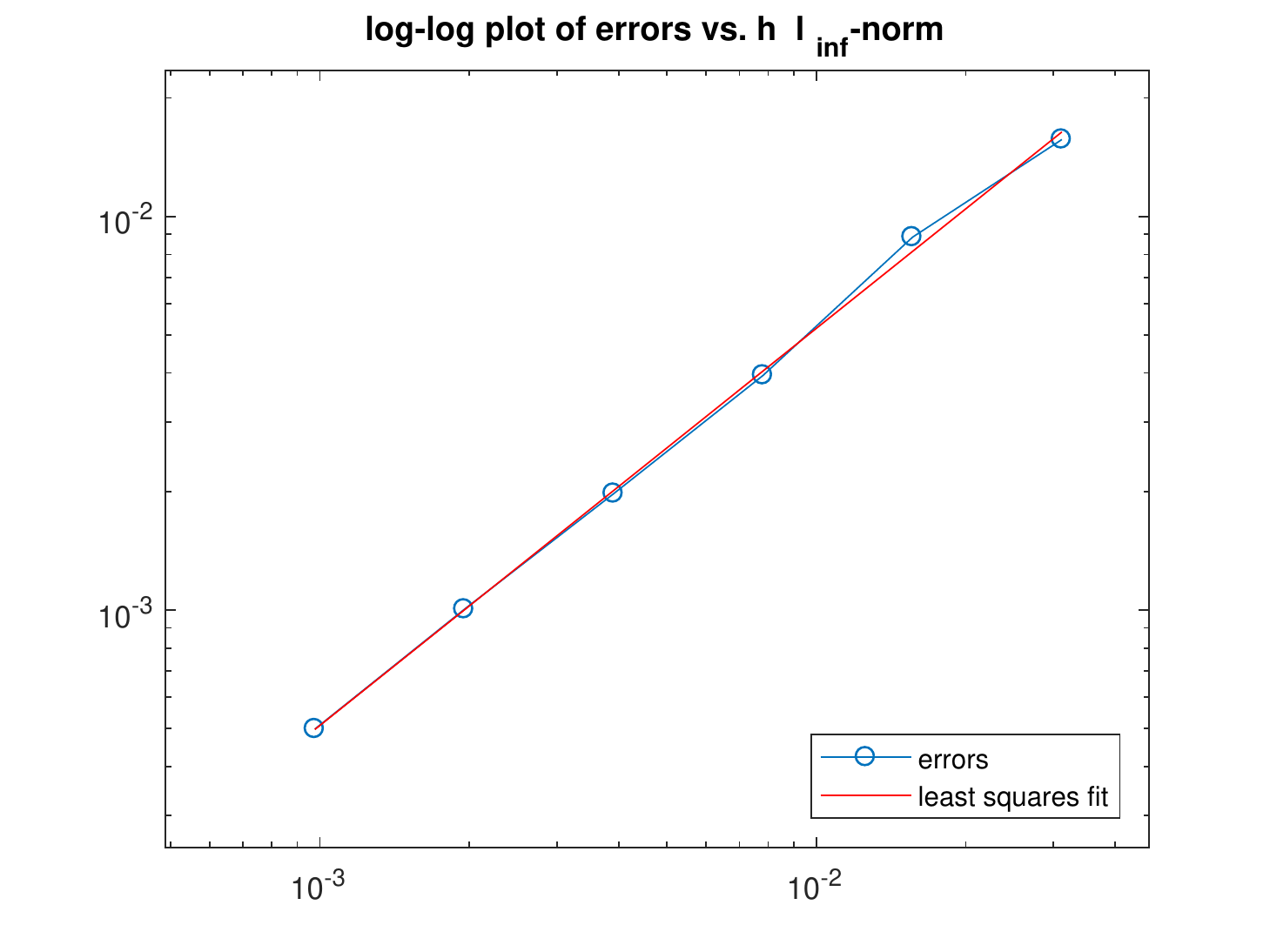}
\centering
\begin{tabular}{c c c c c c c c}
  \hline
    Cells &  $h$  &   $\|u - U\|_{l_h^{1}}$ & Order  & $\|u - U\|_{l_h^{2}}$ & Order &  $\|u - U\|_{l_h^{\infty}}$ & Order \\
    \hline
    $64  \times 64$ &  $0.016$  & $5.156\times10^{-2}$ & -- & $5.573\times10^{-2}$ & -- & $ 1.339\times10^{-1}$ & --  \\
    $128 \times 128$ & $0.007$ & $2.046\times10^{-2}$ &  $1.333$ & $2.493\times10^{-2}$ & $1.161$ & $6.509\times10^{-2}$  & $1.040$ \\
    $256 \times 256$ & $0.004$ & $1.309\times10^{-2}$ & $0.644$ & $1.467\times10^{-2}$ & $0.765$ & $3.761\times10^{-2}$  & $0.791$\\
    $512 \times 512$ & $0.002$ & $6.090\times10^{-3}$ & $1.103$ & $7.034\times10^{-3}$ & $1.061$ & $1.835\times10^{-2}$  & $1.036$\\
   \hline
	 LSF $E(h)$ & & $2.911 * h^{0.989}$  & &  $3.051 * h^{0.972}$  & &  $6.534 * h^{0.939}$ & \\
 \hline
   \label{T2DLin}
  \end{tabular}
\caption{In the table are shown
errors between the numerical approximations 
($U$) and exact solutions ($u$) in $l_h^1$, $l_h^2$ 
and $l_h^{\infty}$ norms to problem (\ref{linear1}) 
with initial condition $u(x,0) = \sin(\pi(x+y))$, 
advanced from $t=0$ to $t=1$ along with CFL condition 
$0.67$. 
}
\label{REFnorm2}
\end{figure}

First, let us consider the following initial value problem
\begin{equation}
   \displaystyle\frac{\partial u}{\partial t} 
   + \frac{\partial u}{\partial x}
   + \frac{\partial u}{\partial y}= 0,
   \quad 
   \hbox{ in the computational domain } 
   \quad (x,y,t) \in \left[0,1\right] \times \left[0,1\right] \times \left[0,1\right],
   \label{linear1}
\end{equation}
and initial condition,
\begin{equation}
u(x,y,0) = \sin(\pi(x+y)).
\label{linear2}
\end{equation}
\begin{figure}[ht!]
\centering
\includegraphics[width=0.45\linewidth]{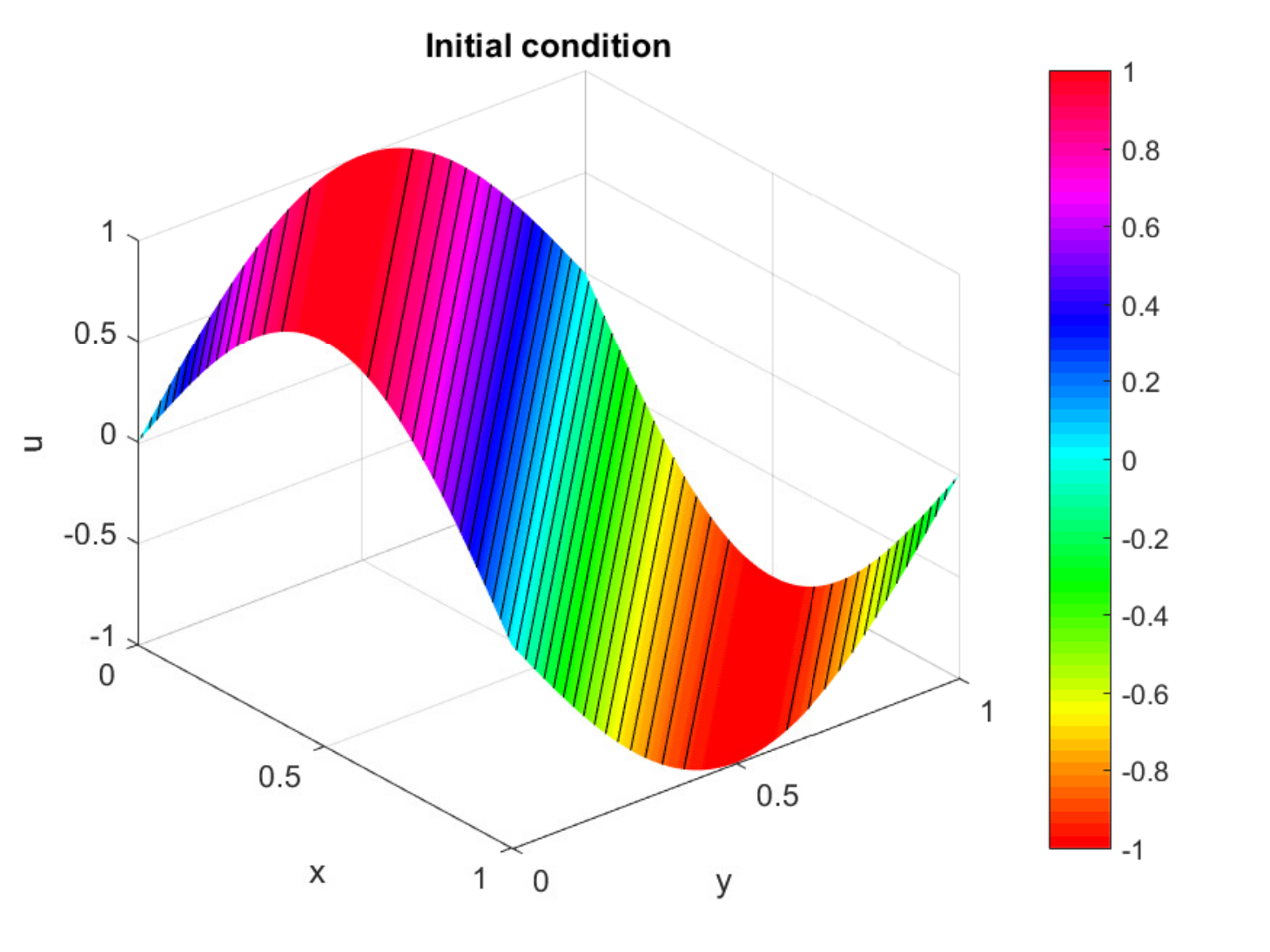}
\includegraphics[width=0.45\linewidth]{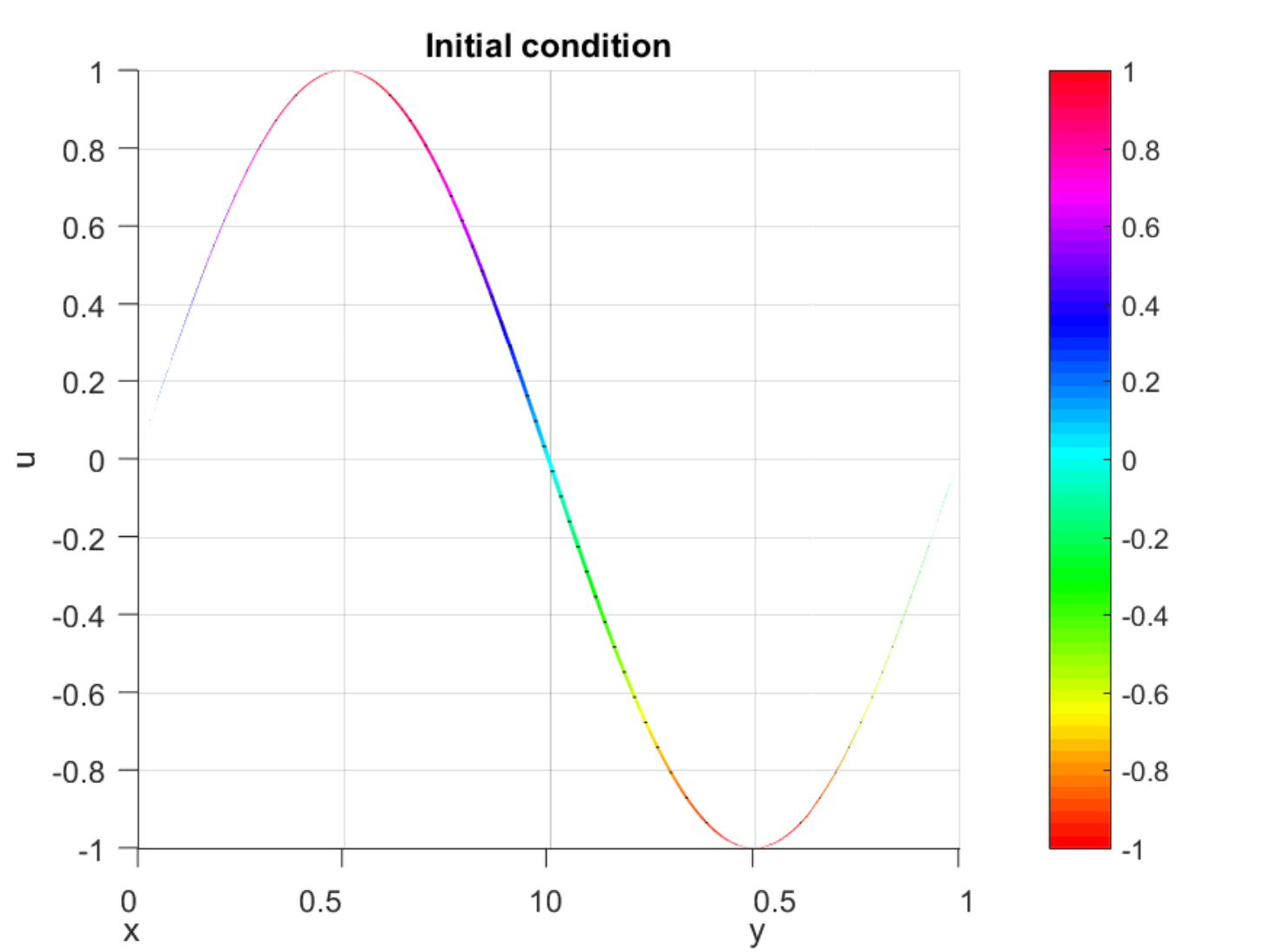}
\caption{Initial condition for problem 
(\ref{linear1})-(\ref{linear2}). On the left 
(resp. right) picture we show a ``3D-plot's view 
angle'' (resp. a oblique projection over the plane $x=y$).}
\label{intcond}
\end{figure} 
It is easy exercise to show that the exact solution to problem
(\ref{linear1})-(\ref{linear2}) is $u(x,y,t) = \sin(\pi(x+y-2t))$. 
The solution will be advanced from $t=0$ to $t=1$ and we notice 
that at this time the solution is merely traslated by one period 
$2\pi$, with respect to (\ref{linear2}) in the oblique $x=y$ direction.
The approximation computed with our scheme to problem 
(\ref{linear1})-(\ref{linear2}) is shown in the Figure \ref{gridref128} 
(left frame) along with the exact solution on the right frame.
In Figure \ref{REFnorm2} we observe numerical convergence rates to 
(\ref{linear1})-(\ref{linear2}) in the $l_h^1$-norm (left), 
in the $l_h^2$-norm (middle) and in the $l_h^{\infty}$-norm 
(right) computed by scheme (\ref{ReltnTnp2})-(\ref{proy}).

\subsubsection{An oblique Riemann problem for 
inviscid 2D Burgers's equation}

\begin{figure}[ht!]
\centering
\includegraphics[width=0.2\textwidth]{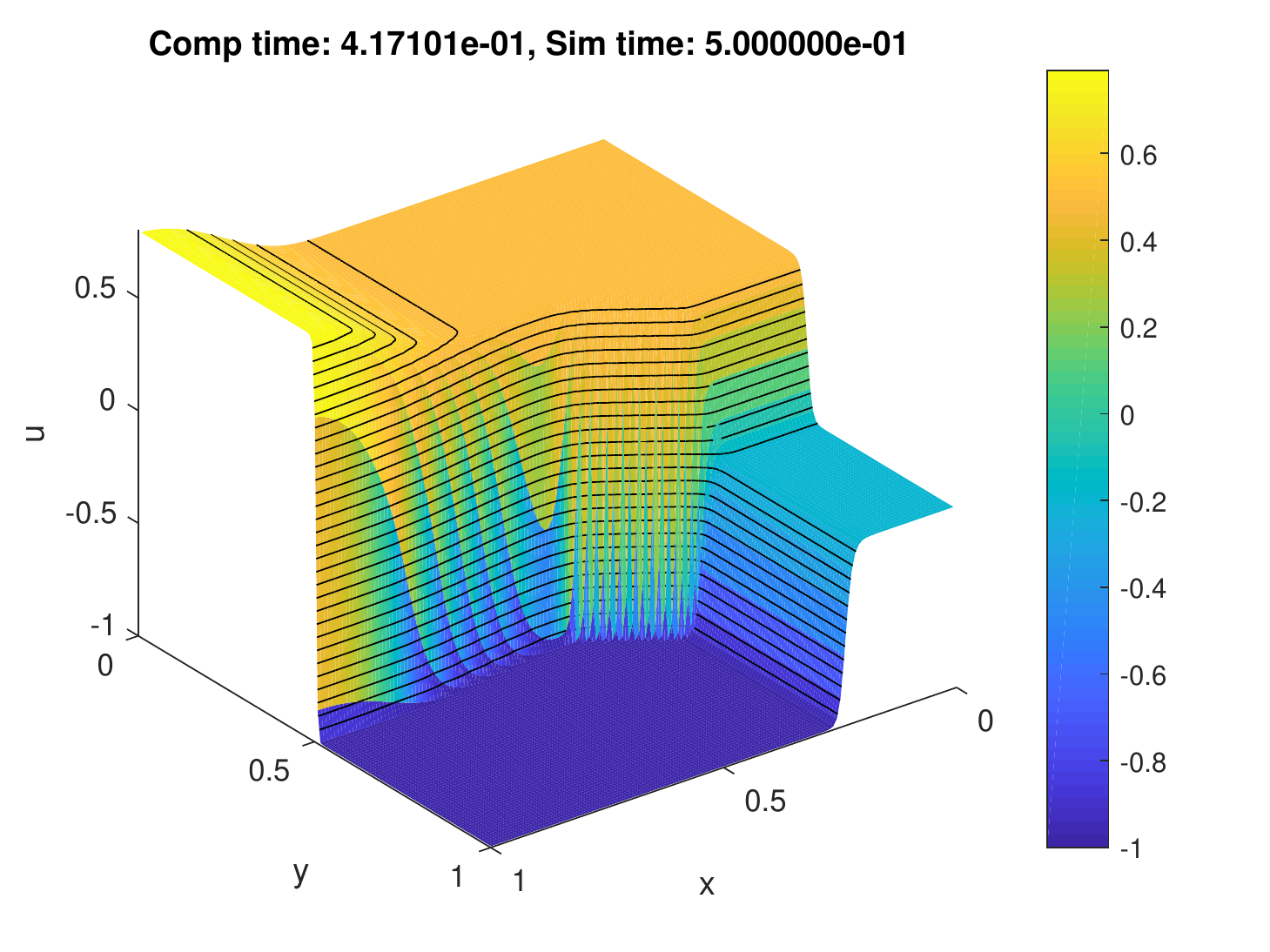}
\includegraphics[width=0.2\textwidth]{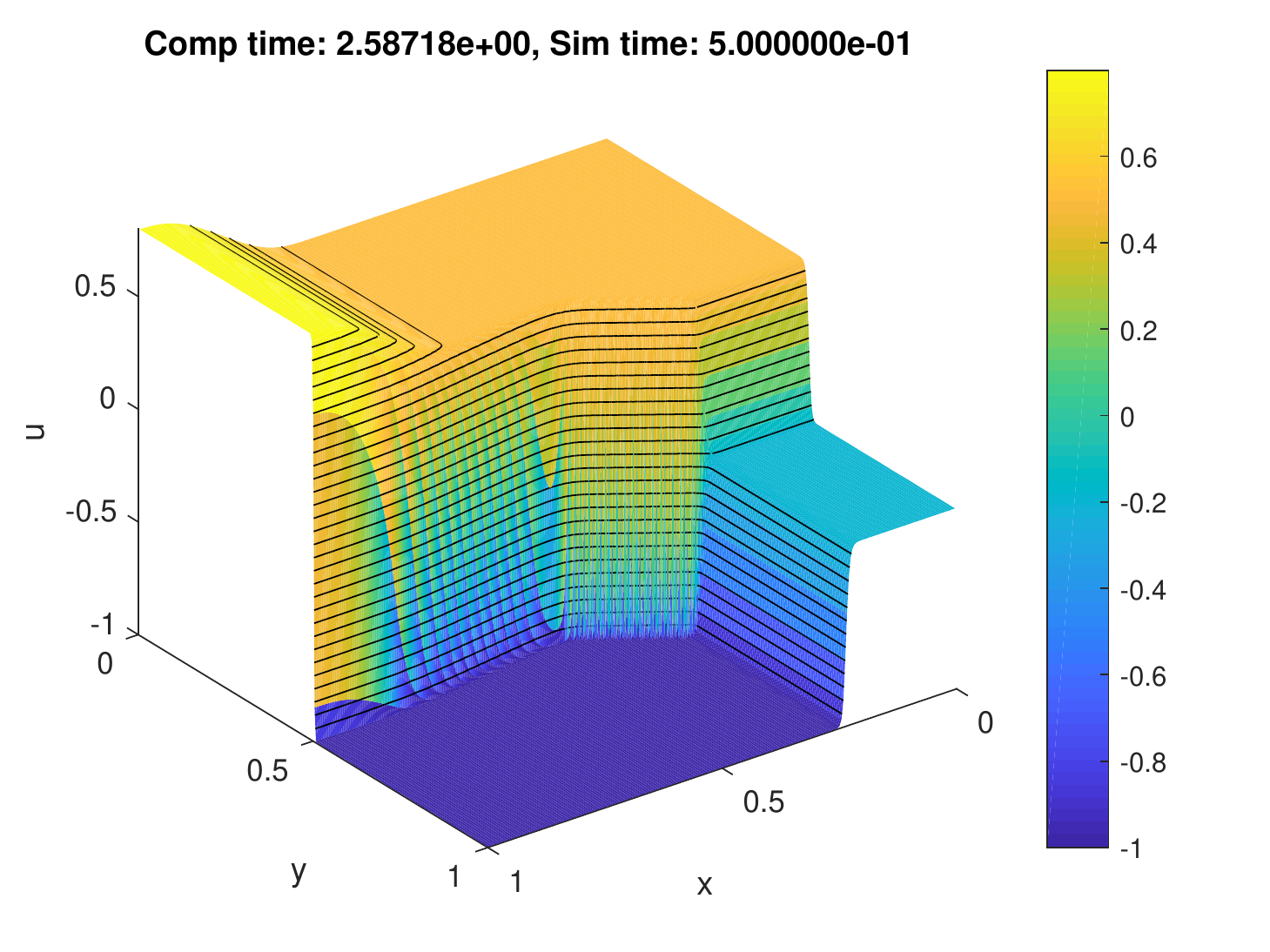}
\includegraphics[width=0.2\textwidth]{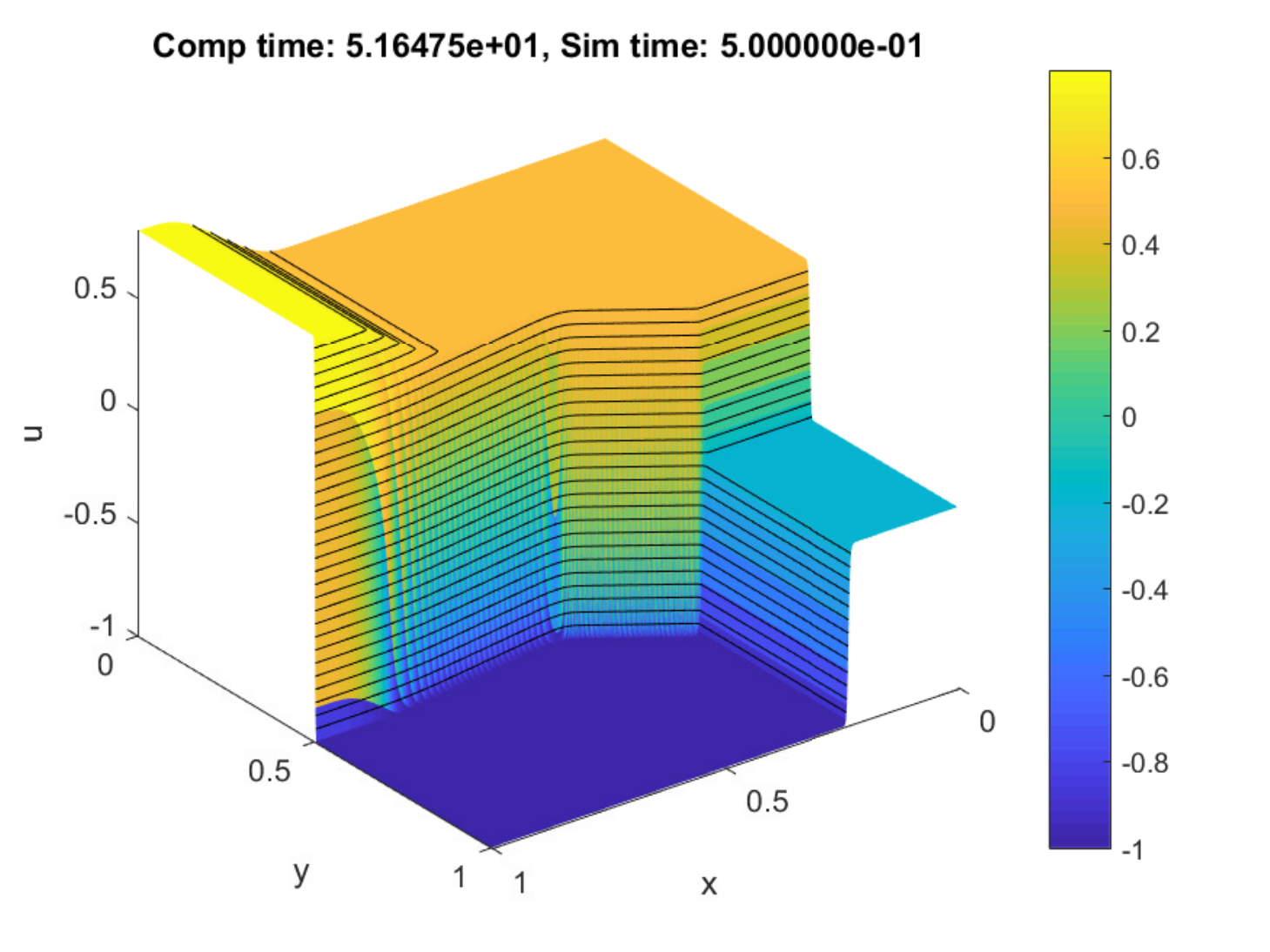} 
\includegraphics[width=0.2\textwidth]{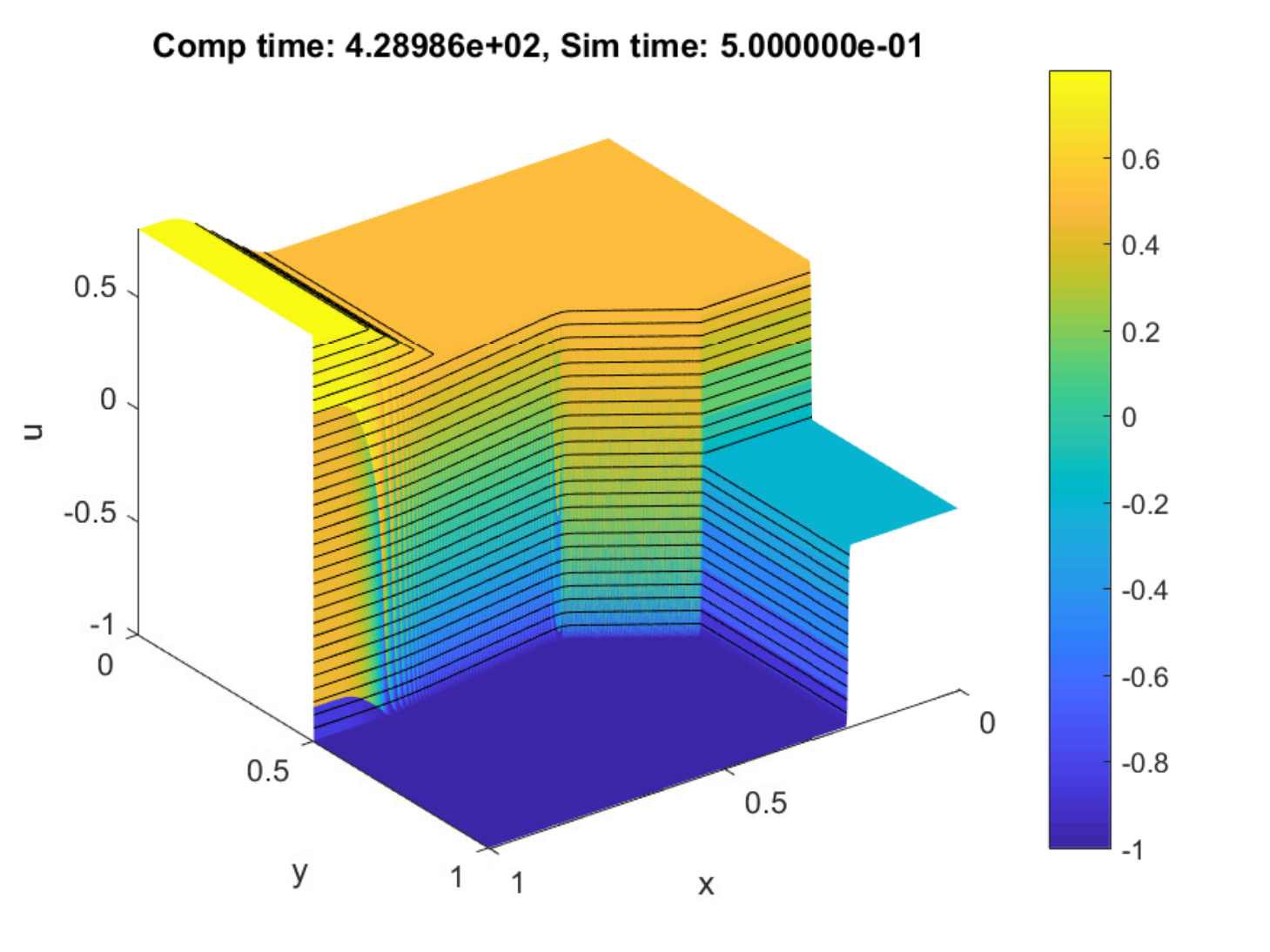}\\
\includegraphics[width=0.2\textwidth]{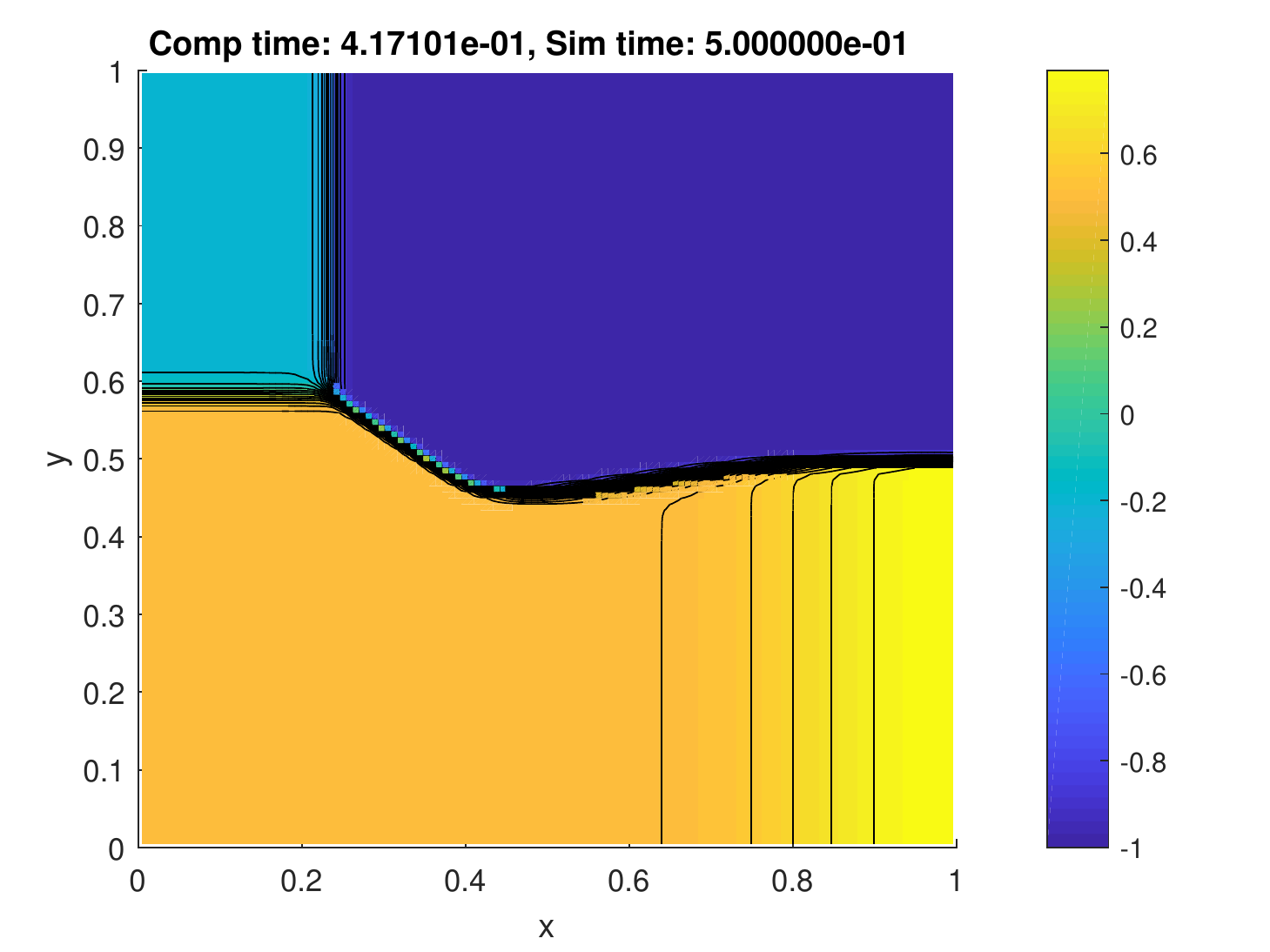}
\includegraphics[width=0.2\textwidth]{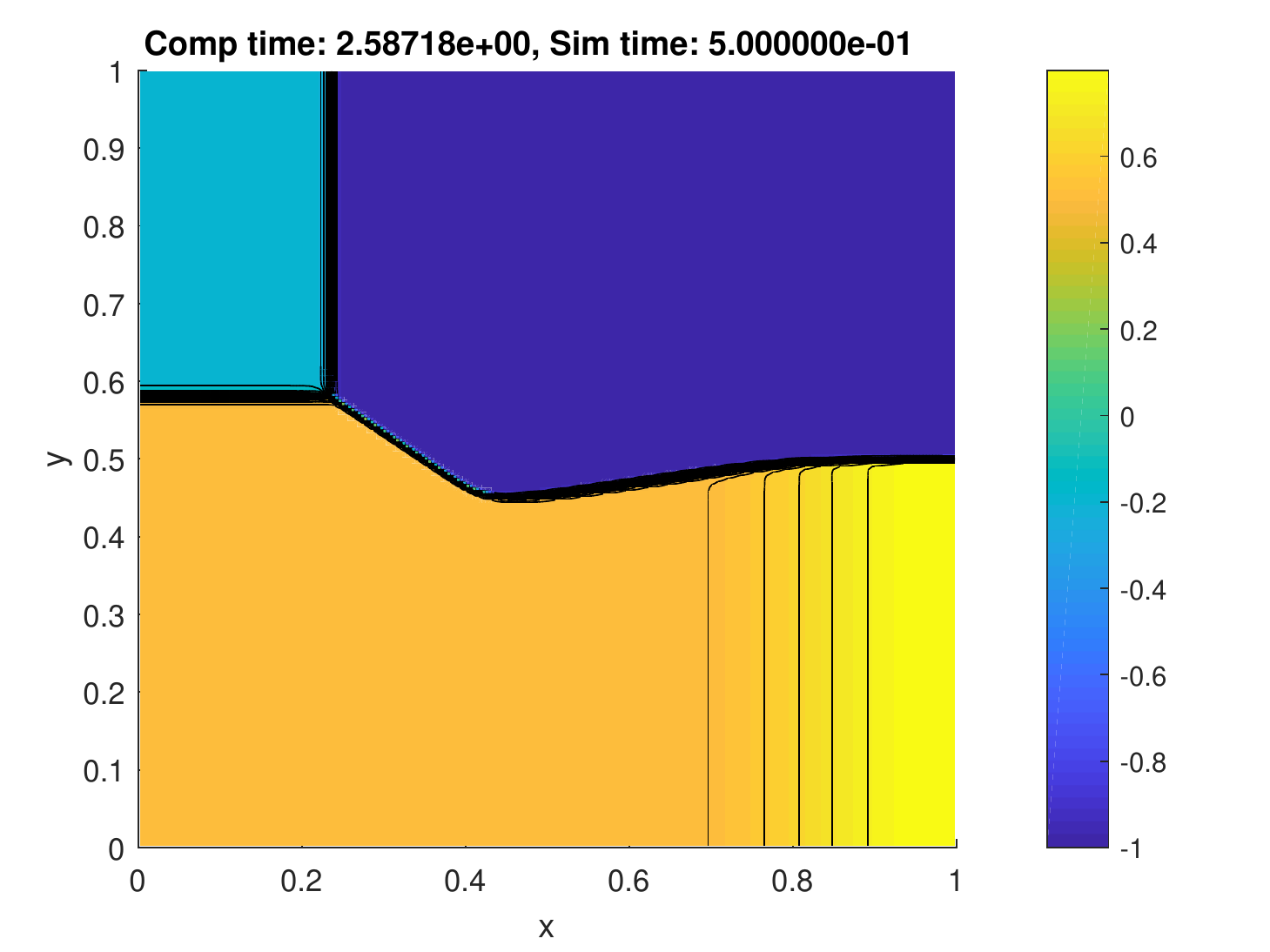}
\includegraphics[width=0.2\textwidth]{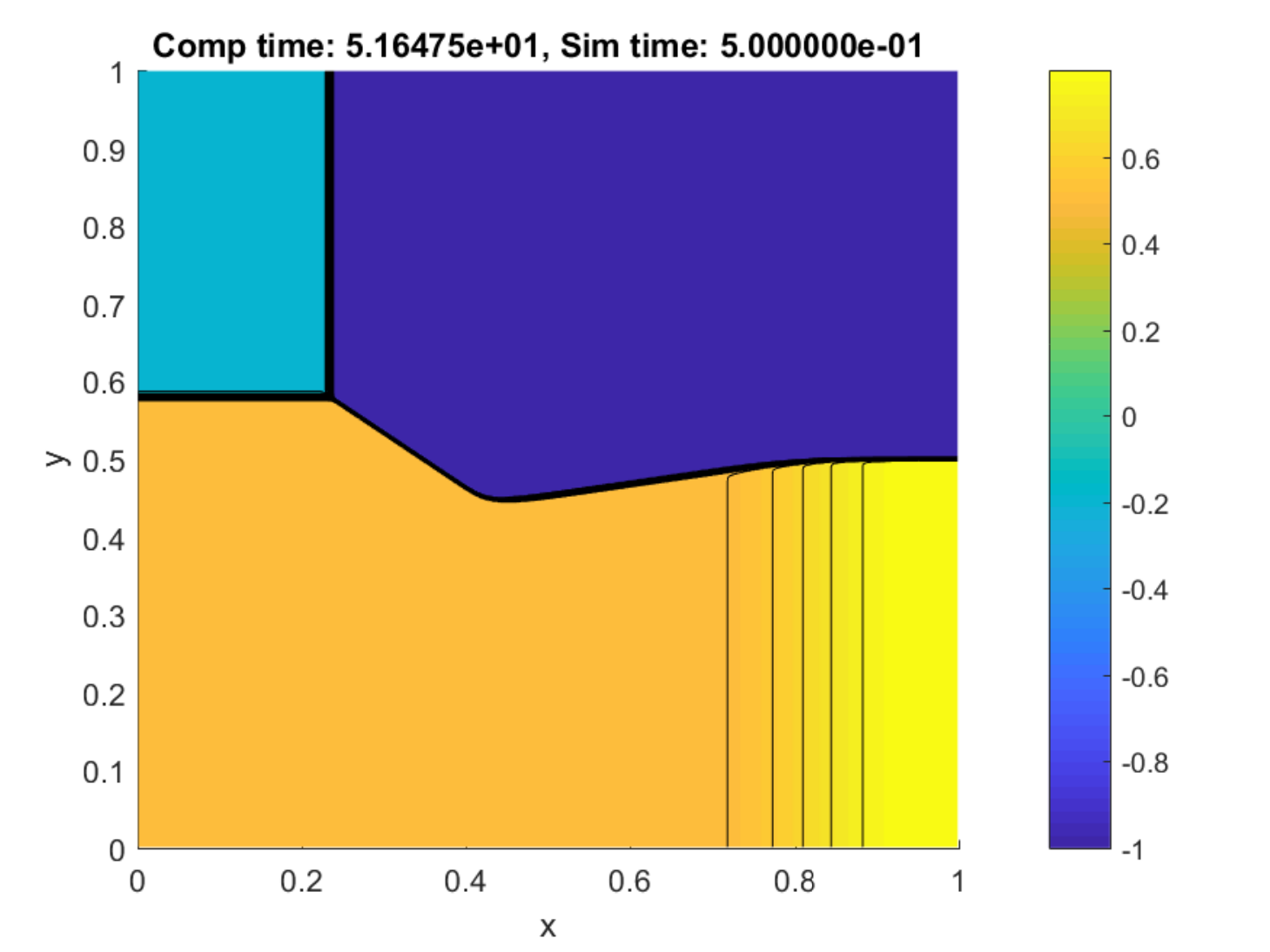}
\includegraphics[width=0.2\textwidth]{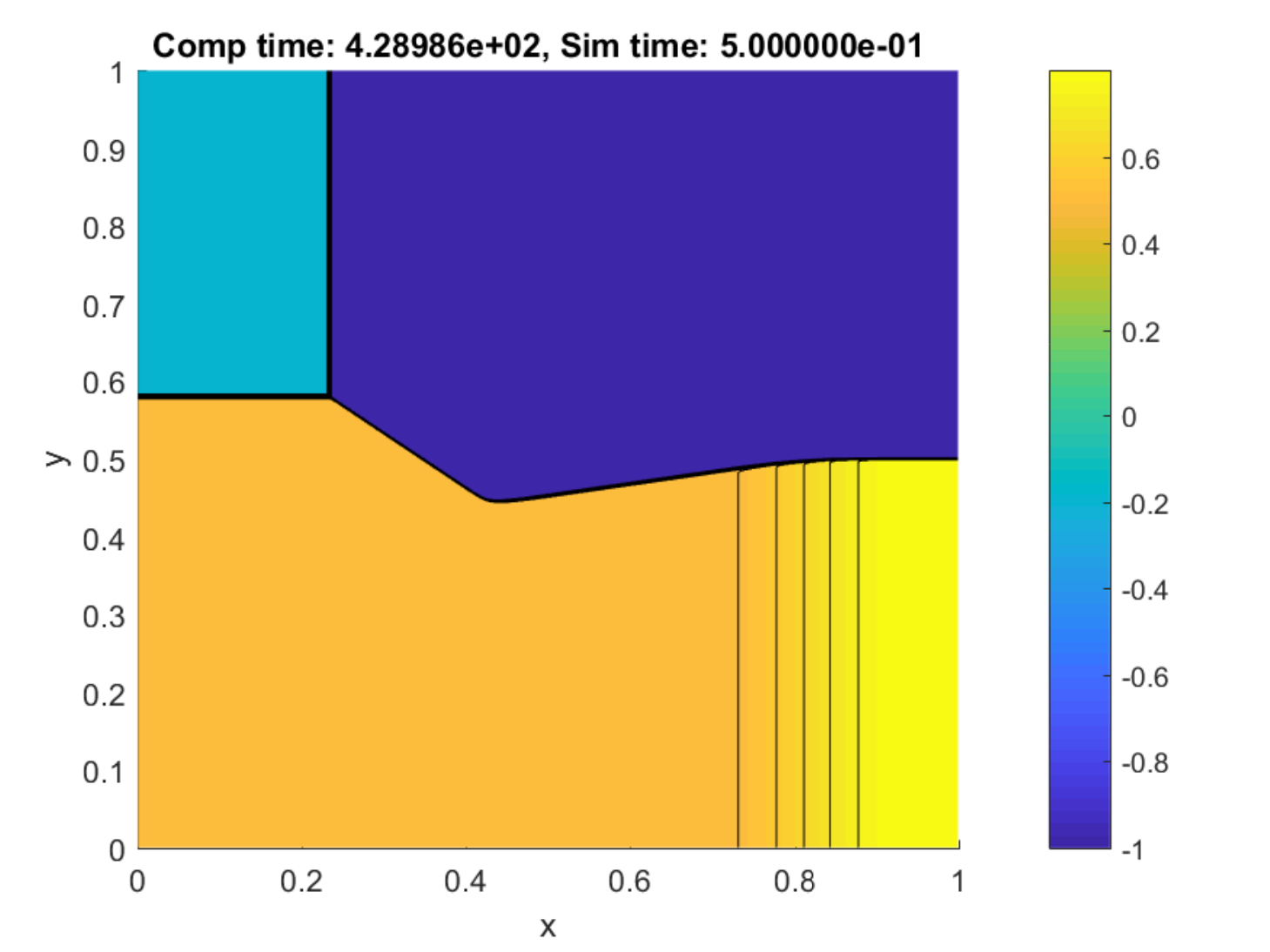}
\caption{
In top Figure are shown 3D numerical 
approximations with 128 x 128 to 
1024 x 1024 size grid respectively, in down 
respectively projections over $xy$-plane.
}
\label{BuObl}
\end{figure}

We consider the 2D initial value problem for the 
inviscid Burgers's equation 
as proposed in \cite{GSJET98},
\begin{equation}
   \displaystyle\frac{\partial u}{\partial t} 
   + \frac{\partial }{\partial x}\left(\displaystyle\frac{u^2}{2}\right)
   + \frac{\partial }{\partial y}\left(\frac{u^2}{2}\right)= 0,
   \label{convex1}
 \end{equation}
\noindent where $(x,y,t) \in \left[0,1\right] \times \left[0,1\right] \times \left[0,0.5\right]$, 
 and with the oblique Riemann data,
  \bigskip
  \begin{equation}
    u(x,y,0) = \begin{cases}
                -1.0, & x > 0.5, \,\, y>0.5, \\
                -0.2, & x < 0.5, \,\, y>0.5, \\
                 0.5, & x < 0.5, \,\, y<0.5, \\
                 0.8, & x > 0.5, \,\, y<0.5, \\
                \end{cases}
     \label{convex2CI}
  \end{equation}
in conjuntion with exact boundary condition on the inflow 
portions of $\partial \Omega$. The correct entropic numerical 
solution is advanced from $t=0$ to $t = \frac{1}{12}$ as 
in \cite{ICBP08} (see Figure \ref{BuObl}).
\begin{figure}[ht!]
\centering
\includegraphics[width=0.36\textwidth]{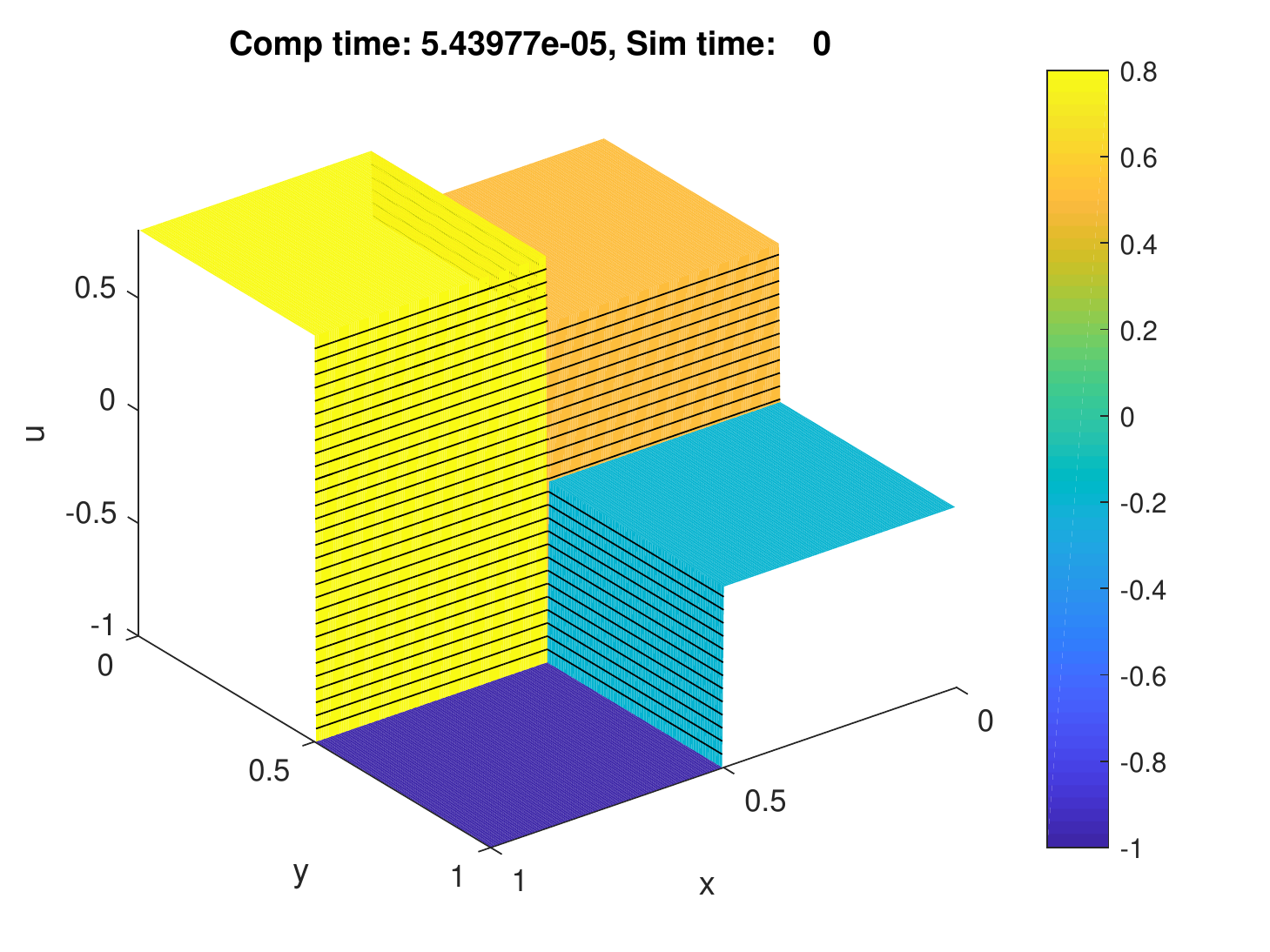}
\includegraphics[width=0.36\textwidth]{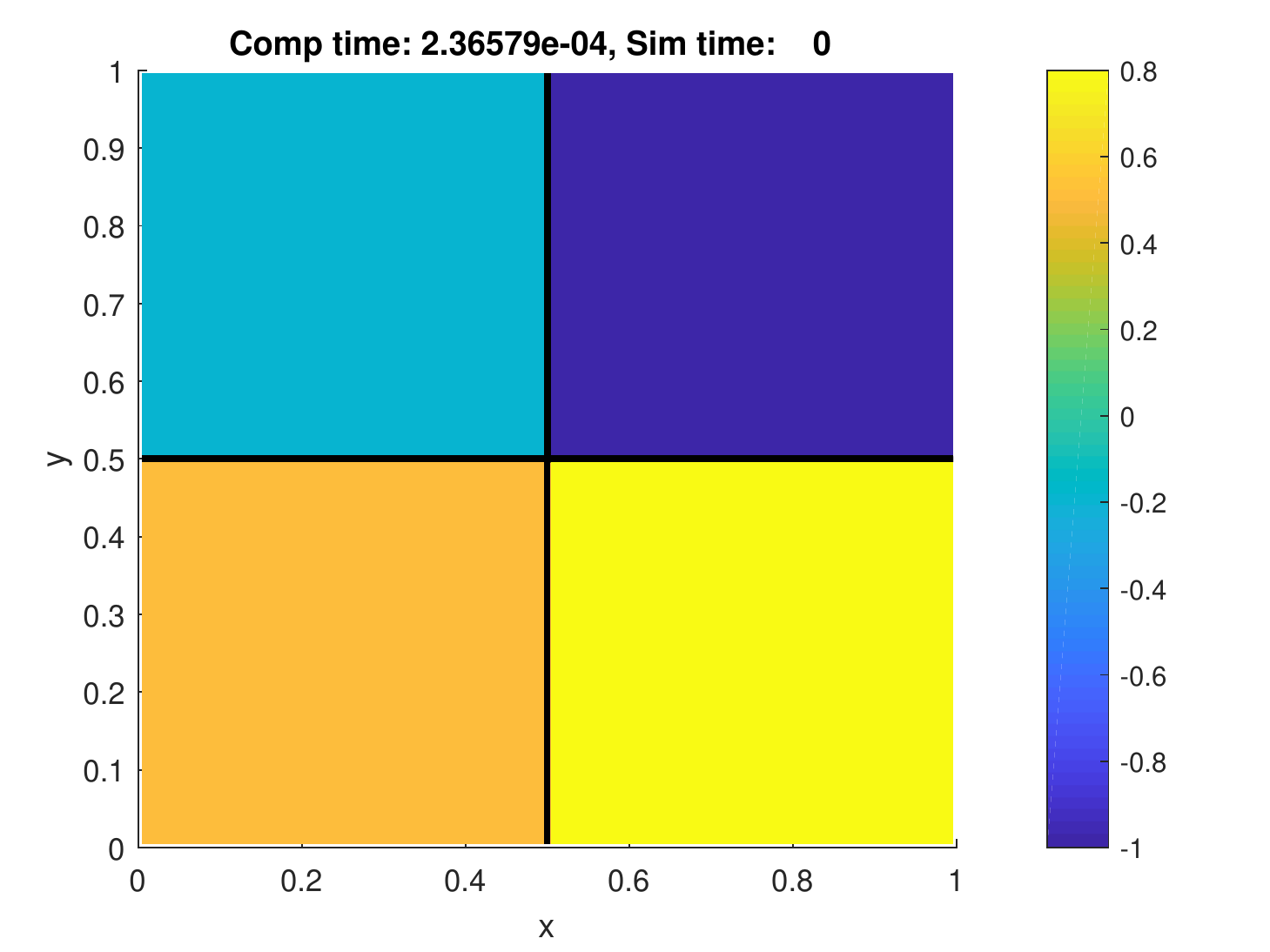}
\caption{Initial Condition}
\end{figure}

\subsubsection{A Buckley-Leverett's problem with gravity}

We consider the reservoir flow model for two-phase water-oil 
immiscible incompressible fluid with gravity \cite{ICBP08},
\vspace{-4mm}
\begin{equation}
   - \nabla \cdot \left[{\bf K} \lambda_{tot}(S_w) \nabla p  \right] = q_{tot}, 
   \label{FlowtwoinmiscPre}
\end{equation}
\begin{equation}
    \displaystyle\frac{\partial \phi \, S_w}{\partial t}
    + \frac{\partial (u_{tot}f(S_w))}{\partial x}
   + \frac{\partial (v_{tot}g(S_w))}{\partial y}= q_{w},
 \label{Flowtwoinmisc}
\end{equation}
where ${\bf K}$ is the absolute permeability tensor, 
$\lambda_{tot}$ is the total mobility, $p$ is the 
thermodynamic pressure, $\phi$ is the porosity, 
$S_w \in \left[ 0,1 \right]$, $S_w$ is the water 
saturation, and ${\bf u}_{tot}=(u_{tot},v_{tot})$ is 
the total velocity (i.e., ${\bf u}_{tot}=u_w+u_o$). 
The pressure equation (\ref{FlowtwoinmiscPre}) as 
written is elliptic in the absence of compressibility
and reads
$- \nabla \cdot \left[{\bf K} \lambda_{tot}(S_w) \nabla p  \right] = 0$. 
Because the total mobility depends of saturation, the 
pressure yields fields changes as the displacement
evolves, this is just a statement of Darcy's law 
combined with the conservation of mass. Once the 
pressure is computed from (\ref{FlowtwoinmiscPre}), 
the total velocity is given by Darcy's law: 
$u_{tot} = -{\bf K}\lambda_{tot}(S_w) \nabla p$. The 
equation (\ref{Flowtwoinmisc}) is referred to as the 
saturation equation. Finally, in the absence of 
gravity and capillarity effects the $x$- and 
$y$-direction flux functions $f(S_w)$ and $g(S_w)$ 
are both just the fractional flow function of water, 
i.e., the non-convex Buckley-Leverett flux:
\begin{equation}
 g(S_w) = f(S_w)=\displaystyle
 \frac{S_w^2}{S_w^2 + \frac{\mu_w}{\mu_0}(1-C_g(1-S_w)^2)},
 \label{Flow}
\end{equation}
here $\mu_w$ and $\mu_o$ are the water and oil phase 
viscosities, respectively. For simplicity, in the 
simulations discussed here, we have chosen the following 
values of the parameters: K is the $2 \time 2$ identity 
matrix, $\lambda_{tot}(S_w) = 1$, $\phi=1$, 
$q_{tot} = q_{w} = 0$. Generally, the complet solution 
of the system (\ref{Flowtwoinmisc}) and 
(\ref{FlowtwoinmiscPre}) is obtained by the implicit 
method to the pressure equation (\ref{FlowtwoinmiscPre}) 
and the explicit method by the hyperbolic equation 
such approach is called an Implicit Pressure Explicit 
saturation (IMPES) sequential solver.

In this example, we consider the Buckley-Leverett problem 
with gravity proposed in \cite{ICBP08} under the above
assumptions with ${\bf u}_{tot}=(1,1)$. The equations 
are \mbox{sign}ificantly more challenging when gravitational 
effects are included in the saturation equation, resulting 
in different (non-convex) flux functions in the $x$- and 
$y$-directions. In this case, $f(\cdot)$ once again the
Buckley-Leverett flux (Flow), but for the flux in the
$y$-direction we have, 
\begin{equation}
 g(S_w) = f(S_w)(1-C_g(1-S_w)^2).
 \label{FlowG}
\end{equation}

\begin{figure}[ht!]     
 \begin{minipage}[b]{0.5\linewidth}
 \begin{equation}
   \displaystyle\frac{\partial S_w}{\partial t} 
   + \frac{\partial (f(S_w))}{\partial x}
   + \frac{\partial (g(S_w))}{\partial y}= 0,
   \label{nonconvexWG}
 \end{equation}
  with $(x,y,t) \in \left[-1.5,1.5\right] \times \left[-1.5,1.5\right] \times \left[0,0.5\right]$, 
 and initial condition,
    \begin{equation}
    u(x,y,0) = \begin{cases}
                1 ,  & x^2 + y^2 < 0.5, \\
                0 , & \hbox{ otherwise}. \,\,
                \end{cases}
     \label{WithGCI}
    \end{equation}
   \end{minipage}
   \begin{minipage}[t]{0.4\linewidth}
     \centering
        \includegraphics[scale=0.5]{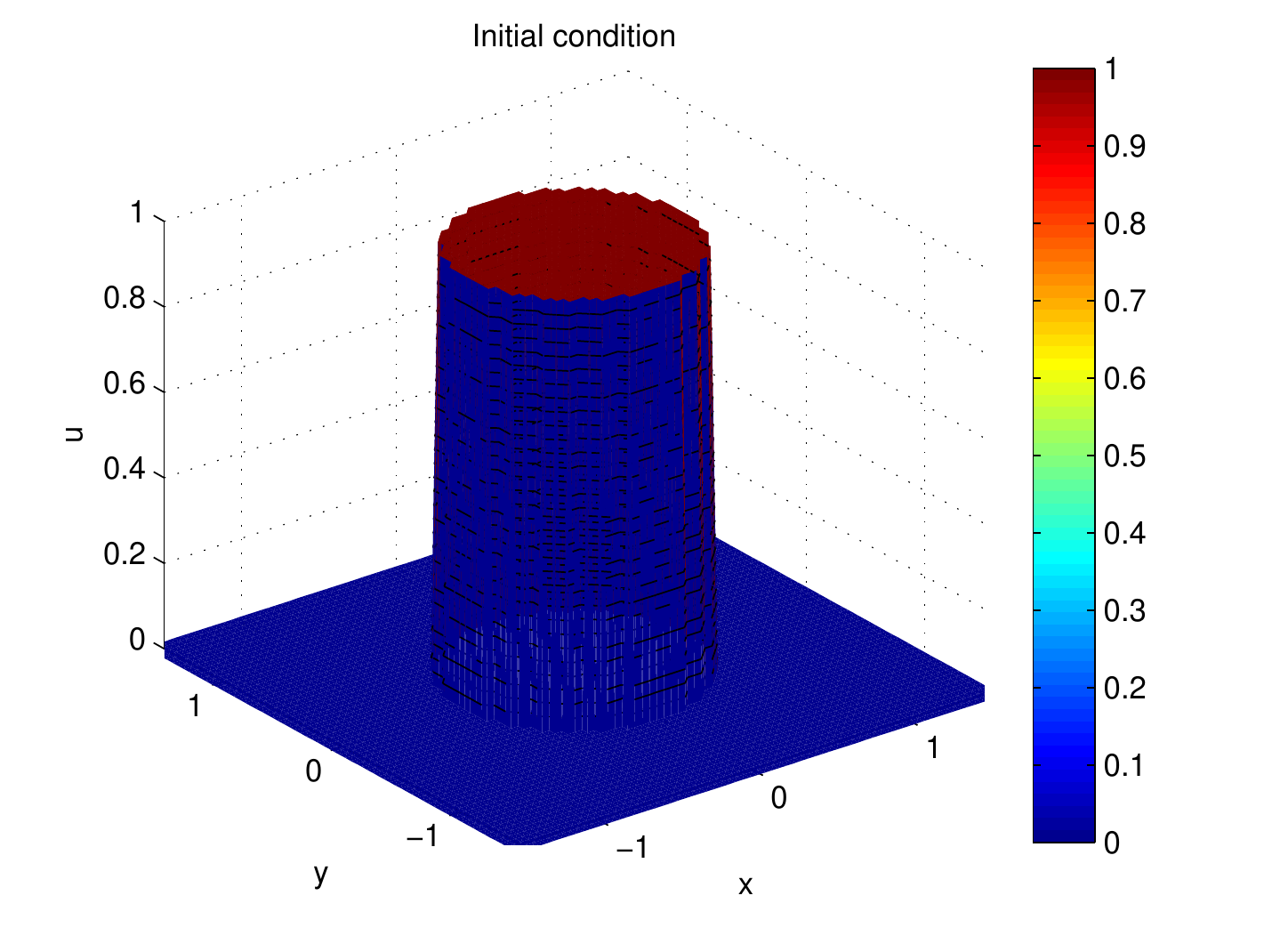}
	\end{minipage}
	\caption{Initial condition flow with gravity.}
  \end{figure} 

\begin{figure}[ht!]
\centering
\includegraphics[width=0.2\textwidth]{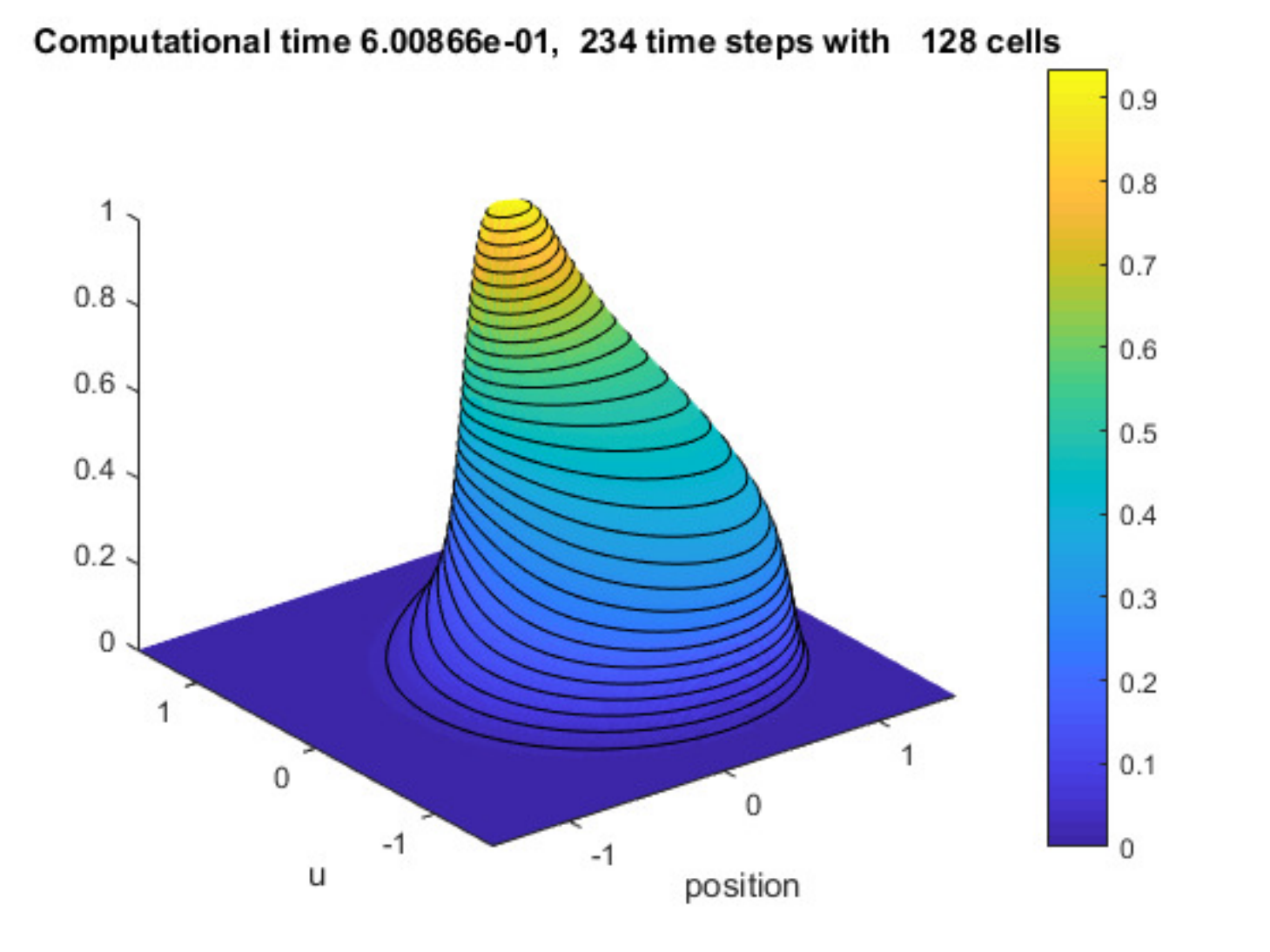}
\includegraphics[width=0.2\textwidth]{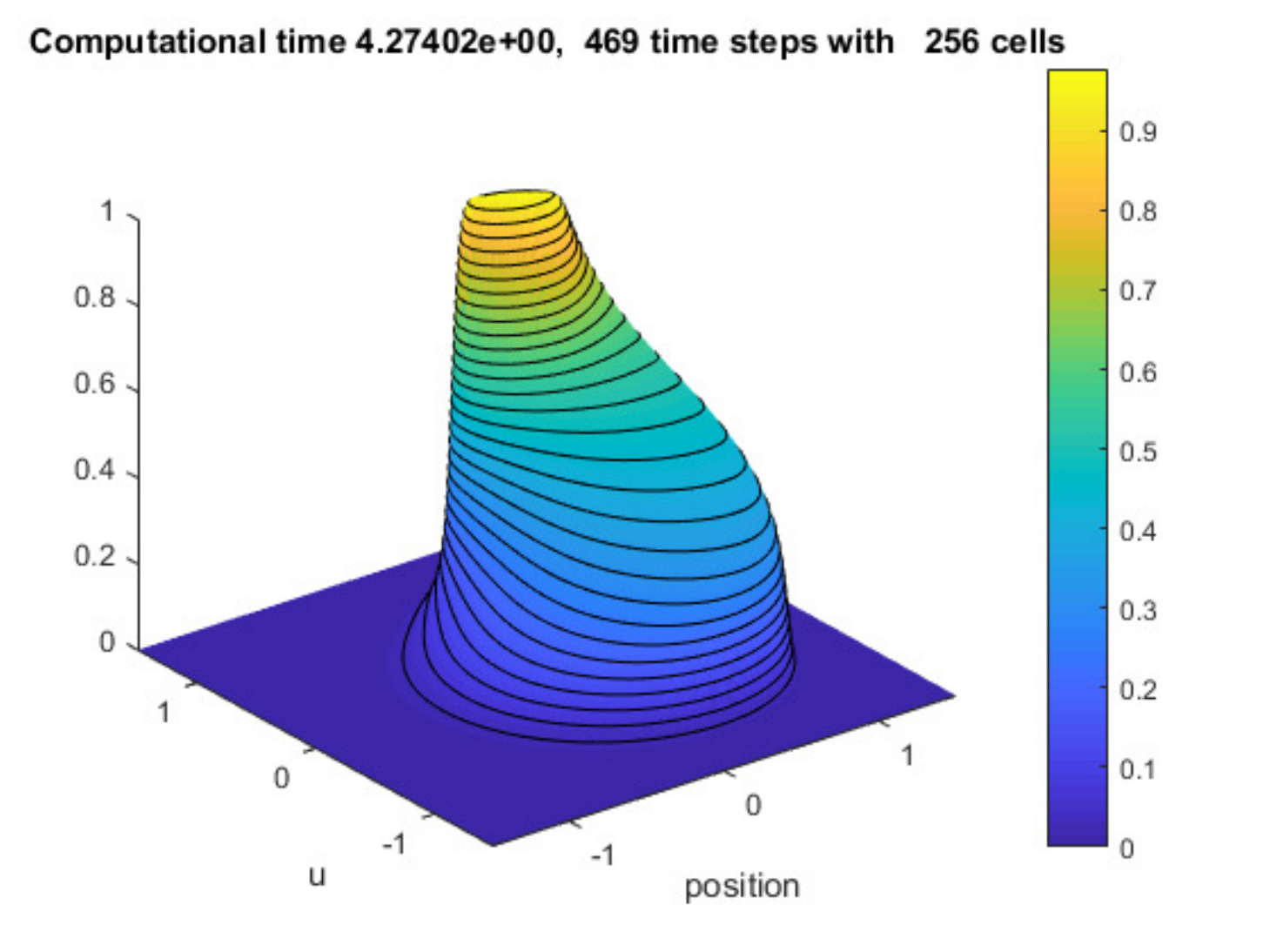}
\includegraphics[width=0.2\textwidth]{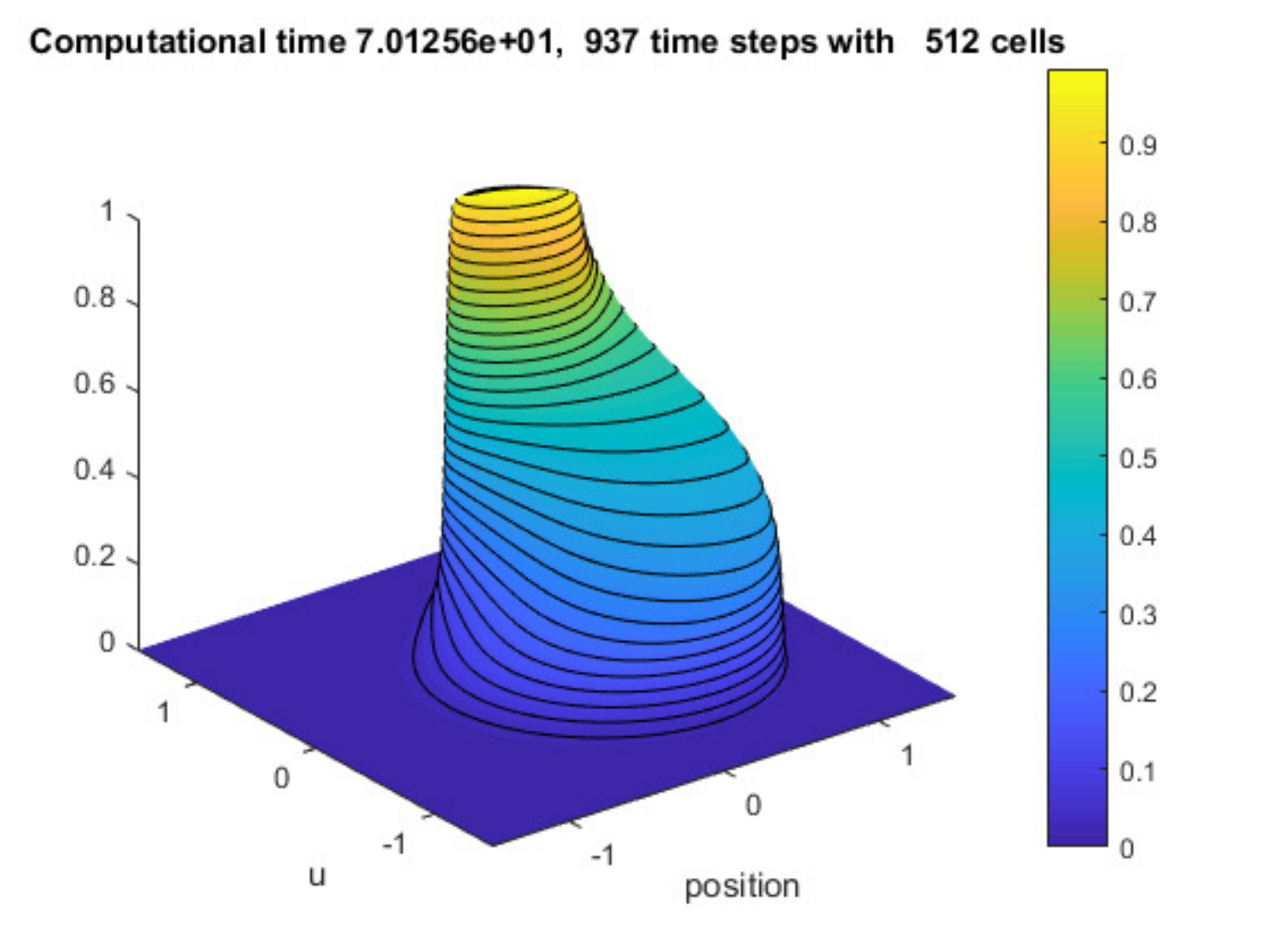}
\includegraphics[width=0.2\textwidth]{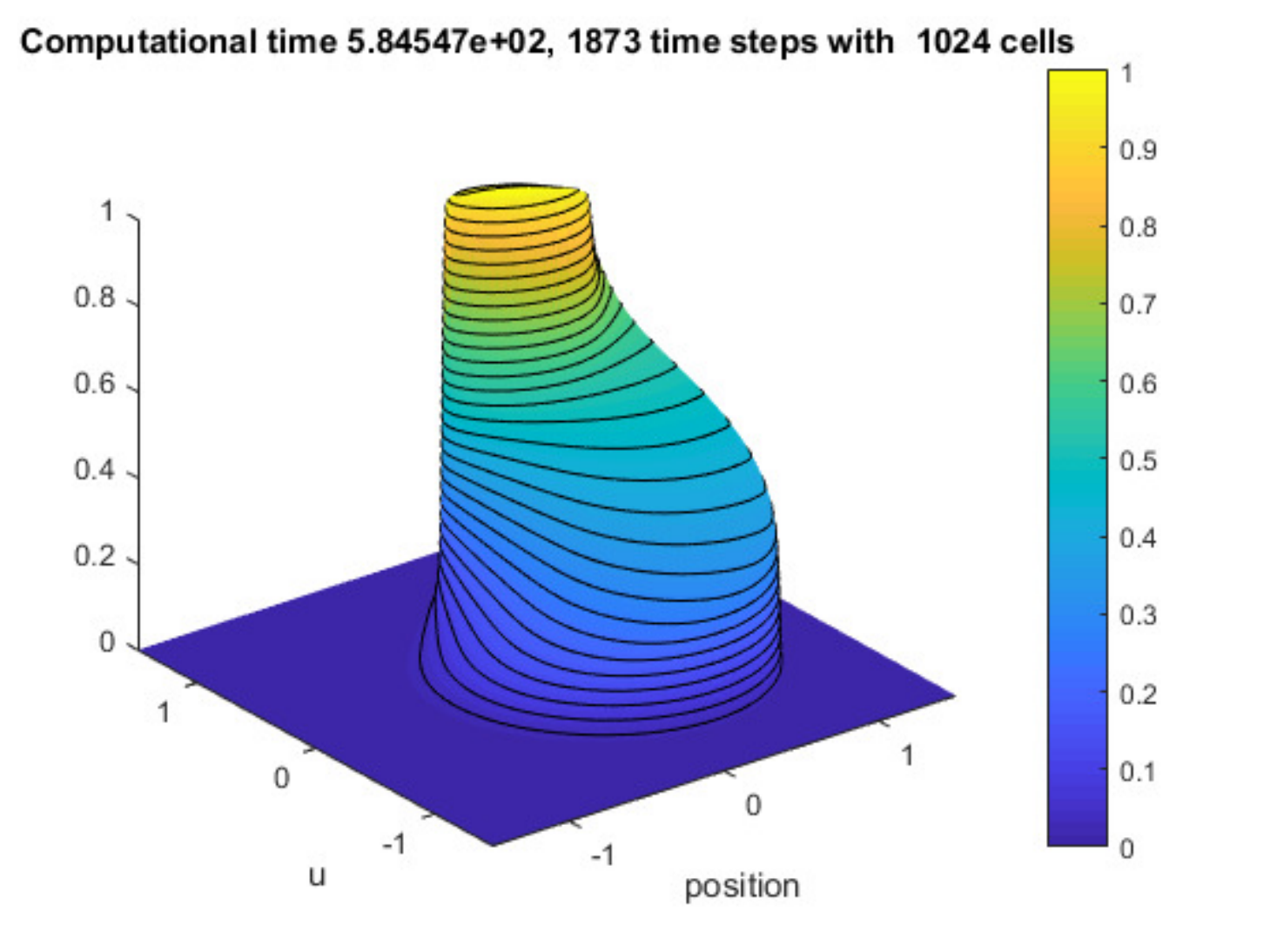} \\
\includegraphics[width=0.2\textwidth]{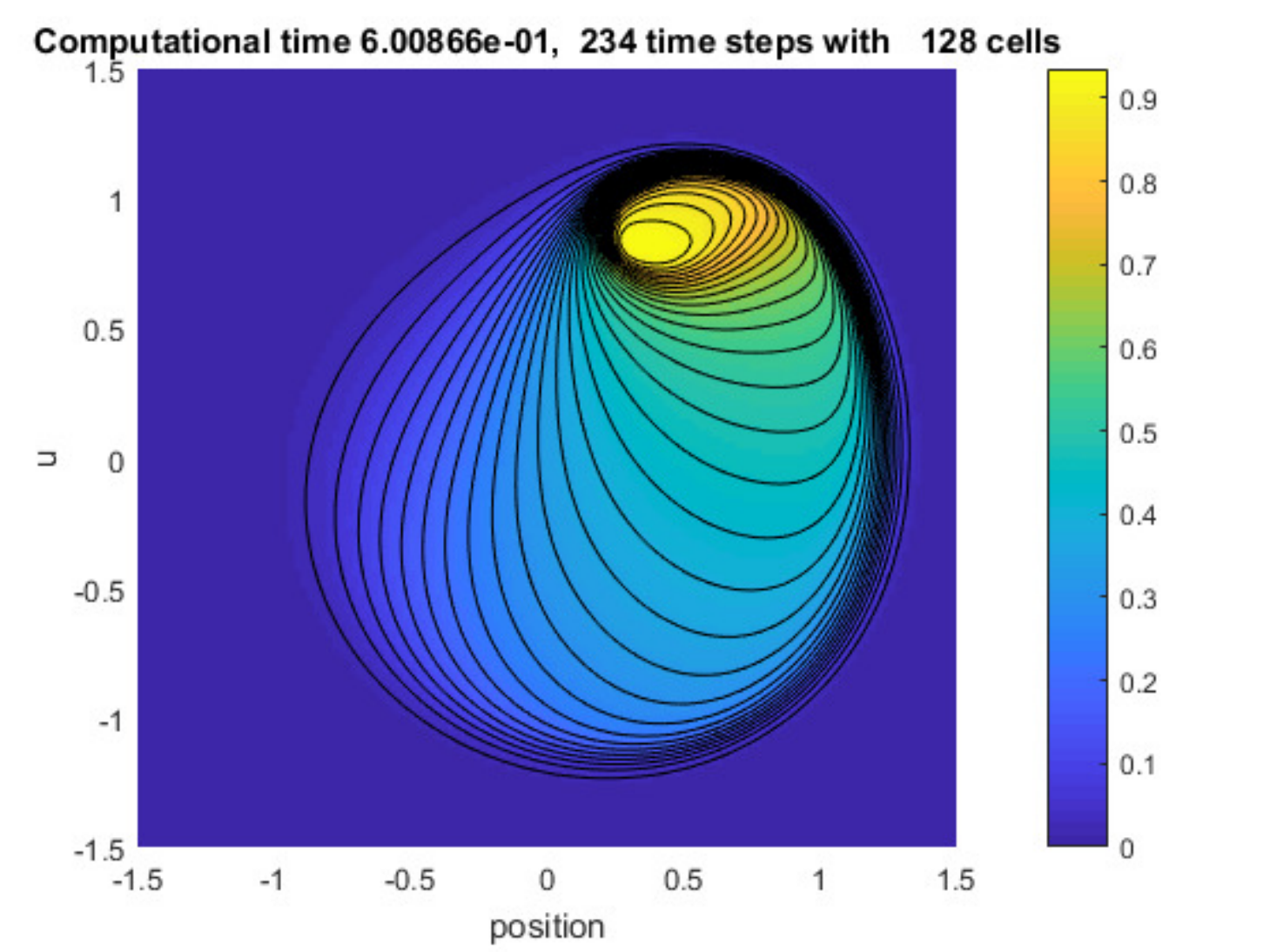}
\includegraphics[width=0.2\textwidth]{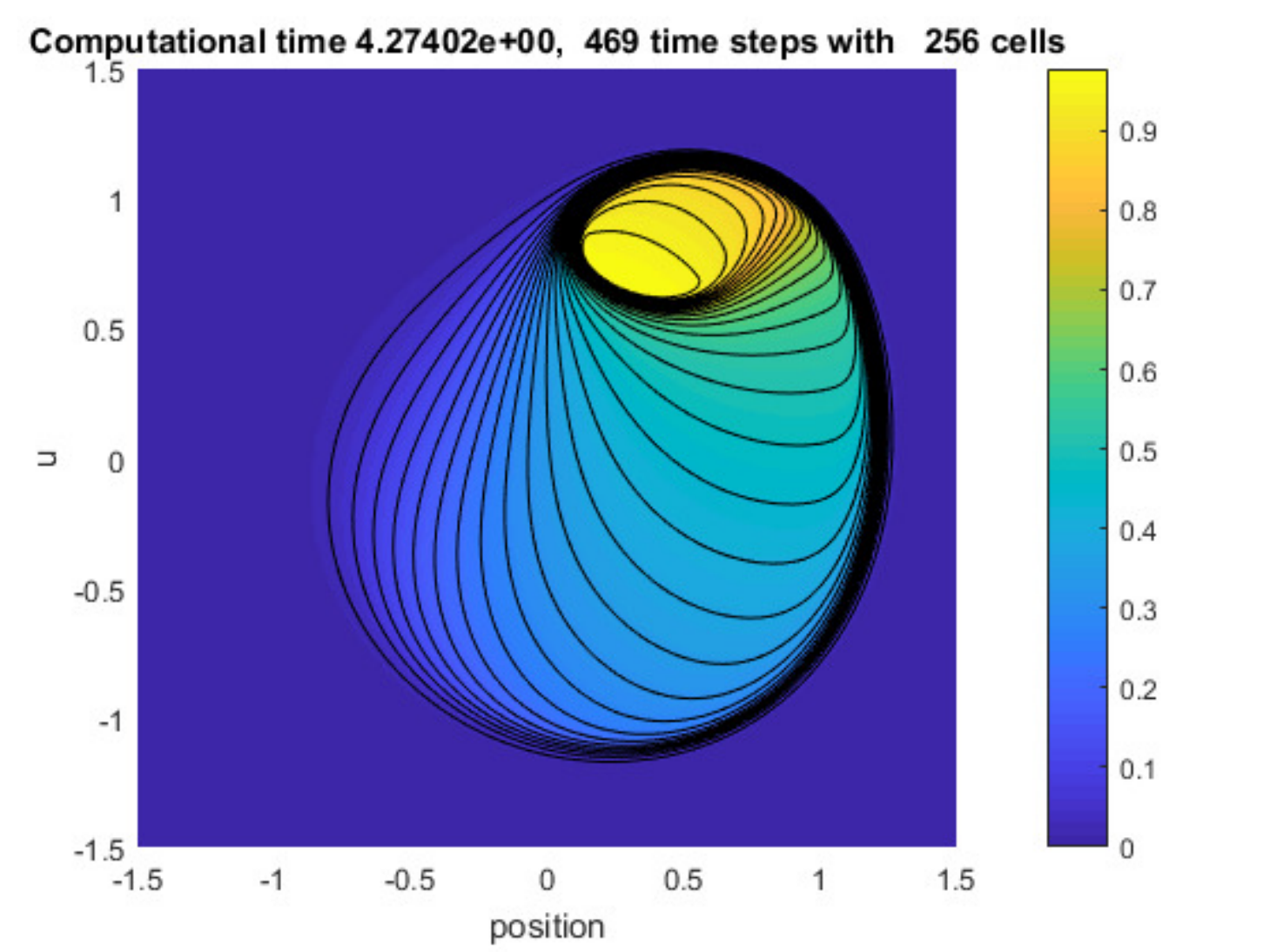}
\includegraphics[width=0.2\textwidth]{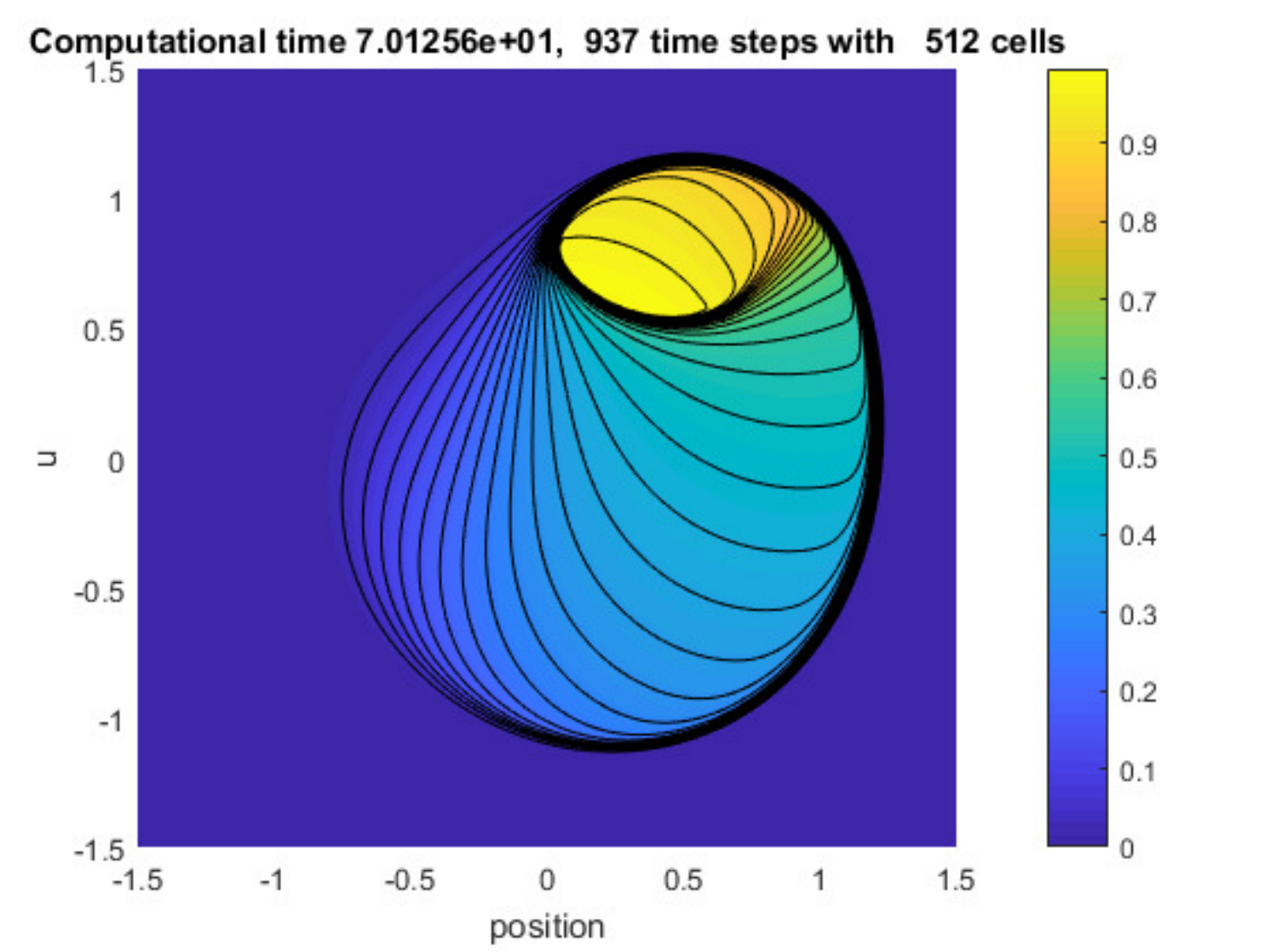}
\includegraphics[width=0.2\textwidth]{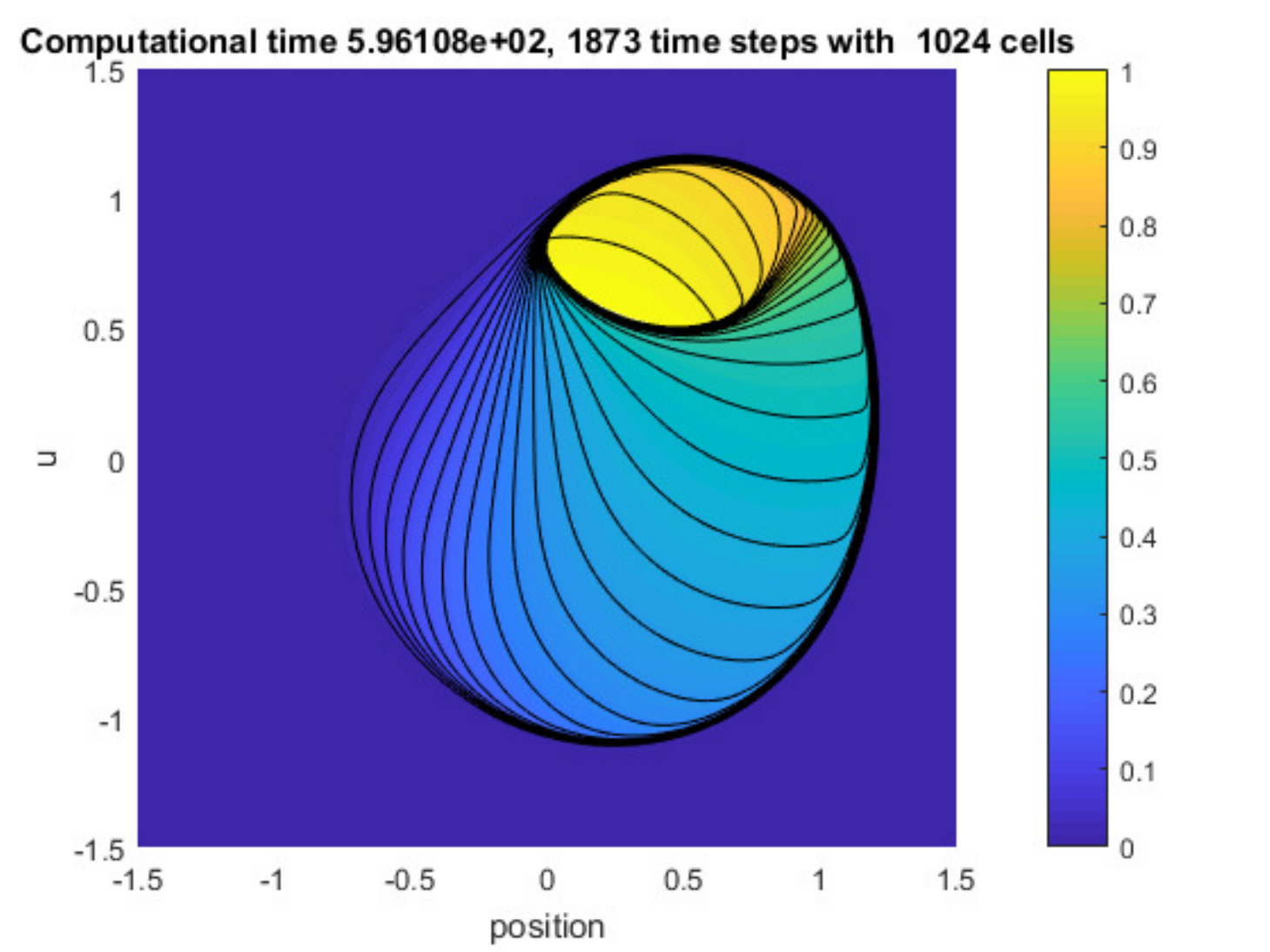}
\caption{
From top to bottom are shown the  Buckley-Leverett's 
numerical solutions computed with the two-dimensional 
scheme with respect to the reservoir waterflooding problem 
with gravity described in (\ref{FlowtwoinmiscPre})-(\ref{WithGCI}).}
\label{BLWithG}
\end{figure}

Finally, we notice that we impose the solid wall (slip) 
boundary condition ${\bf u}_{tot} \cdot n = 0$, everywhere 
on the boundary $\partial_\Omega$, where $n$ is the outward 
unit normal to $\partial \Omega$, upon the system 
(\ref{FlowtwoinmiscPre}) and (\ref{Flowtwoinmisc}). This 
means that there are no inflow boundaries and, hence, no 
boundary conditions on $S_w$. Here we have two situations 
we want to test our Eulerian-Lagrangian scheme: (1) a 
rudimentar test to address the issue of grid orientation 
effects (this anomalous phenomenon is observed when 
computational grid is rotated and substantially different 
numerical solutions are obtained for a same problem) and 
(2) accommodation of no flow boundary condition, exact or 
approximate. 
Finally, we see our numerical solutions shown in Figure 
\ref{BLWithG} for the Buckley-Levertt's problem described 
above (\ref{FlowtwoinmiscPre})-(\ref{WithGCI}) are in a
very good agreement with those computed solutions as in 
\cite{ICBP08}.

\subsubsection{\bf Conservation property verification test 
with scalar analytical solution velocity}

\begin{figure}[ht!]
\centering
\includegraphics[width=0.2\textwidth]{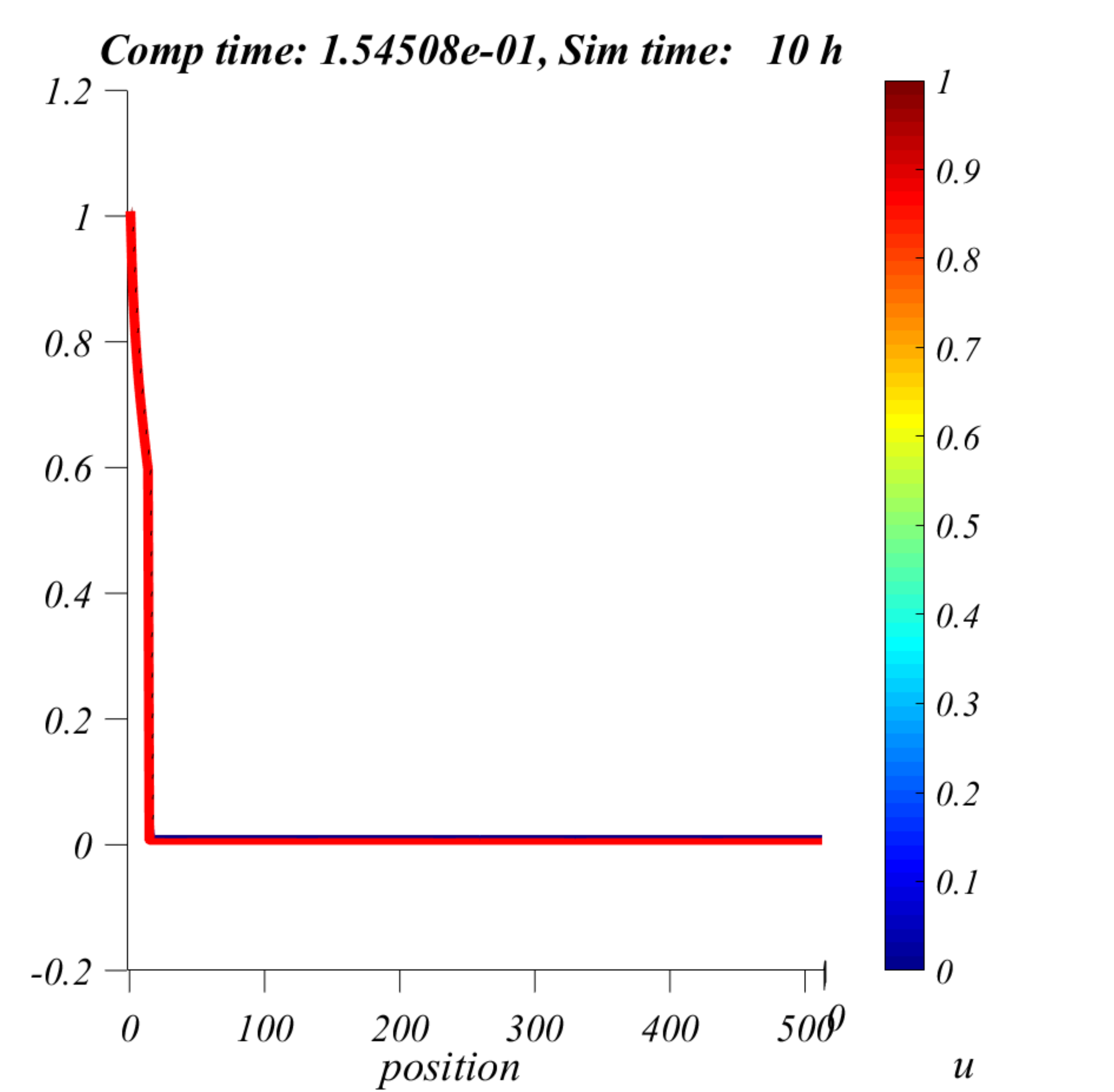}
\includegraphics[width=0.2\textwidth]{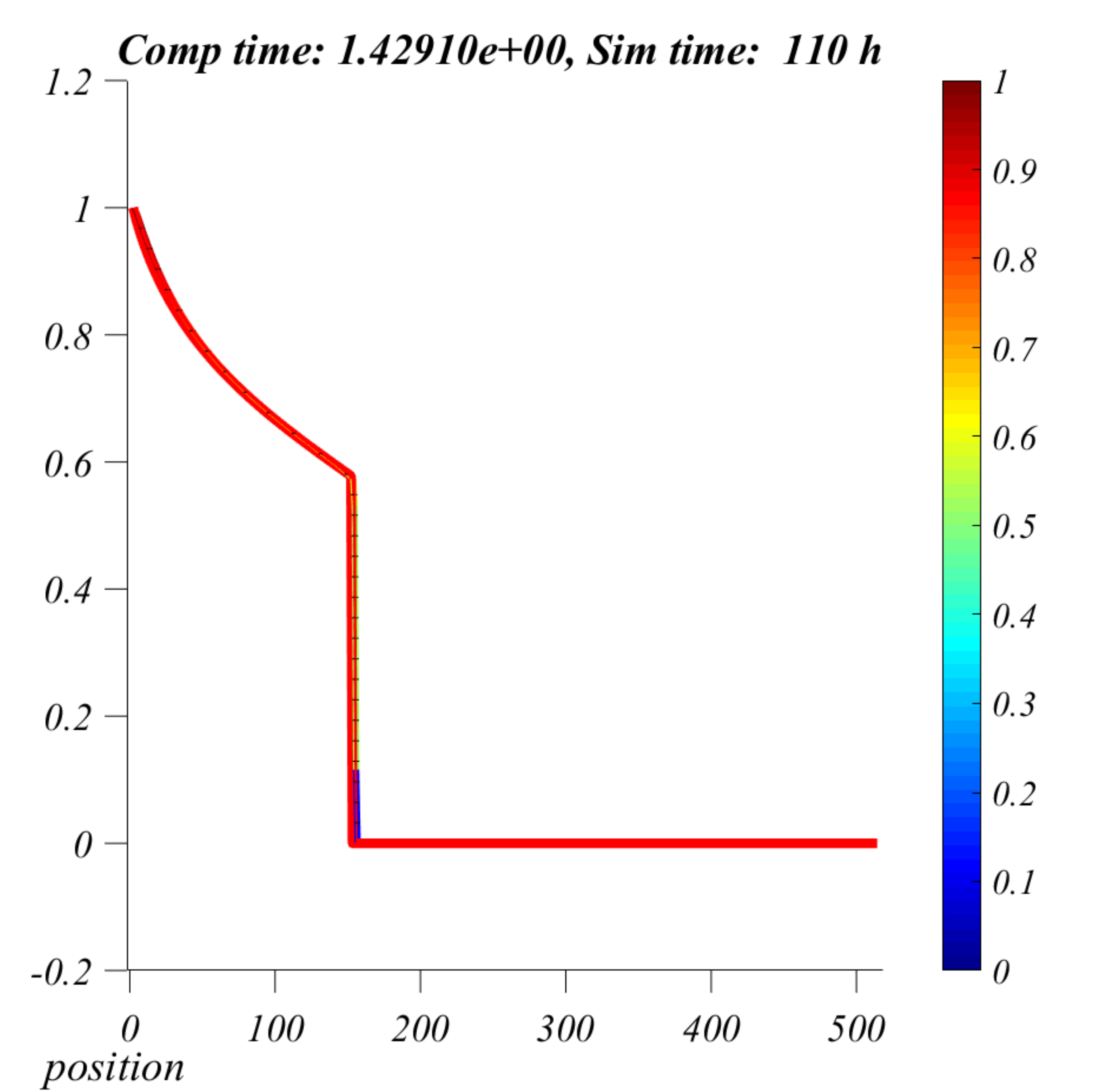}
\includegraphics[width=0.2\textwidth]{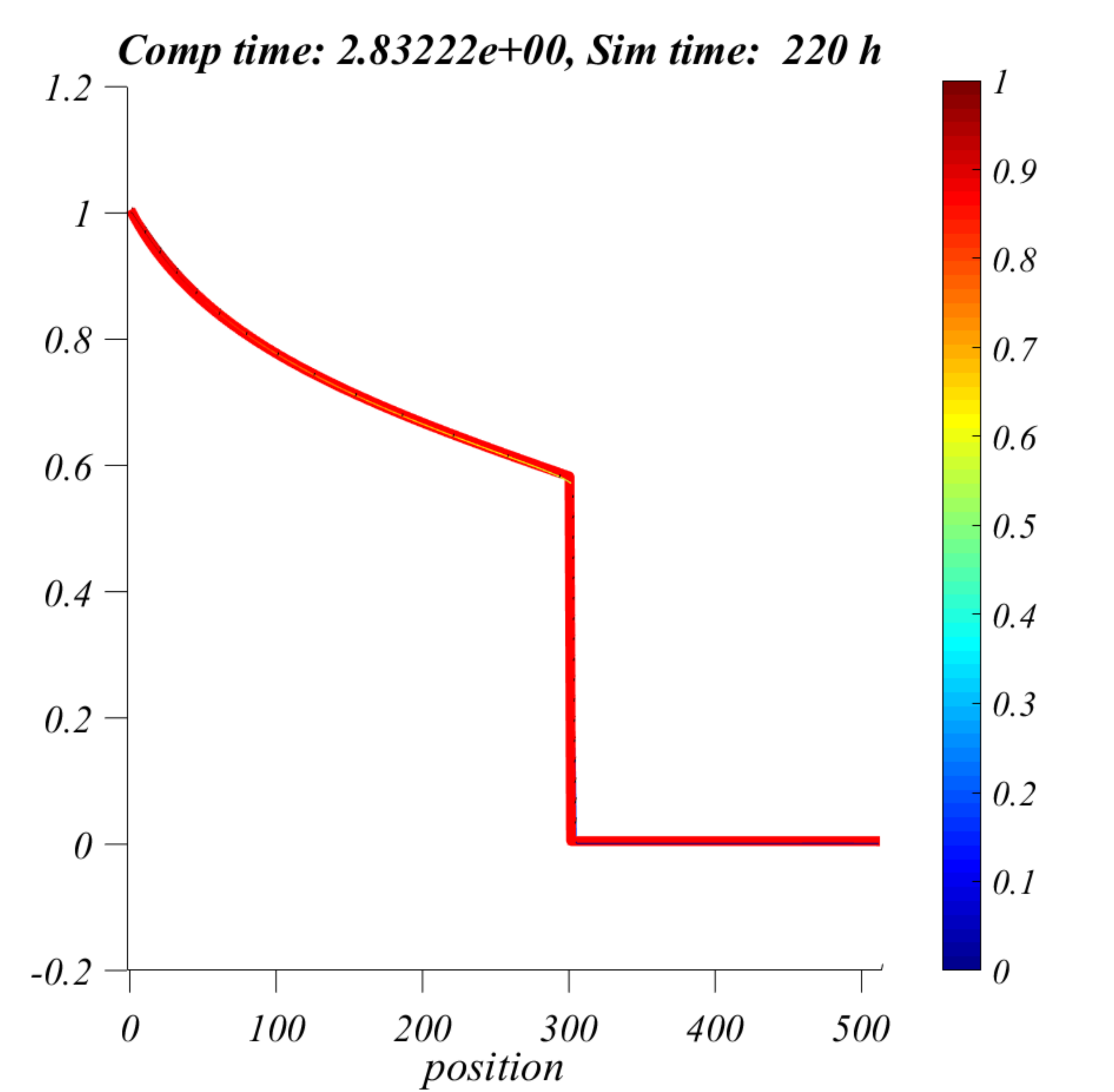}
\includegraphics[width=0.2\textwidth]{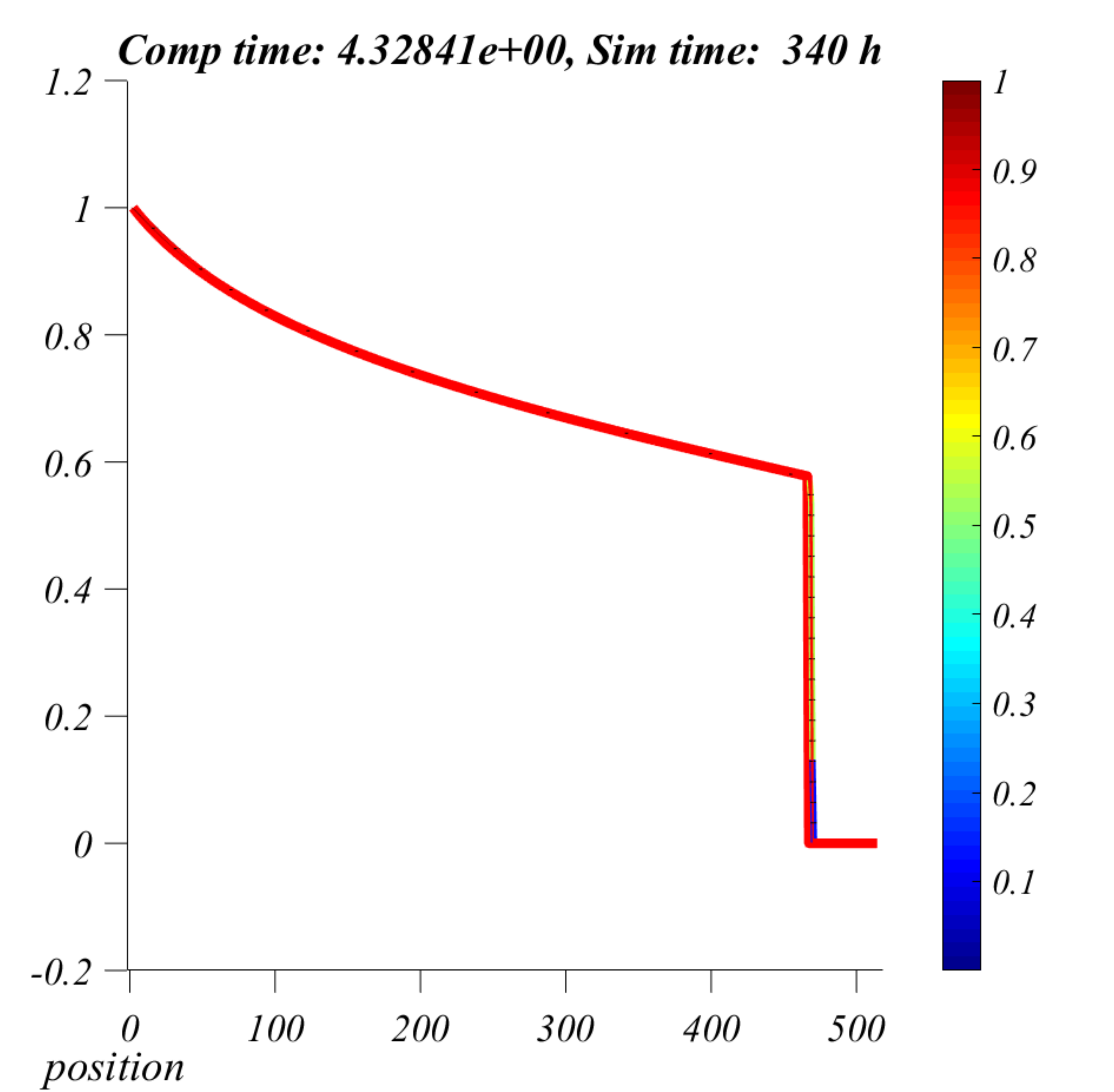} \\
\includegraphics[width=0.2\textwidth]{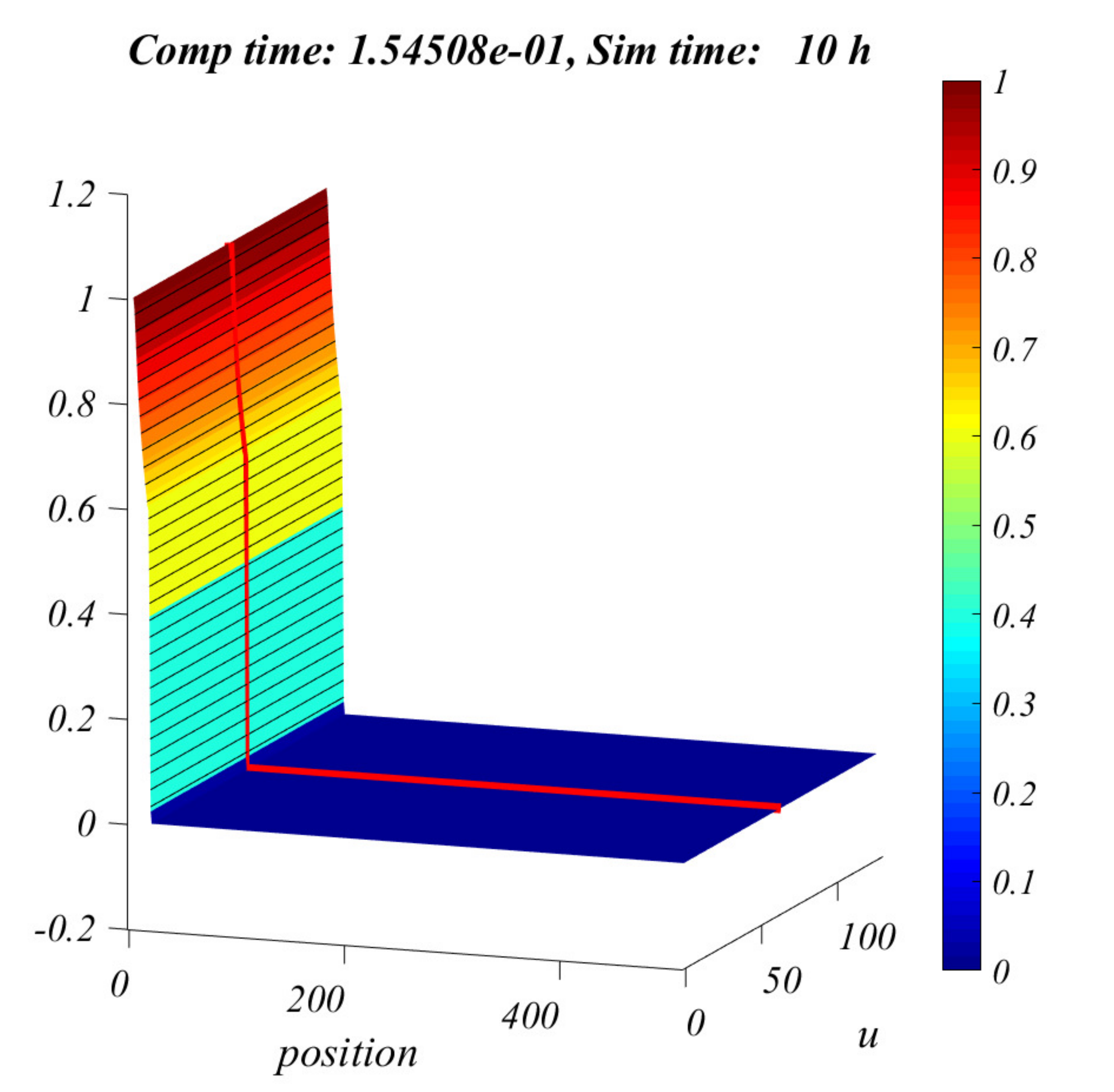}
\includegraphics[width=0.2\textwidth]{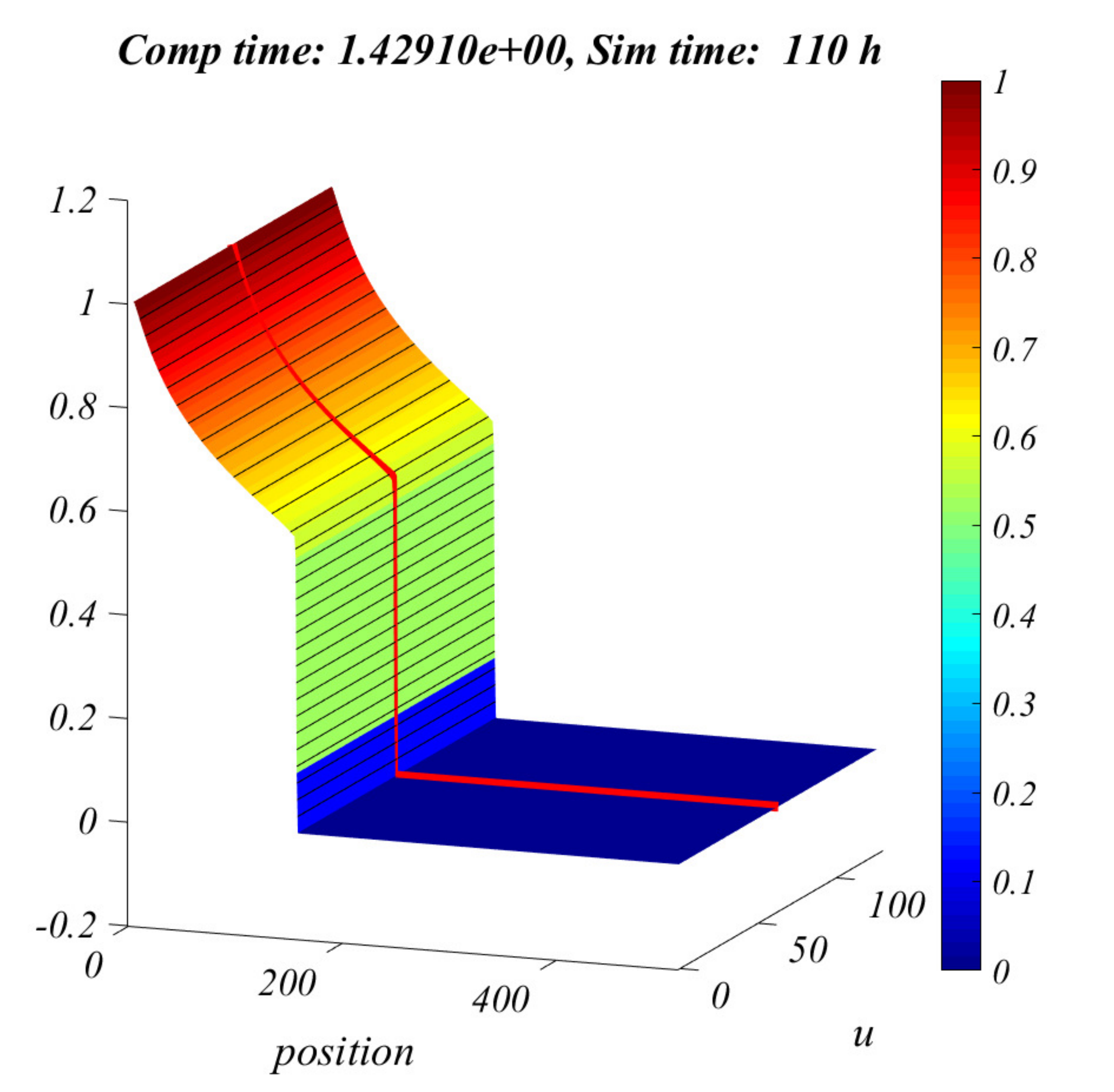}
\includegraphics[width=0.2\textwidth]{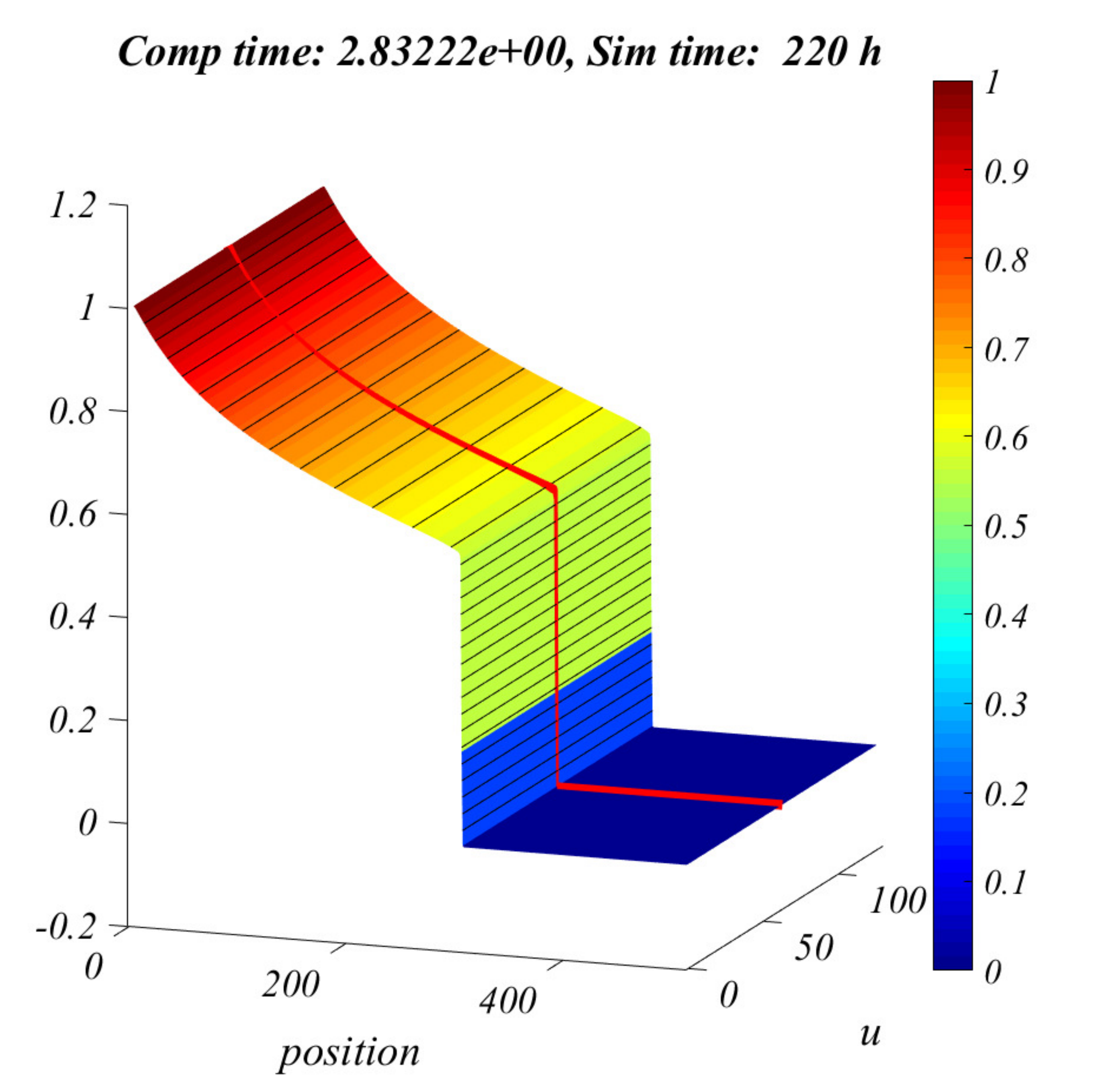}
\includegraphics[width=0.2\textwidth]{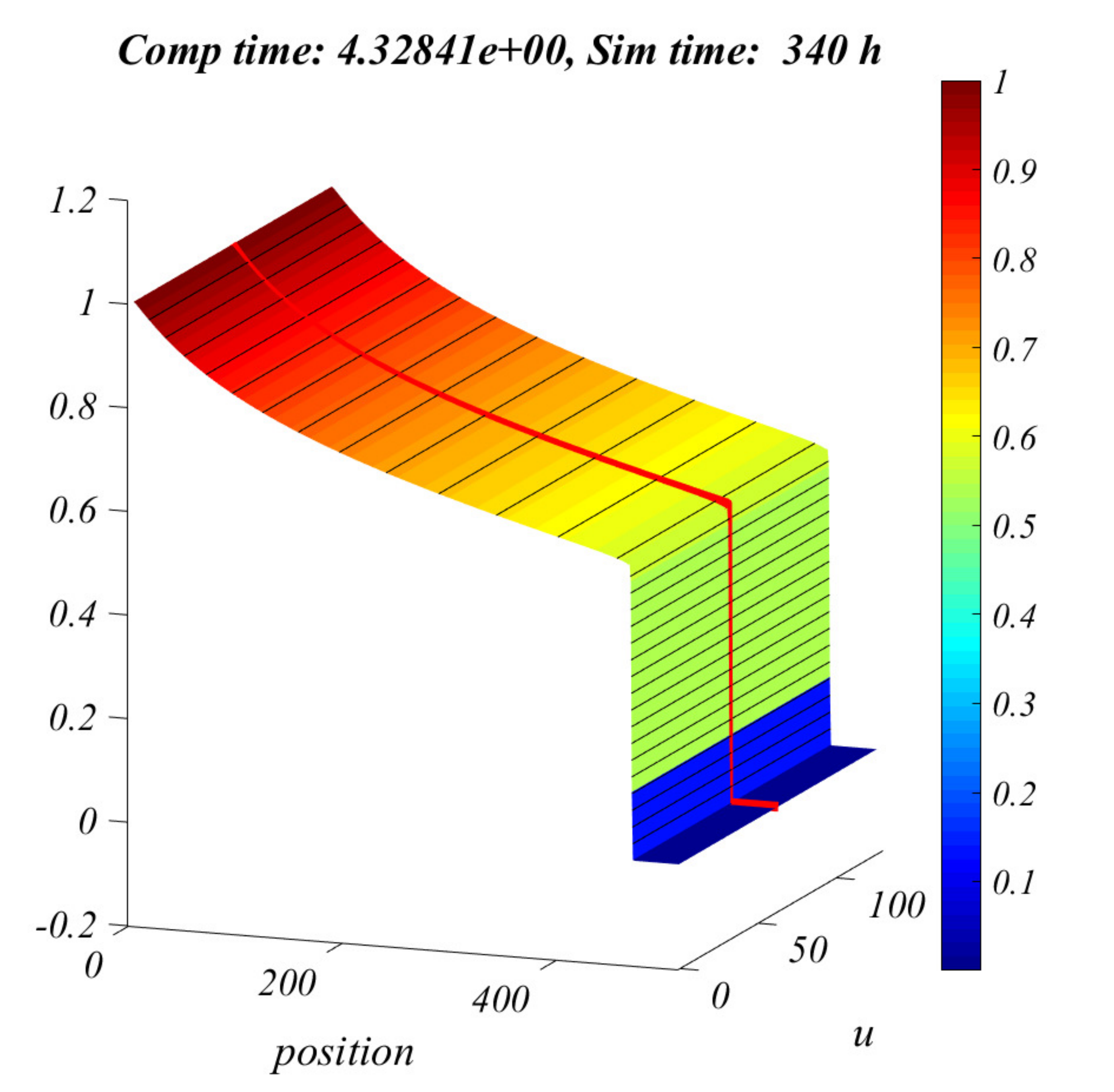}
\caption{Conservation property verification test 
with scalar analytical solution velocity: projections 
over $xu$-plane at times $T$ = 10, 110, 220 and 430, 
respectively (top) and the corresponing 3D views 
(bottom).}
\label{verifytest}
\end{figure}

In Figure \ref{verifytest} is shown 2D numerical 
solutions displayed as time evolves for a 2D symmetric 
Buckley-Leverett problem, in which is possible compare 
with exact analytical solution. Notice in top frames 
are shown the projection over $xu$-plane at times 
$T = 10, 110, 220, 340$ hours, respectively. It is 
clear we might see the correct front velocity of 
the 2D simulation being approximated with our 
proposed scheme when compared with the Buckley-Leverett 
analytical solution (superimposed red lines).

\section{Coupling conservative finite element method 
for Darcy flow problem with a locally conservative 
Lagrangian-Eulerian method for hyperbolic-transport}\label{Coupling}

We combine a novel high-order conservative finite 
element method for Darcy flow problem (Section 
\ref{Ch2}) with locally conservative Lagrangian-Eulerian 
method for hyperbolic-transport (Section \ref{theLE2d}) 
to address conservation properties in multiple scale 
in both, the complexity multiscale heterogeneity 
structures from rock geology appearing in the 
elliptic-pressure-velocity model as well as 
multiscale wave structures resulting from shock 
wave interactions from the hyperbolic-transport model.

We solve saturation and pressure equations 
(\ref{hyper})-(\ref{pres}) in the IMPES sequential 
fashion in a geological domain of $256m \times 64m$ 
considering two representative situations in which 
any one can easly reproduce further latter, namely, 
homogeneous medium and heterogeneous barrier and 
both in the slab geometry as previously described.

\begin{itemize}
    \item {\bf Test 1 :} Let's consider an homogeneous 
medium ($K(x) \equiv 1$) and use the method to solve 
the model (\ref{hyper})-(\ref{pres}) for time from 
$0$ to $220$ using numerical square meshes of norms 
($h_x = h_y = 32, 16, 8, 4, 2,$ and $1$).
    
    \item {\bf Test 2 :} Let's consider an heterogeneous 
barrier medium and use the method to solve the model 
(\ref{hyper})-(\ref{pres}) for time from $0$ to $220$ 
using numerical square meshes of norms ($h_x = h_y = 
32, 16, 8, 4, 2$ and $1$). 
\end{itemize}

\begin{figure}[ht!]
\includegraphics[scale=0.23]{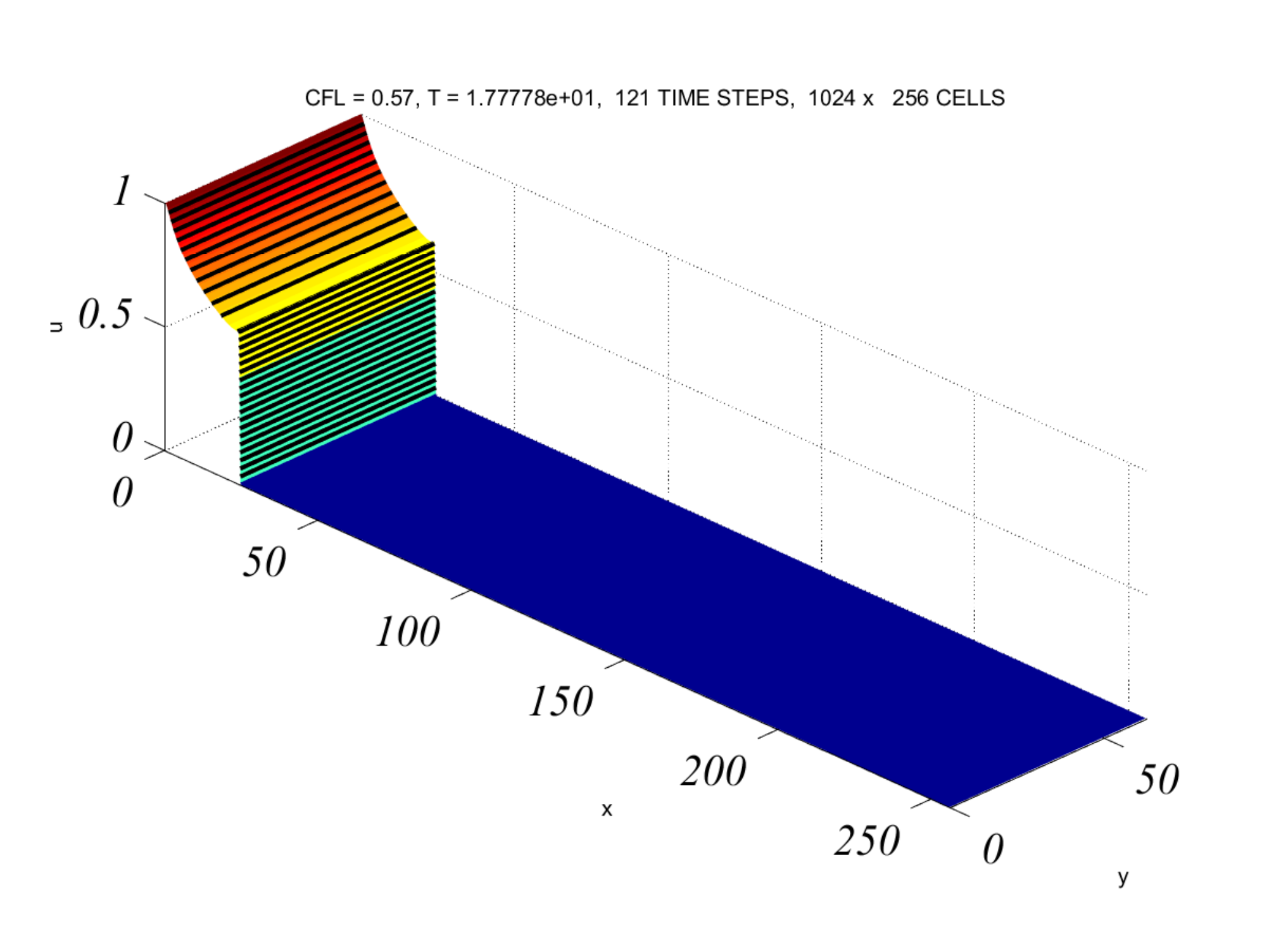}
\includegraphics[scale=0.23]{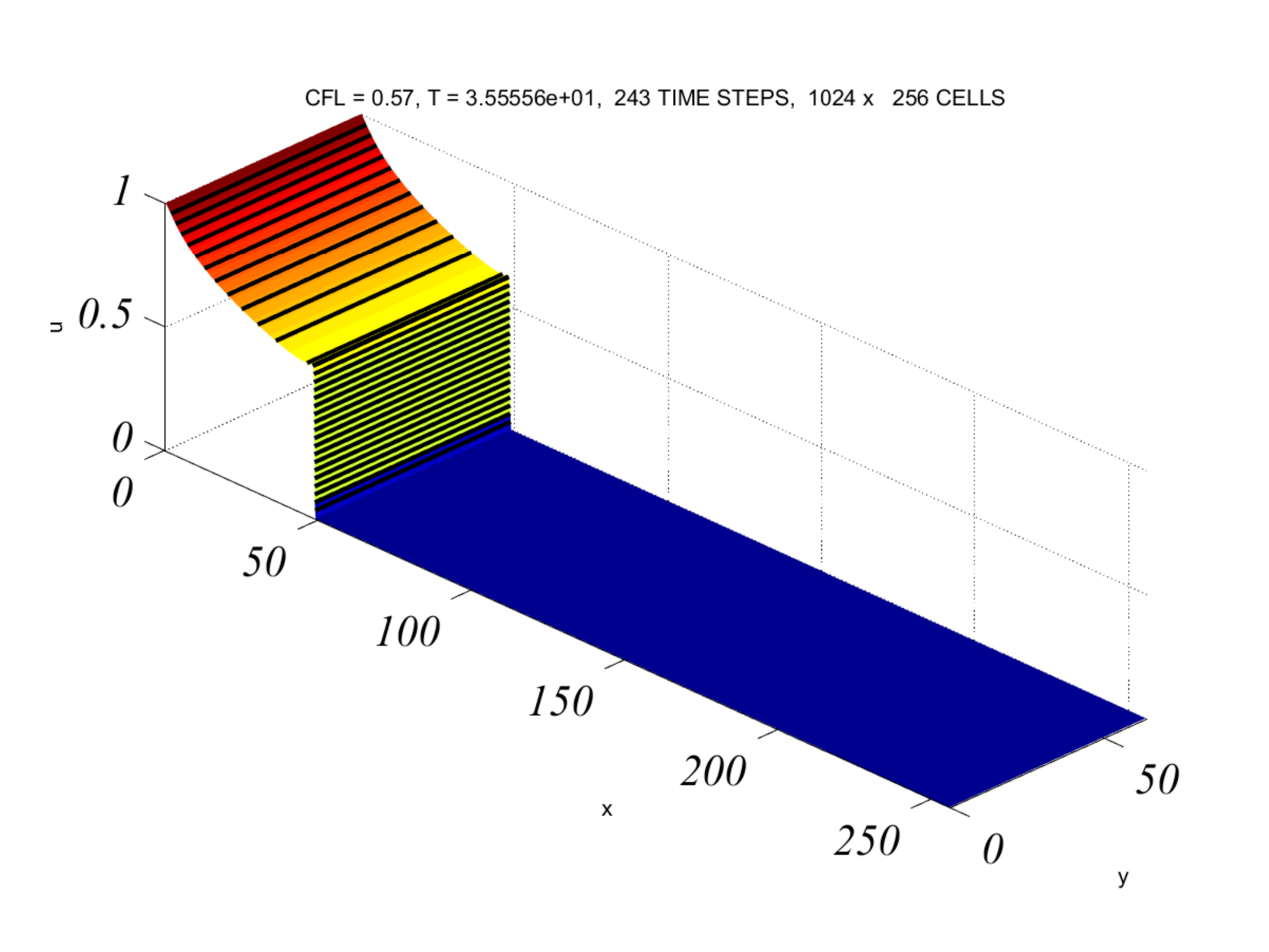}
\includegraphics[scale=0.23]{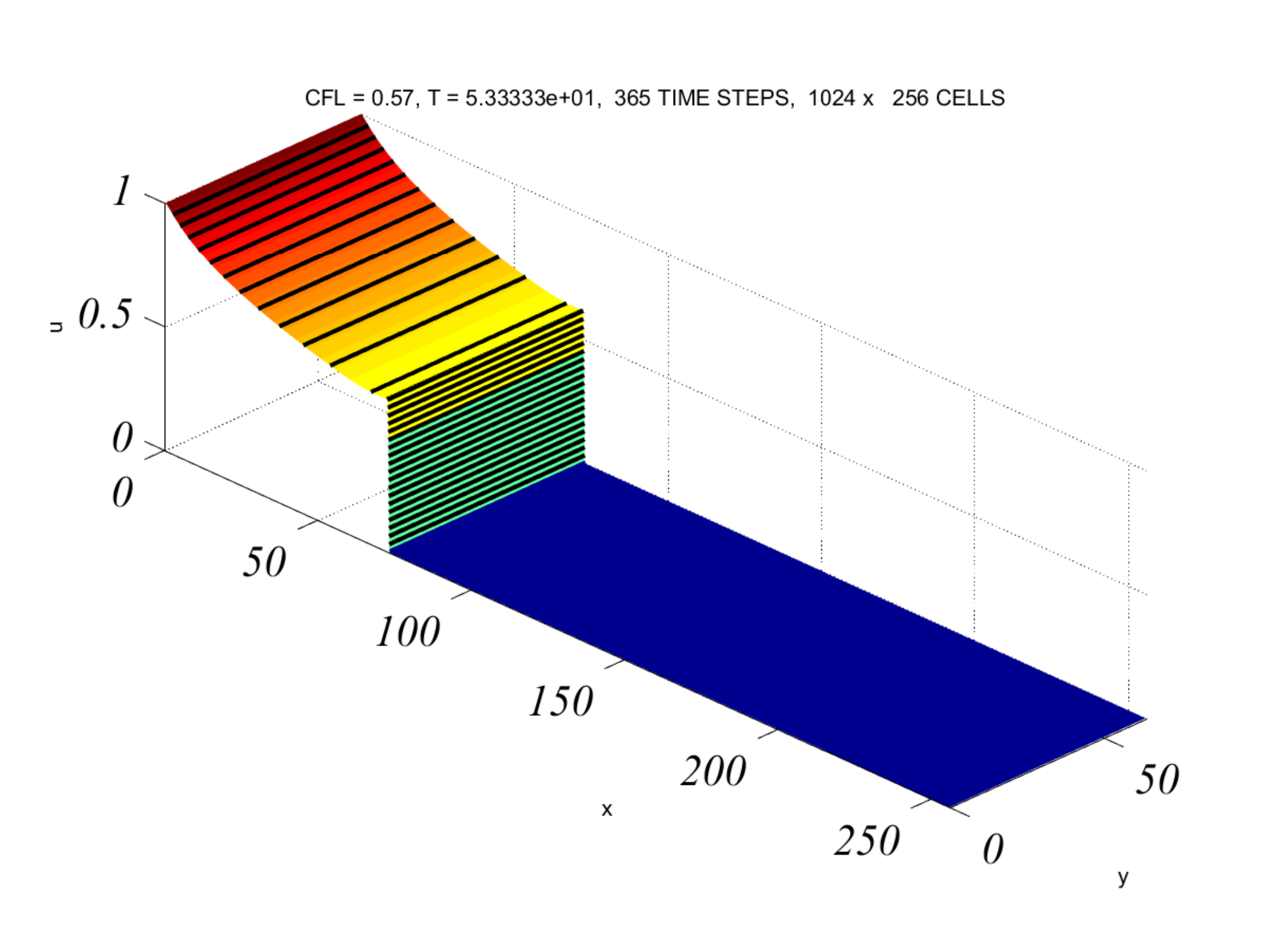} \\
\includegraphics[scale=0.23]{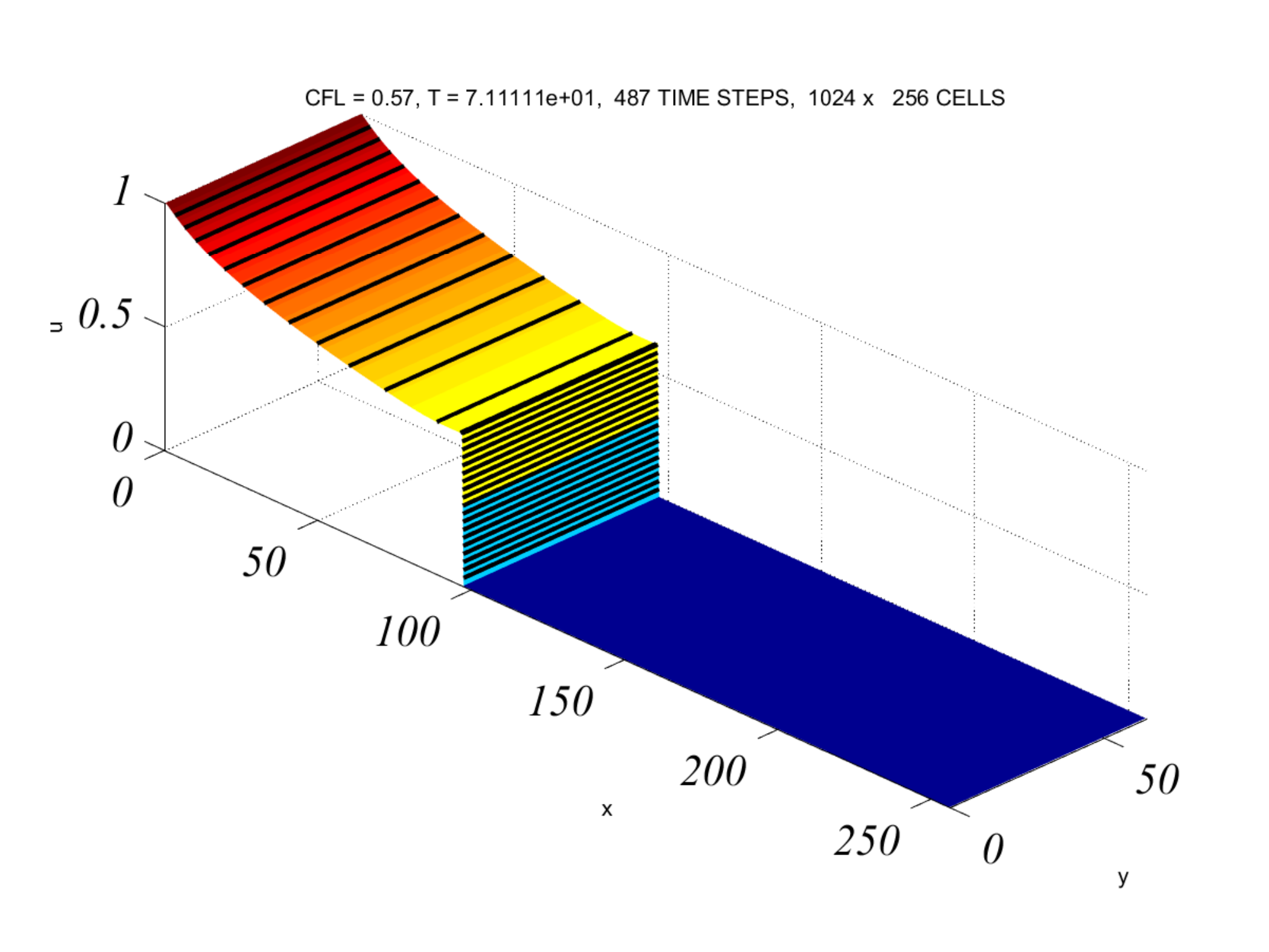}
\includegraphics[scale=0.23]{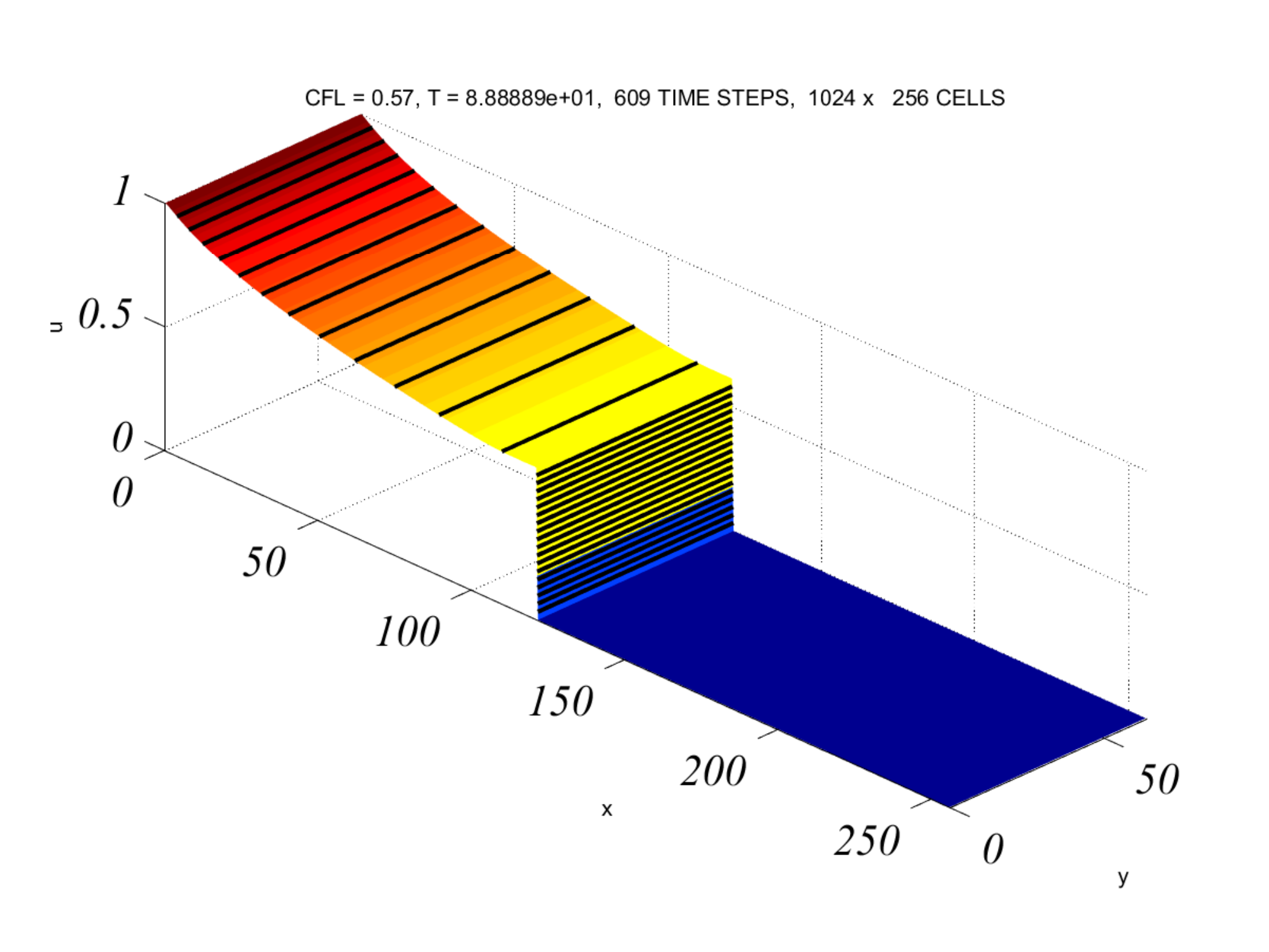}
\includegraphics[scale=0.23]{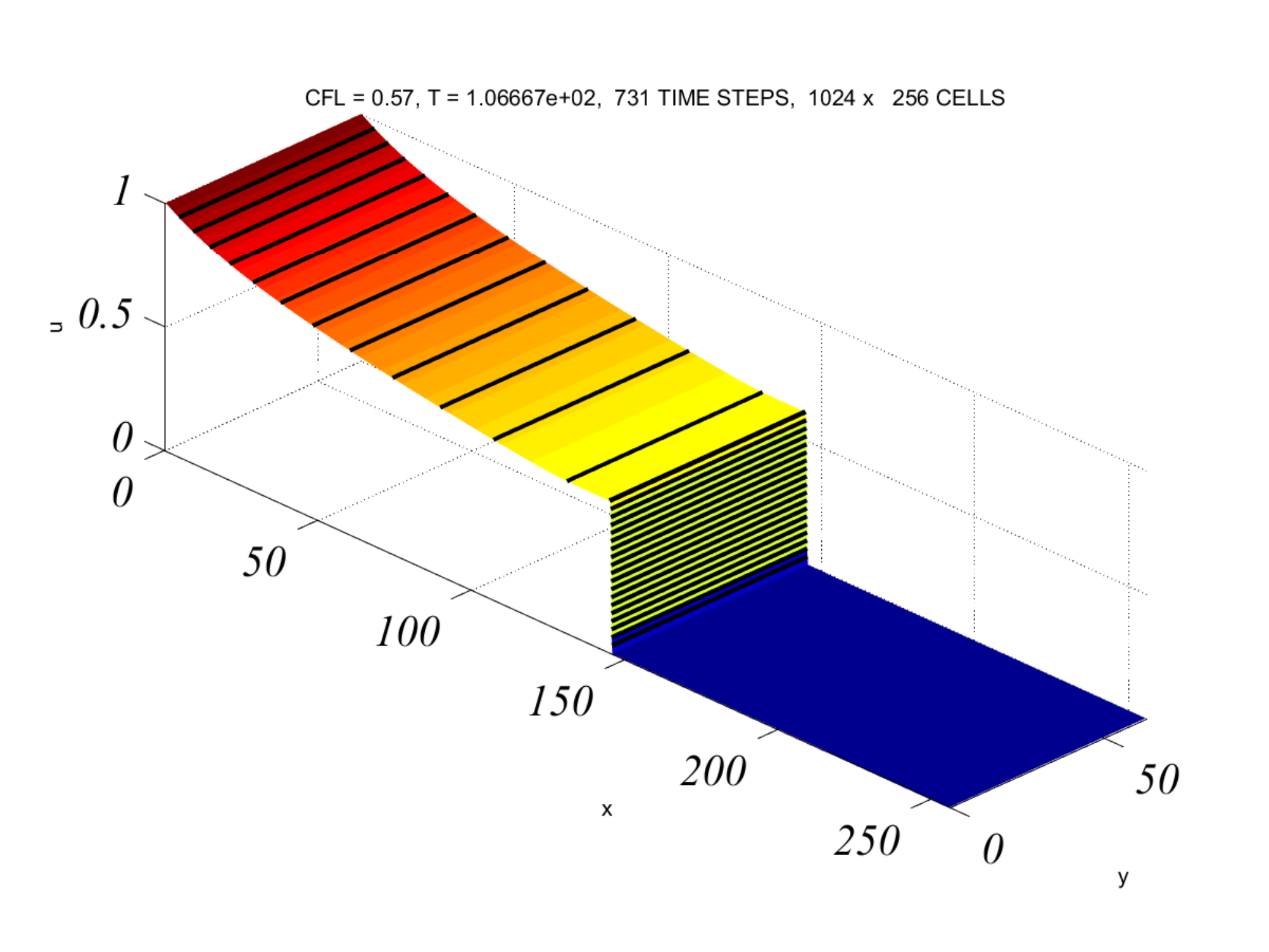} \\
\includegraphics[scale=0.23]{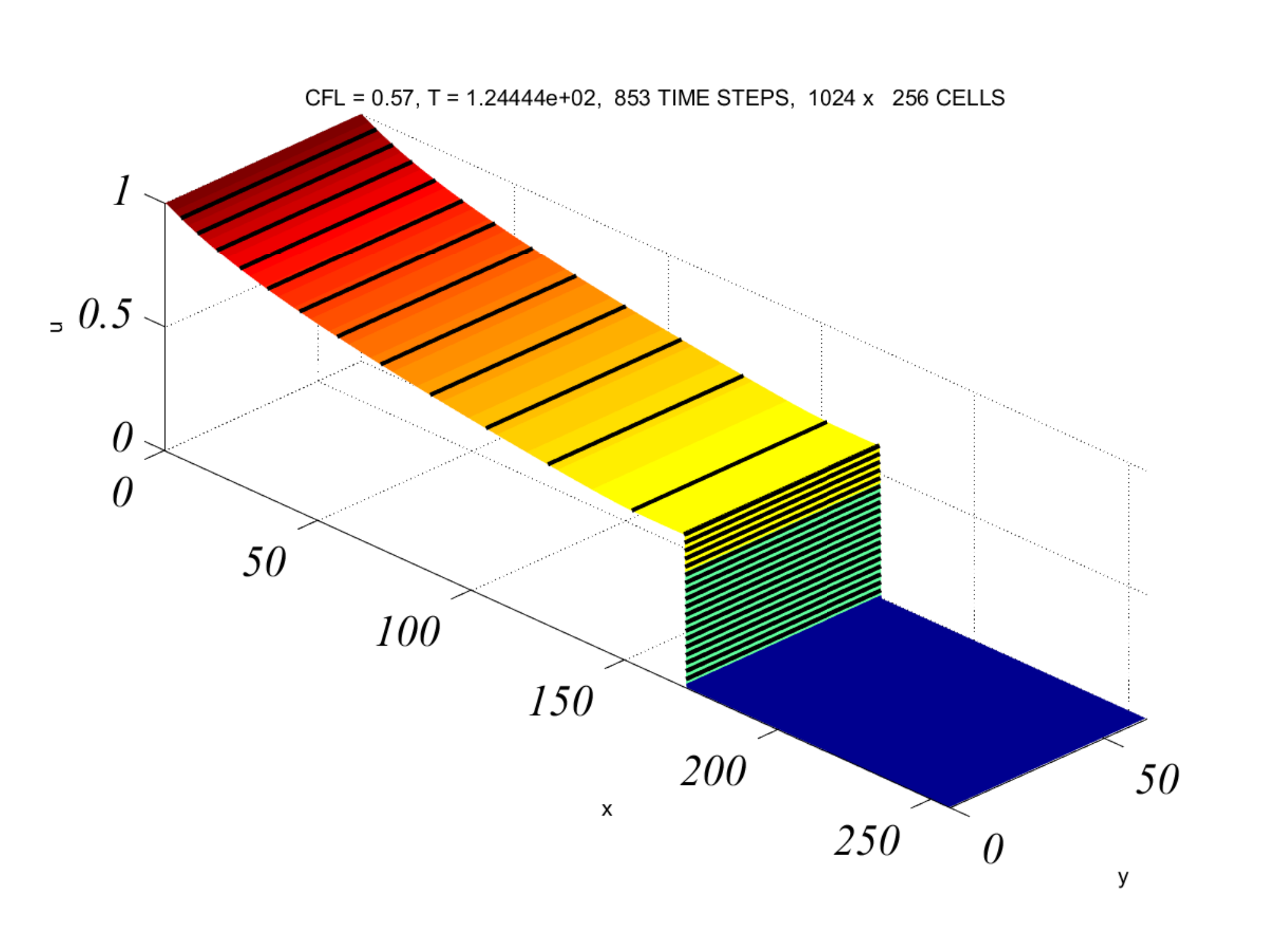}
\includegraphics[scale=0.23]{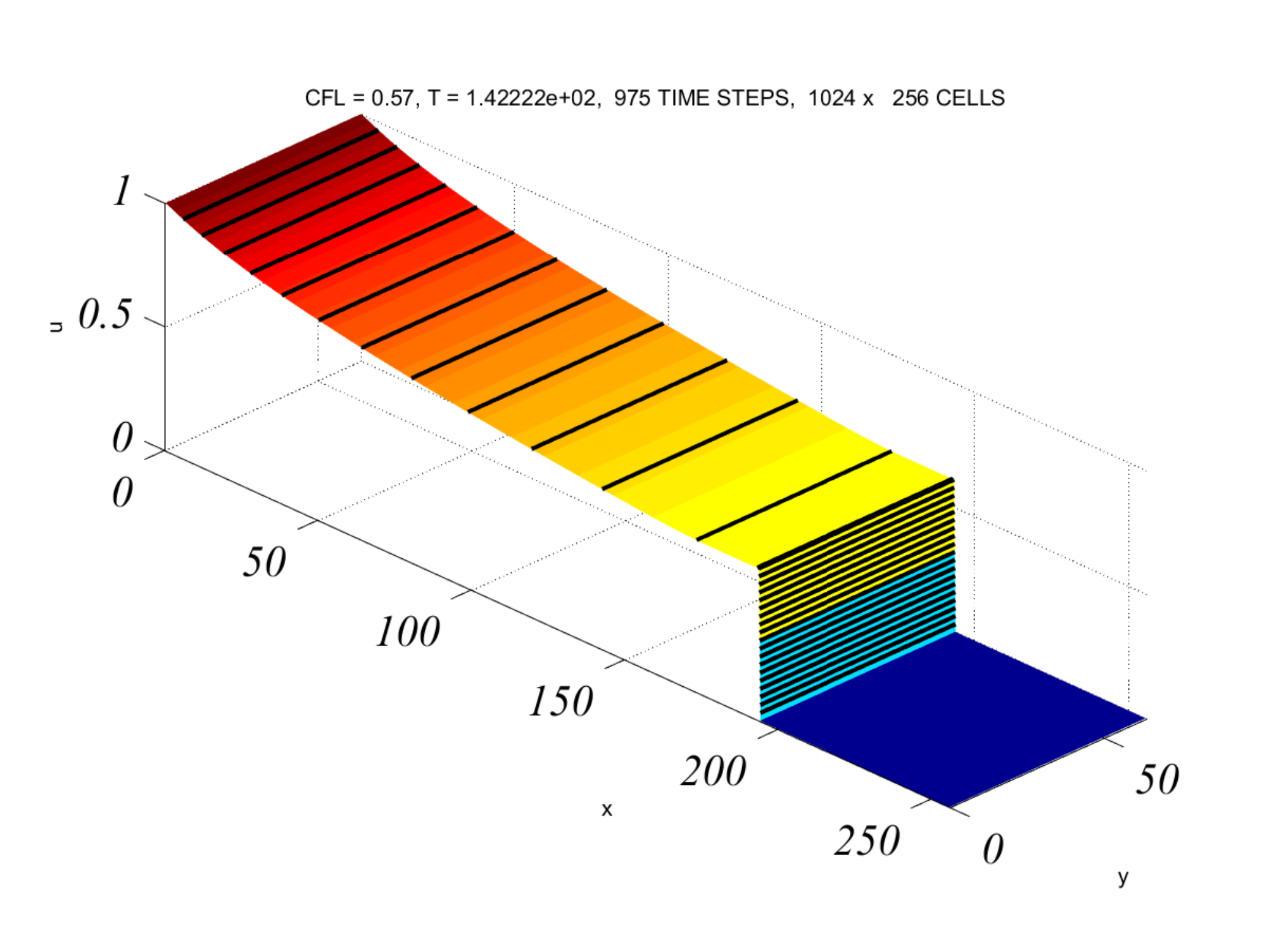}
\includegraphics[scale=0.23]{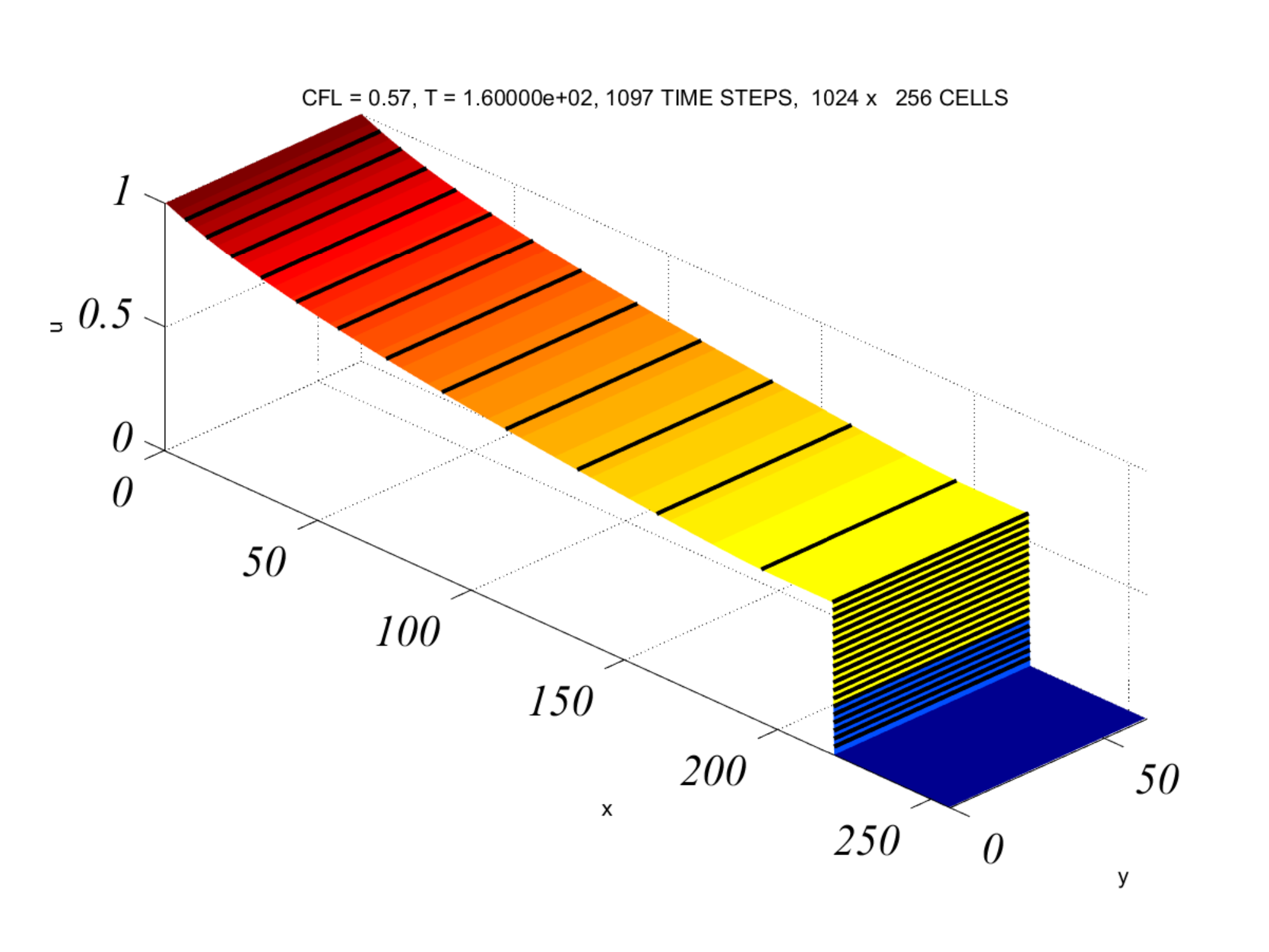}
\caption{Evolution of saturation for \textbf{Test 1} using mesh 
norm $h_x = h_y = 1$. Solutions at times: 
$t = 24,~48,~73,~97,~122,~146,~ 171,~195,~220$.}
\label{fig:Test_1_Sat_Prof}
\end{figure}

\begin{figure}[ht!]
\includegraphics[scale=0.23]{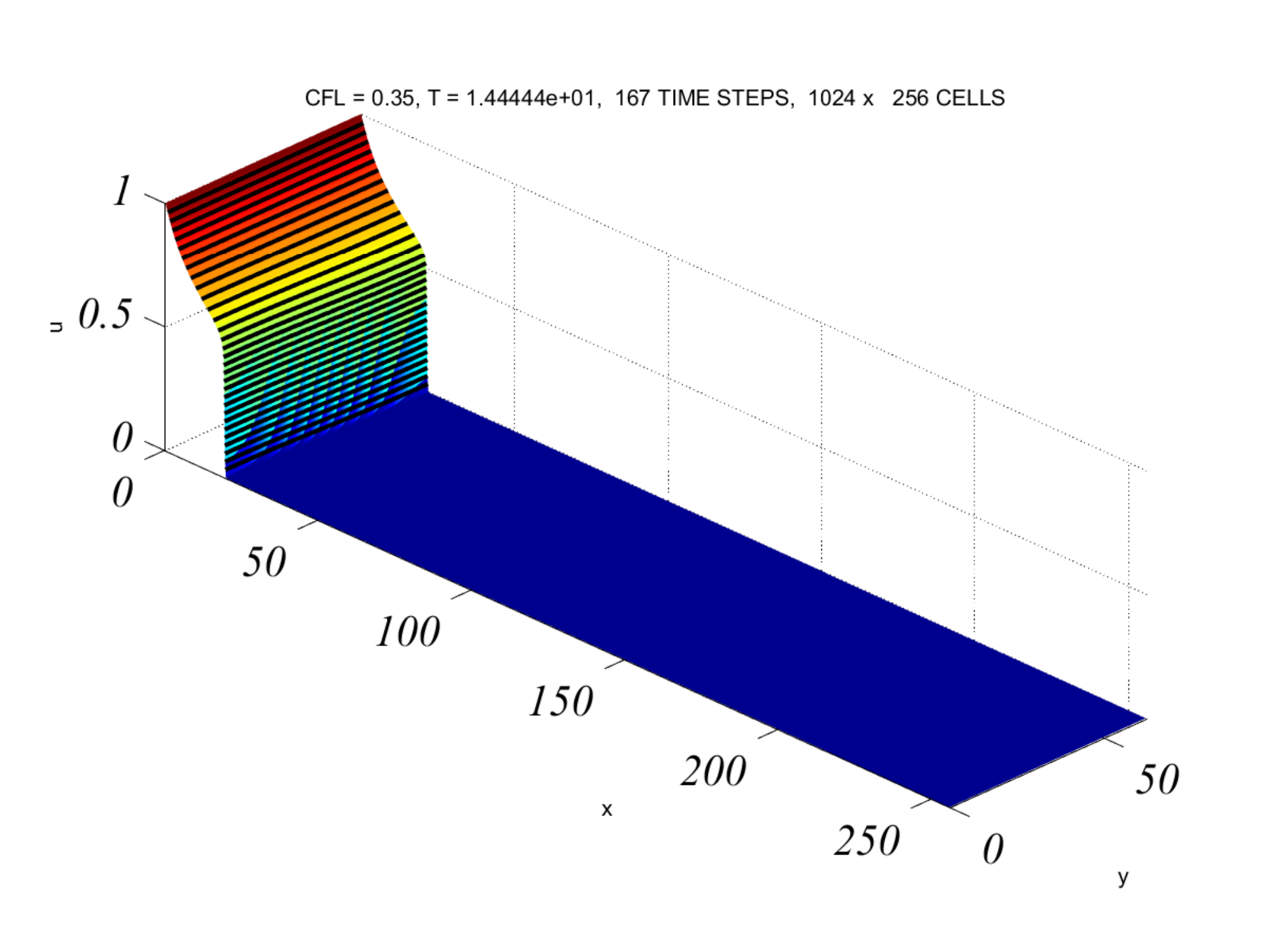}
\includegraphics[scale=0.23]{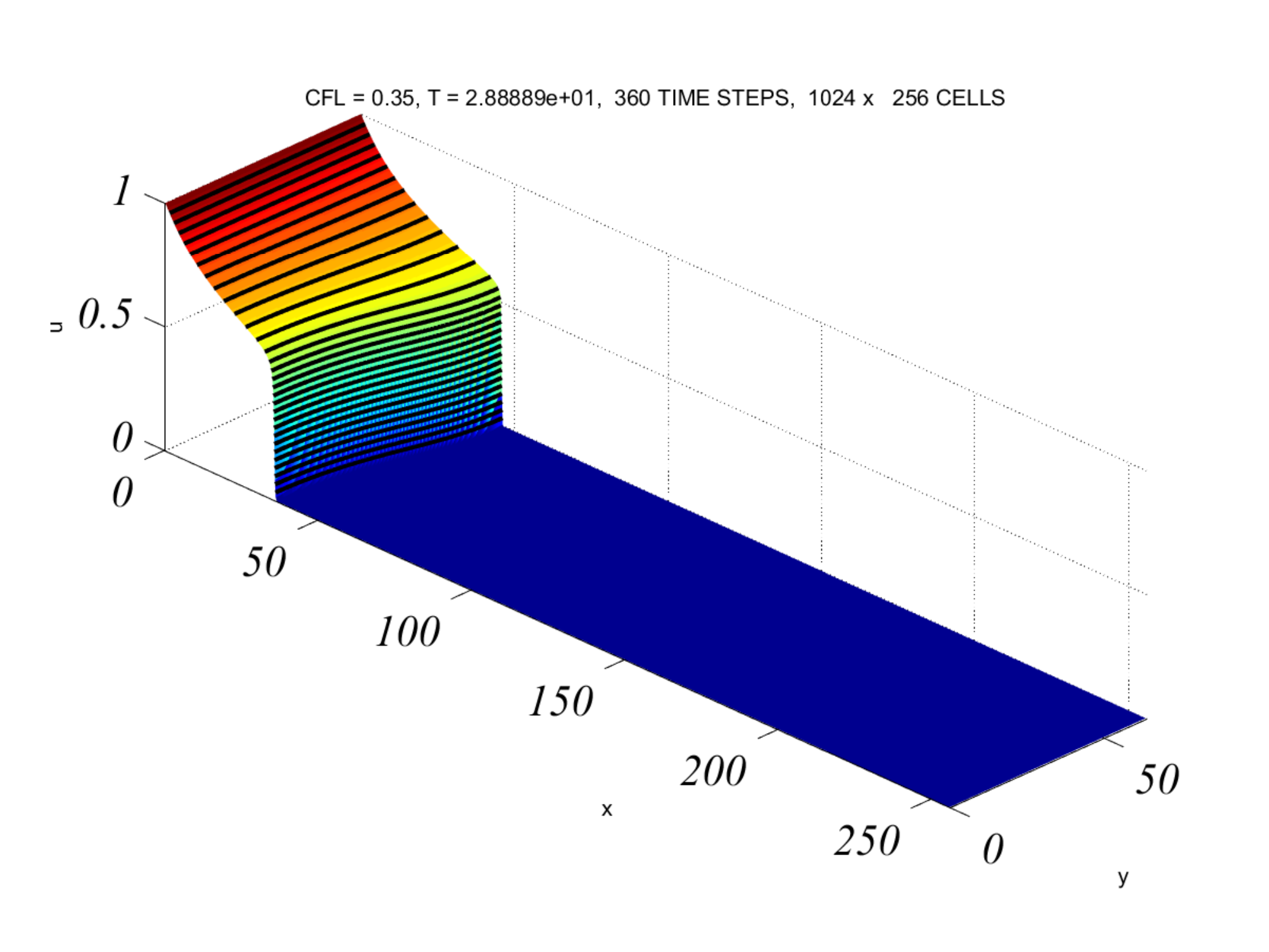}
\includegraphics[scale=0.23]{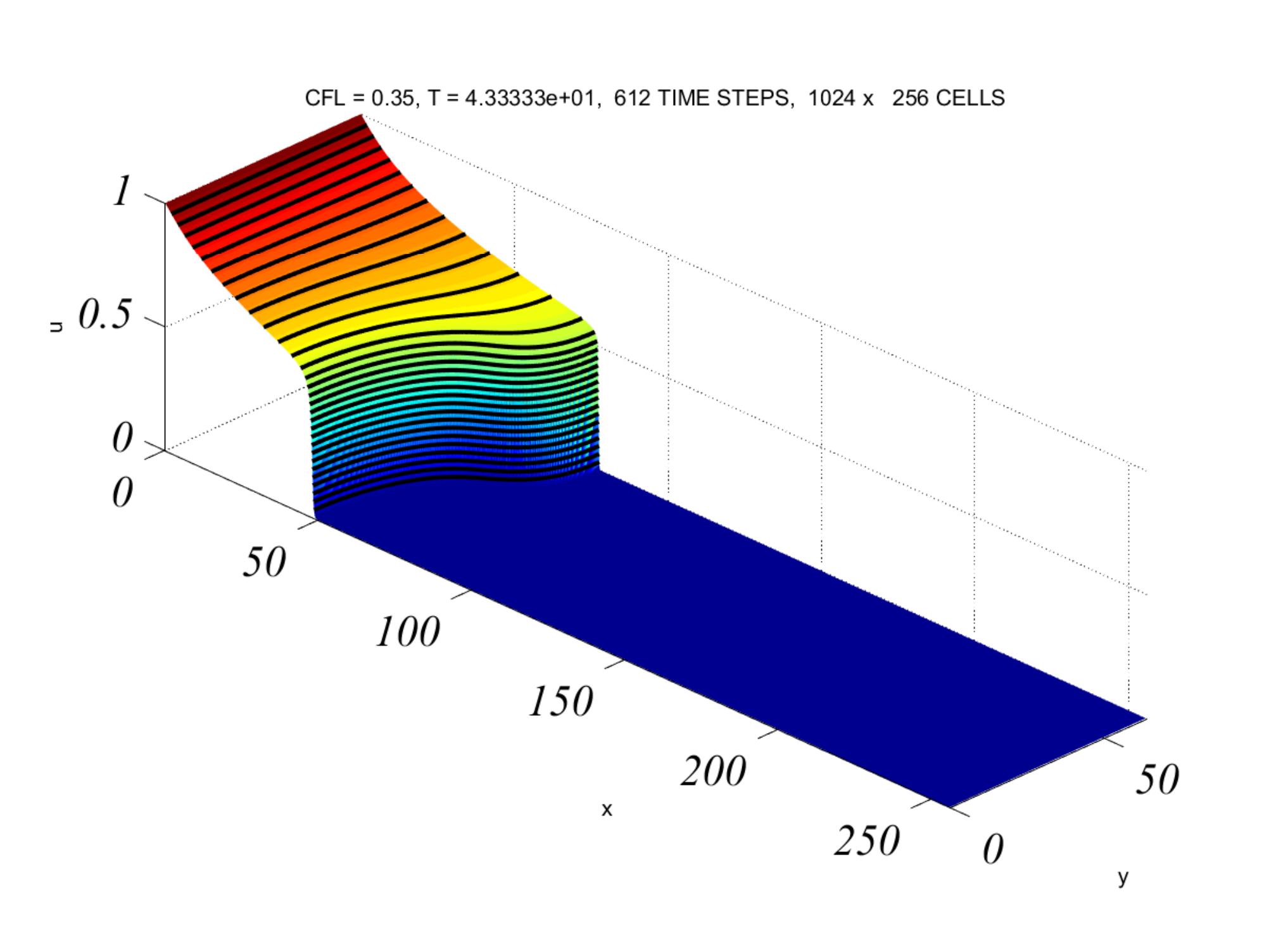}
\includegraphics[scale=0.23]{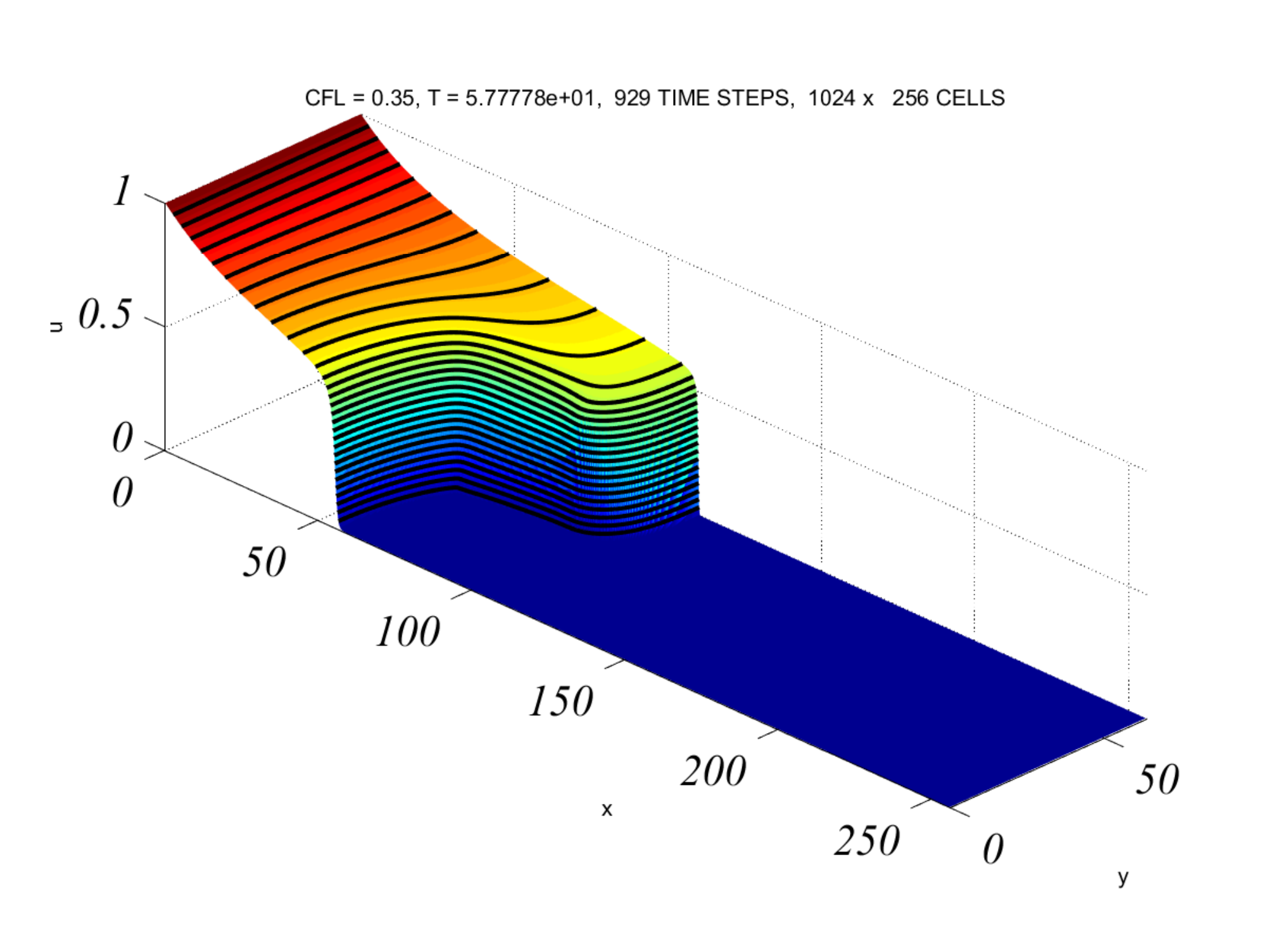}
\includegraphics[scale=0.23]{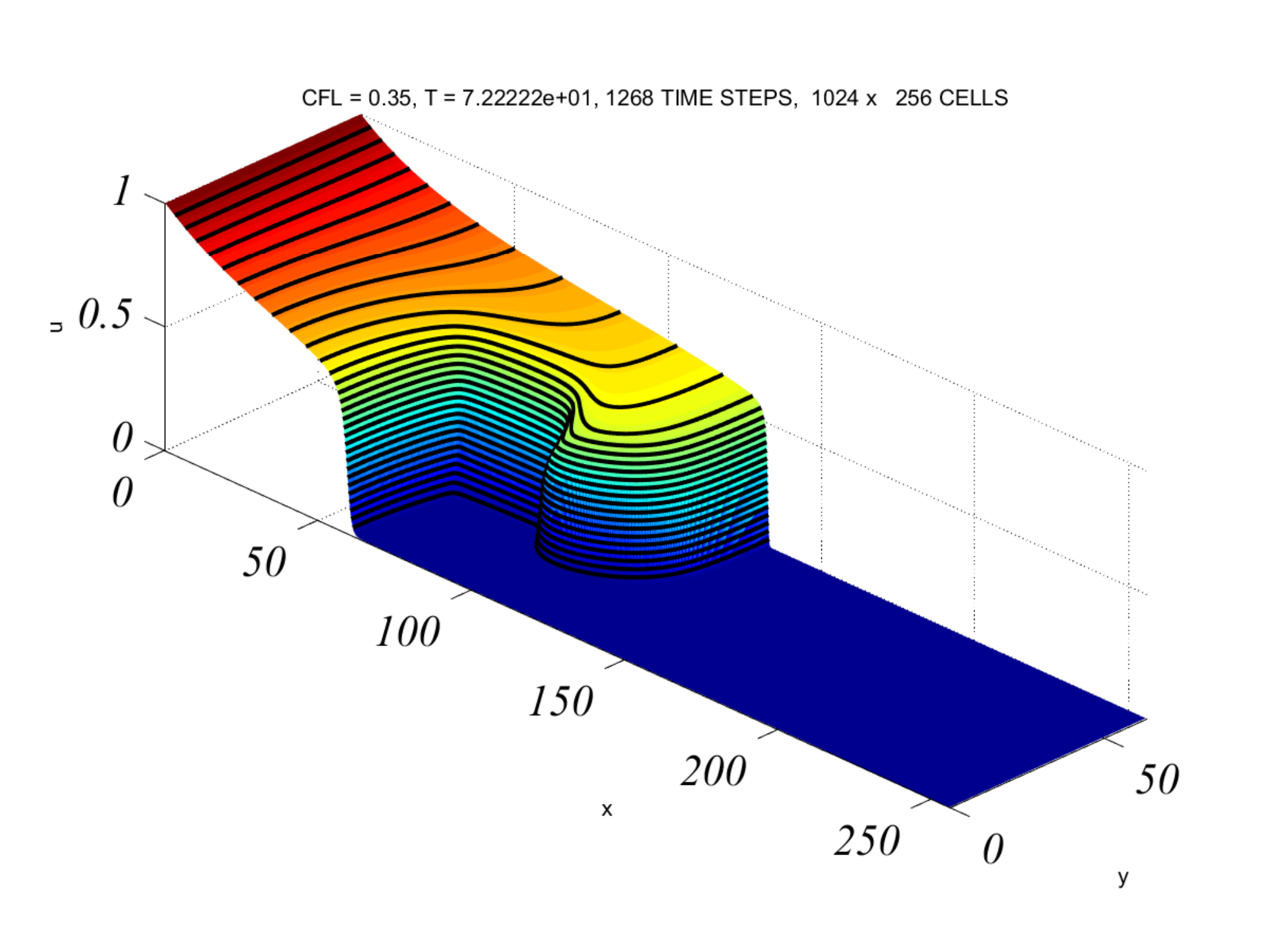}
\includegraphics[scale=0.23]{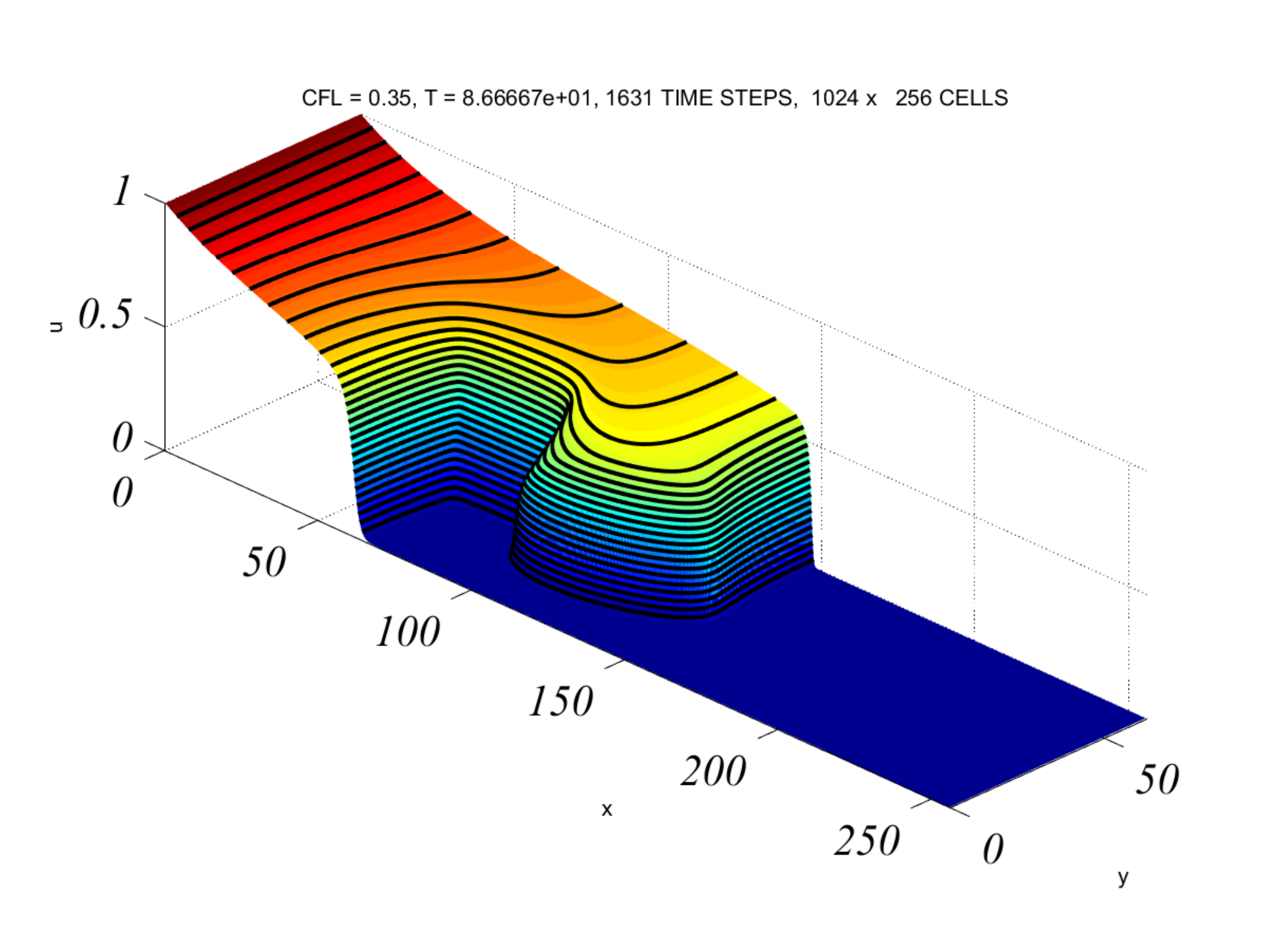}
\includegraphics[scale=0.23]{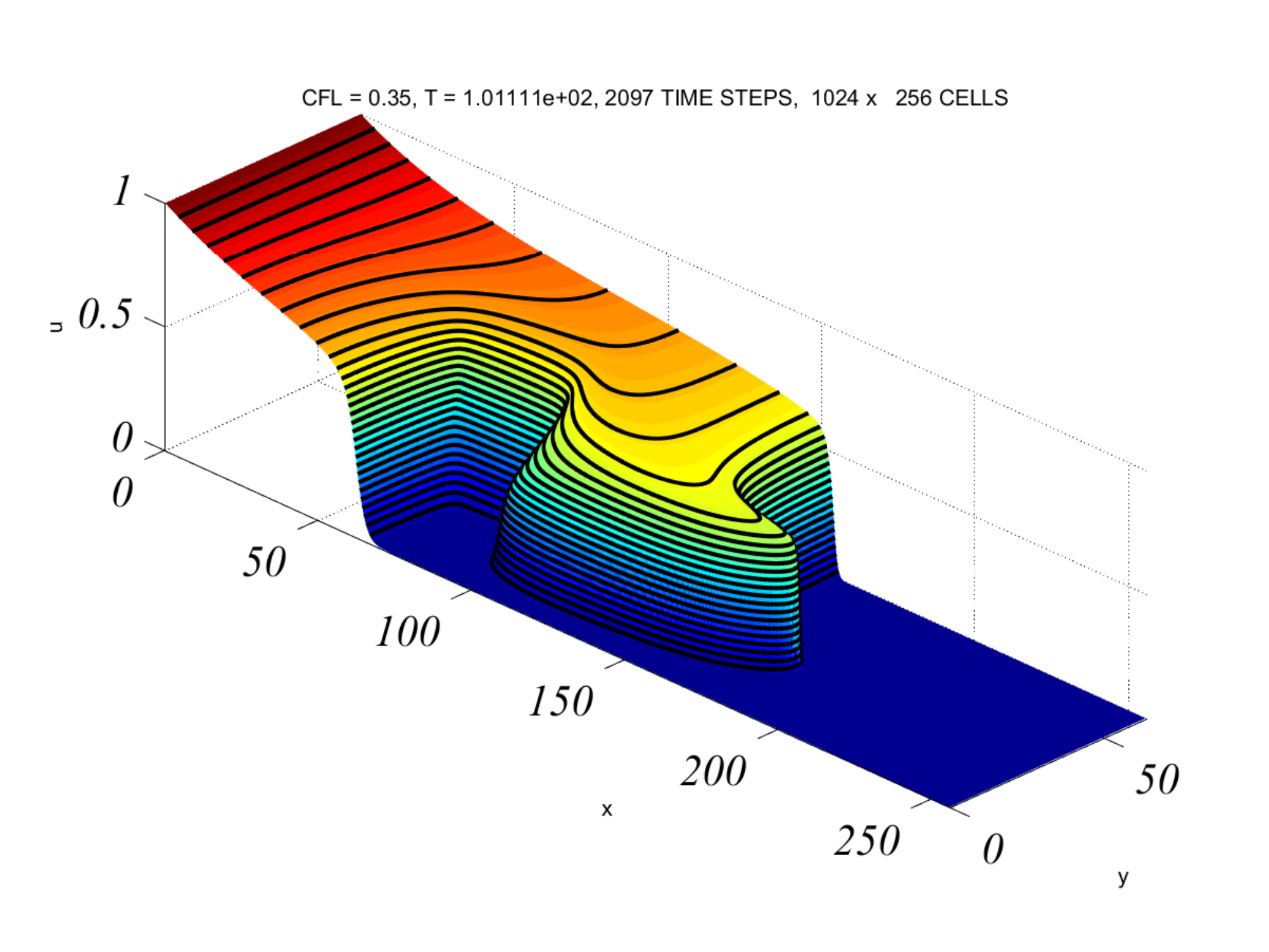}
\hspace{4mm}
\includegraphics[scale=0.23]{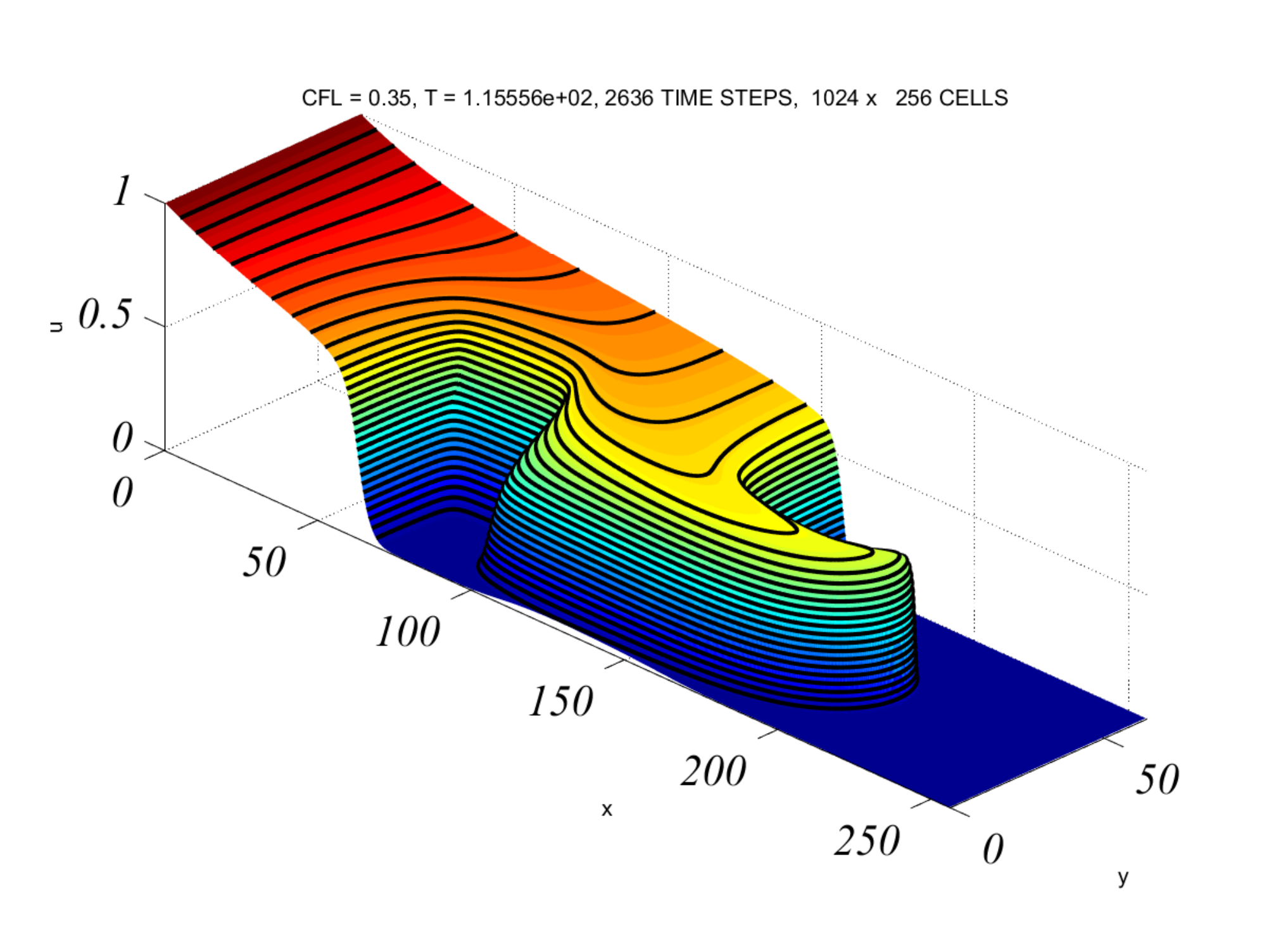}
\hspace{6mm}
\includegraphics[scale=0.23]{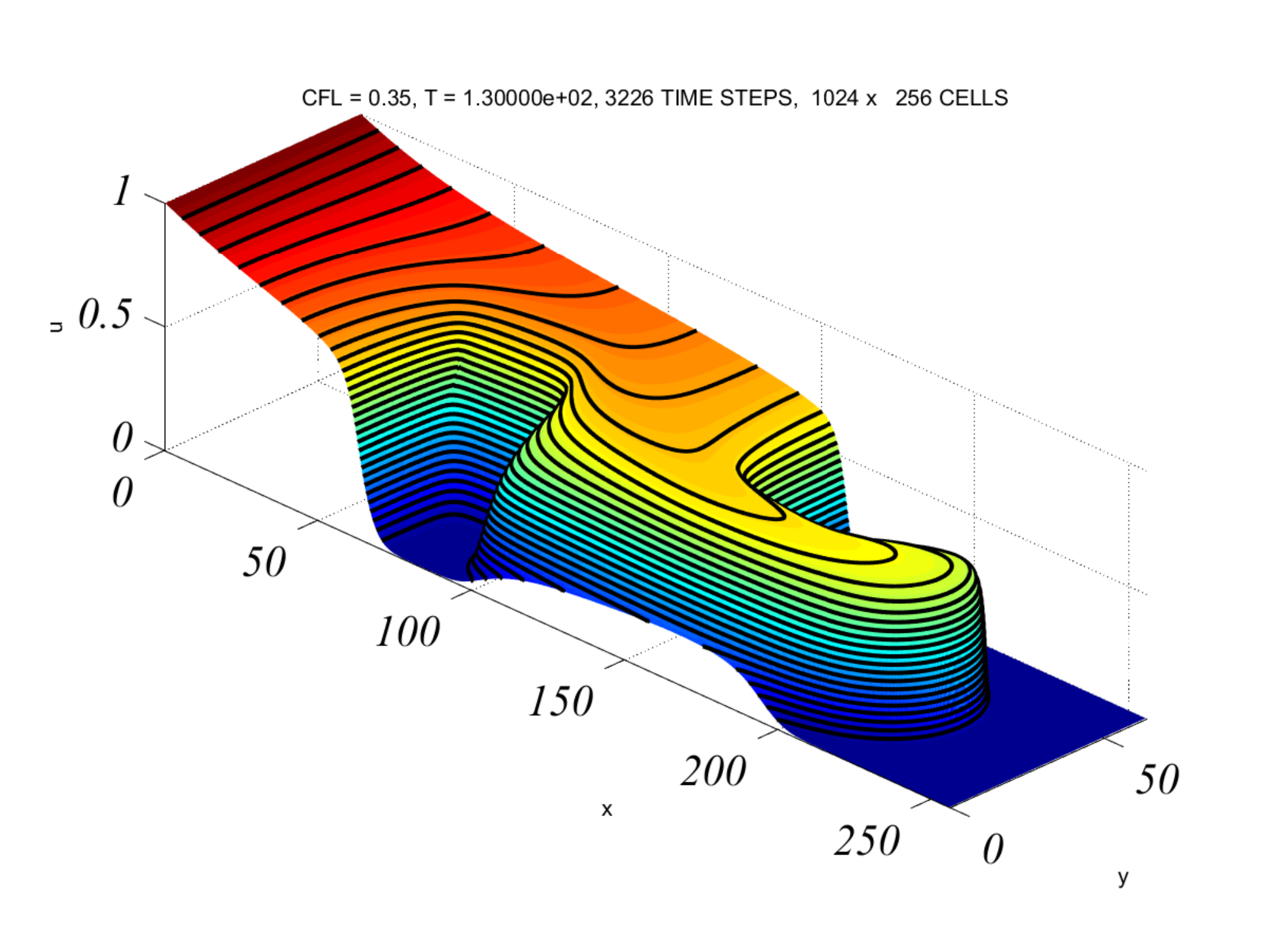}
\caption{Evolution of saturation for \textbf{Test 2} 
using mesh norm $h_x = h_y = 1$. Solutions at times: 
$t = 24,~48,~73,~97,~122,~146,~ 171,~195,~220$.}
\label{fig:Test_2_Sat_Prof}
\end{figure}

\begin{figure}[ht!]
\includegraphics[scale=0.113]{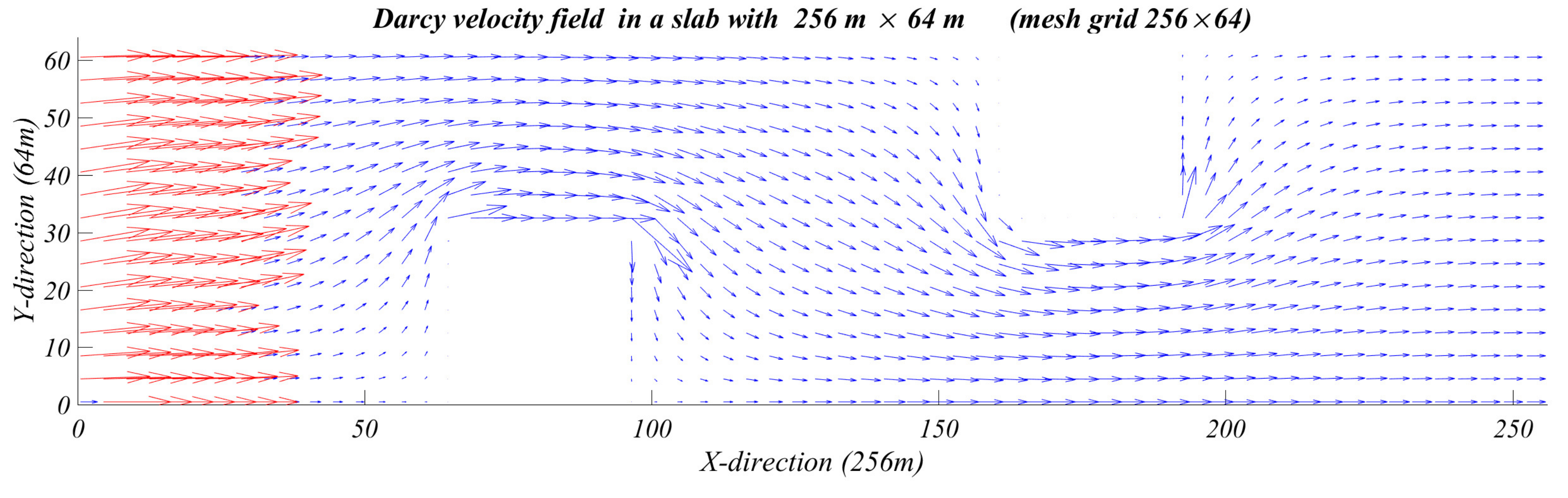}
\includegraphics[scale=0.113]{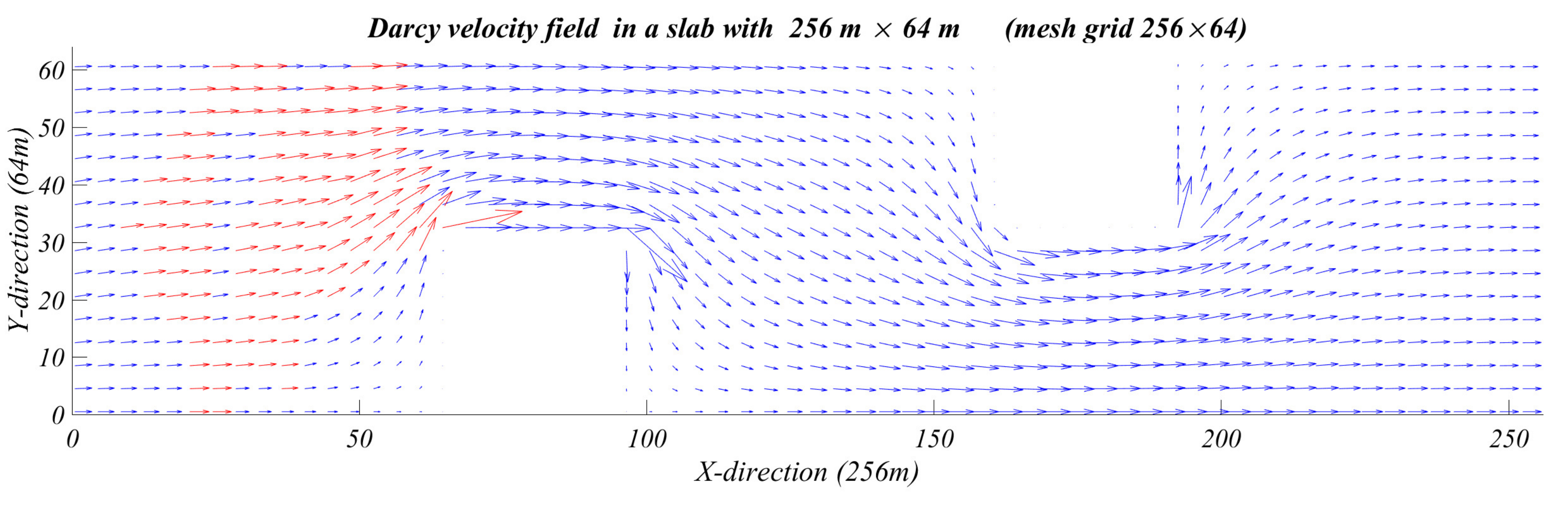}
\includegraphics[scale=0.113]{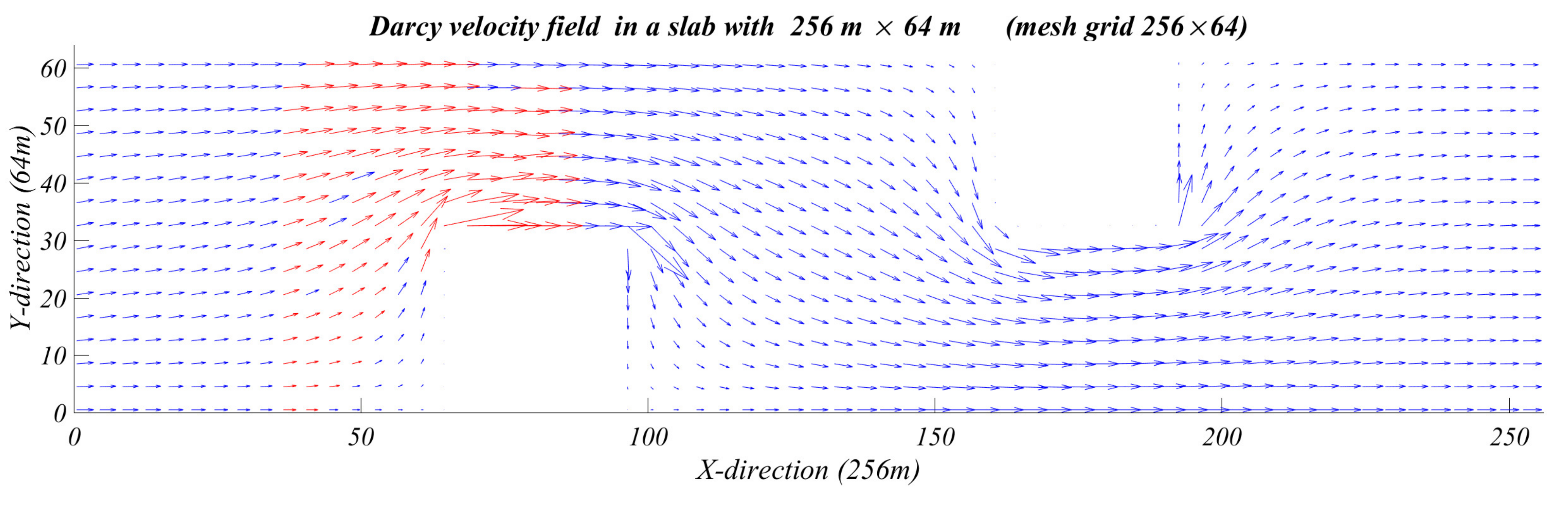} \\
\includegraphics[scale=0.113]{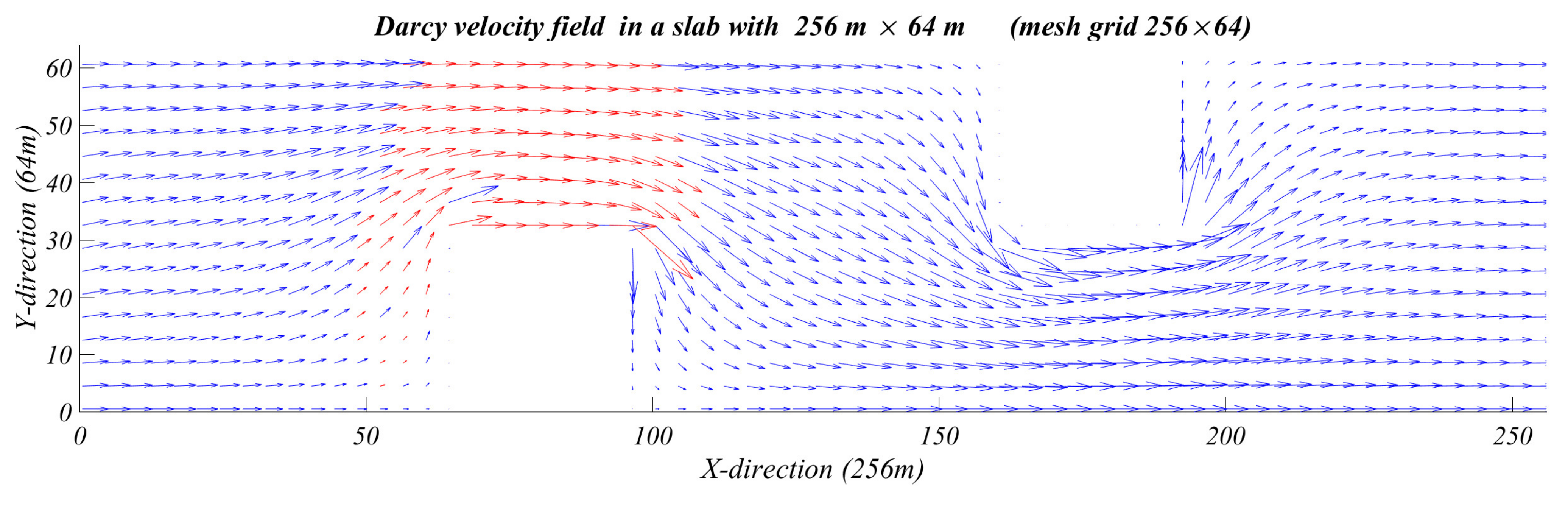}
\includegraphics[scale=0.113]{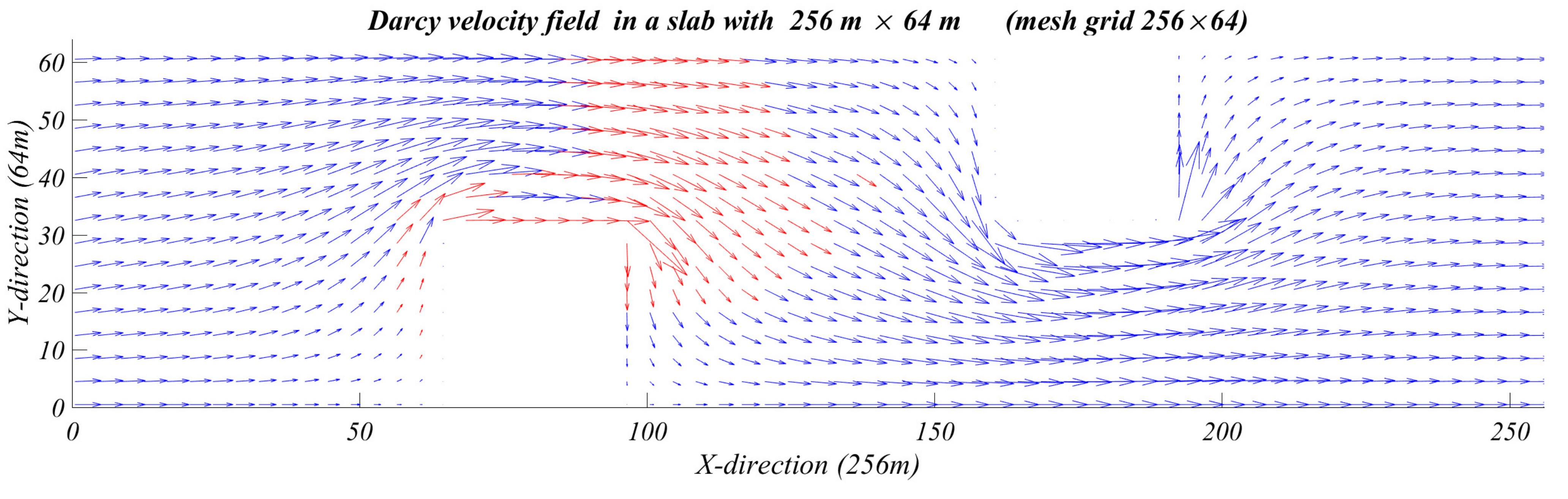}
\includegraphics[scale=0.113]{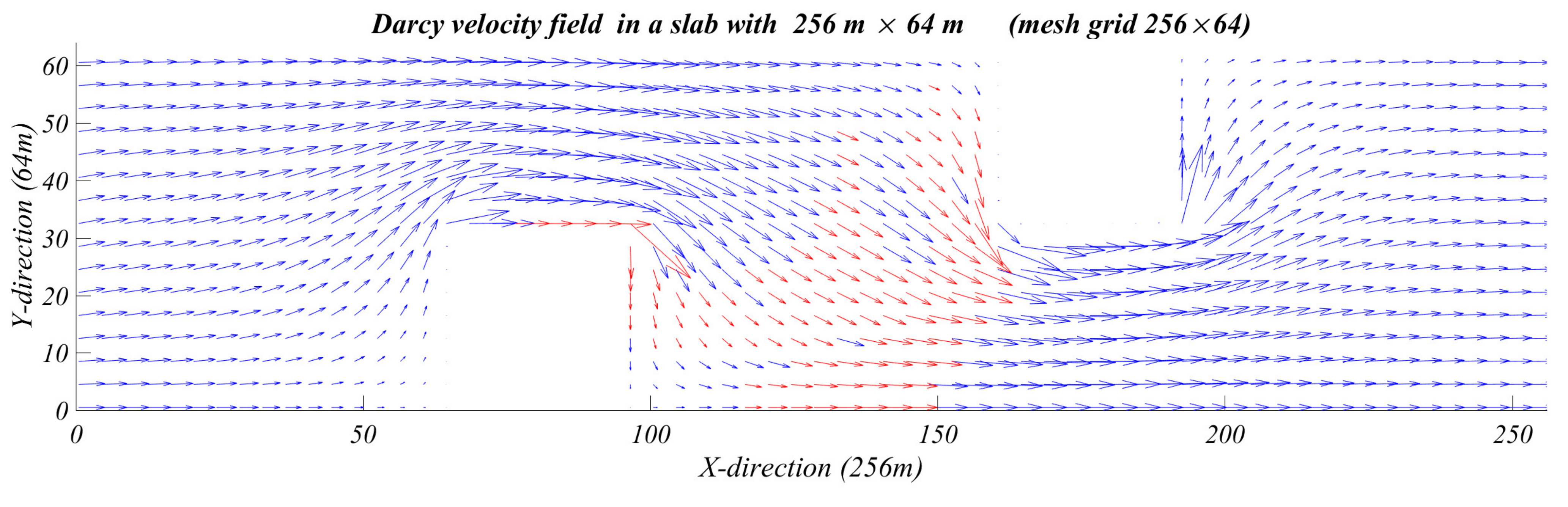} \\
\includegraphics[scale=0.113]{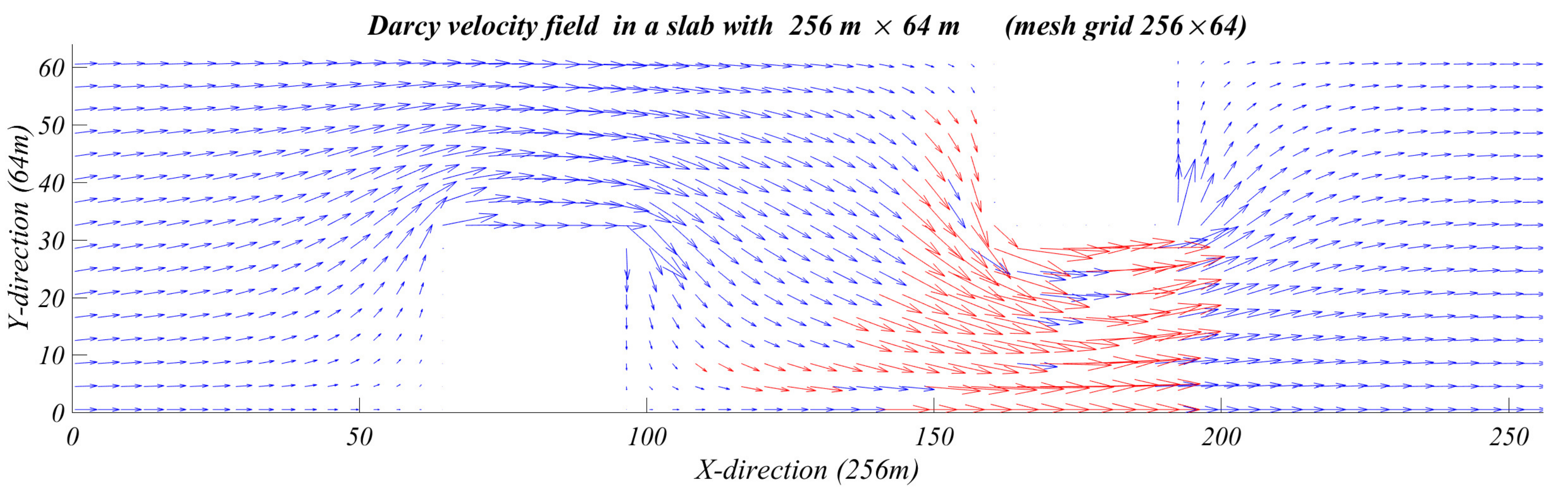}
\includegraphics[scale=0.113]{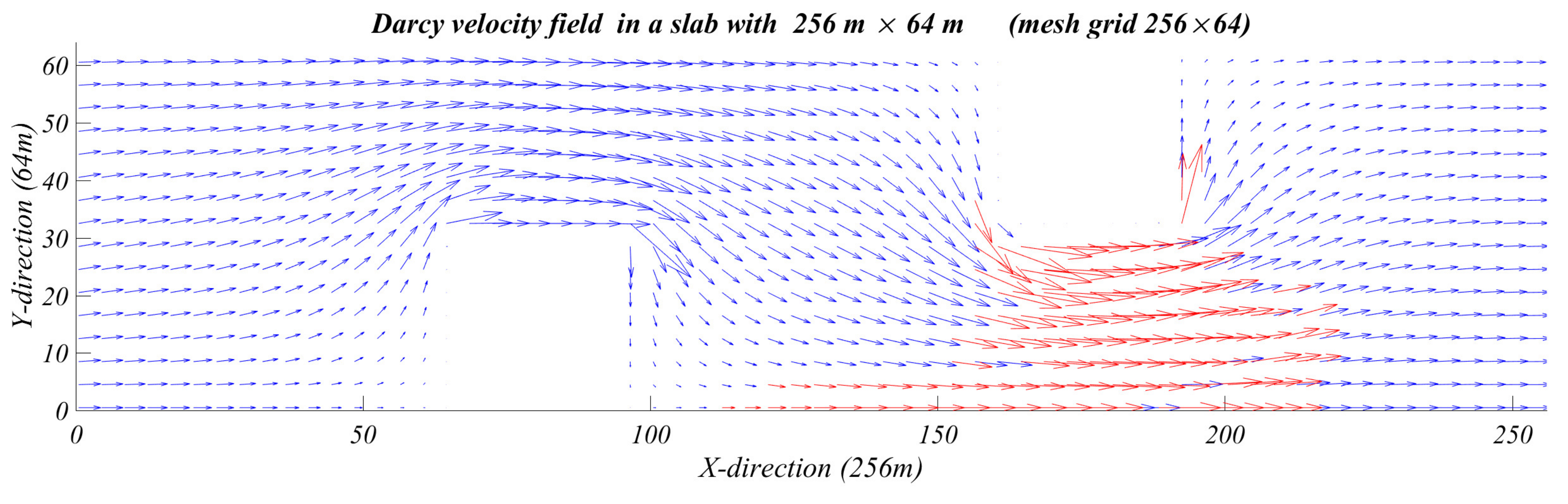}
\includegraphics[scale=0.113]{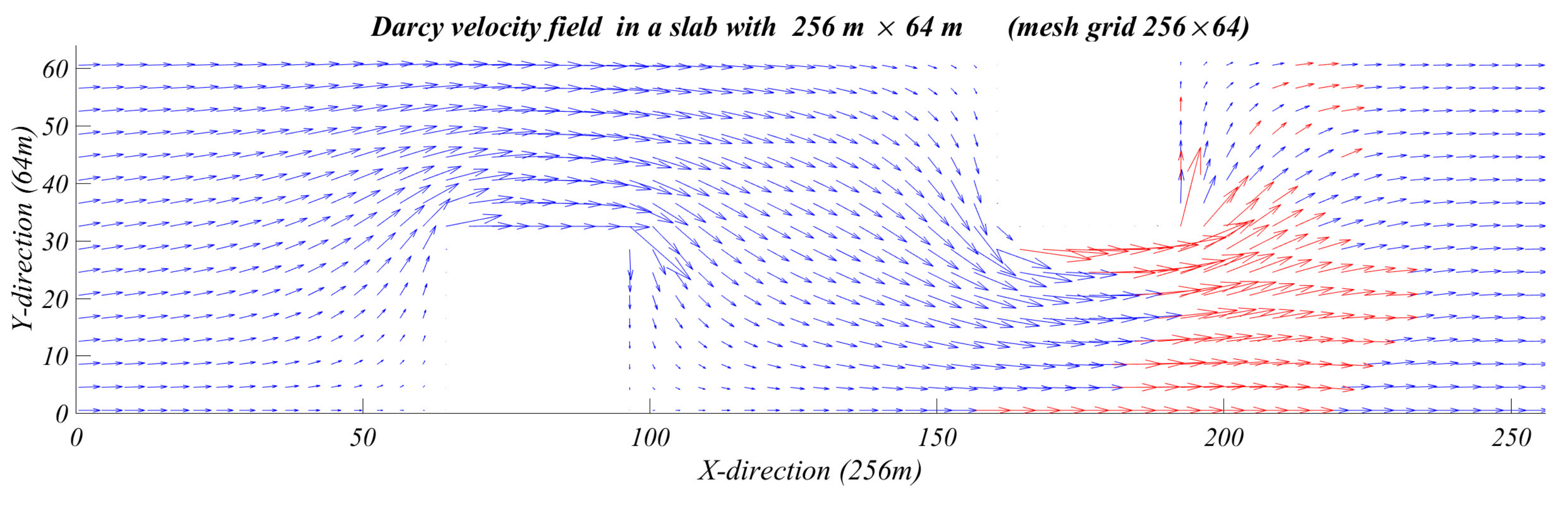}
\caption{Evolution of velocity field for \textbf{Test 2} 
using mesh norm $h_x = h_y = 1$. Solutions at times: 
$t = 24,~48,~73,~97,~122,~146,~ 171,~195,~220$.}
\label{fig:sfig2}
\end{figure}

\begin{figure}[ht!]
\centering 
\includegraphics[scale=0.18]{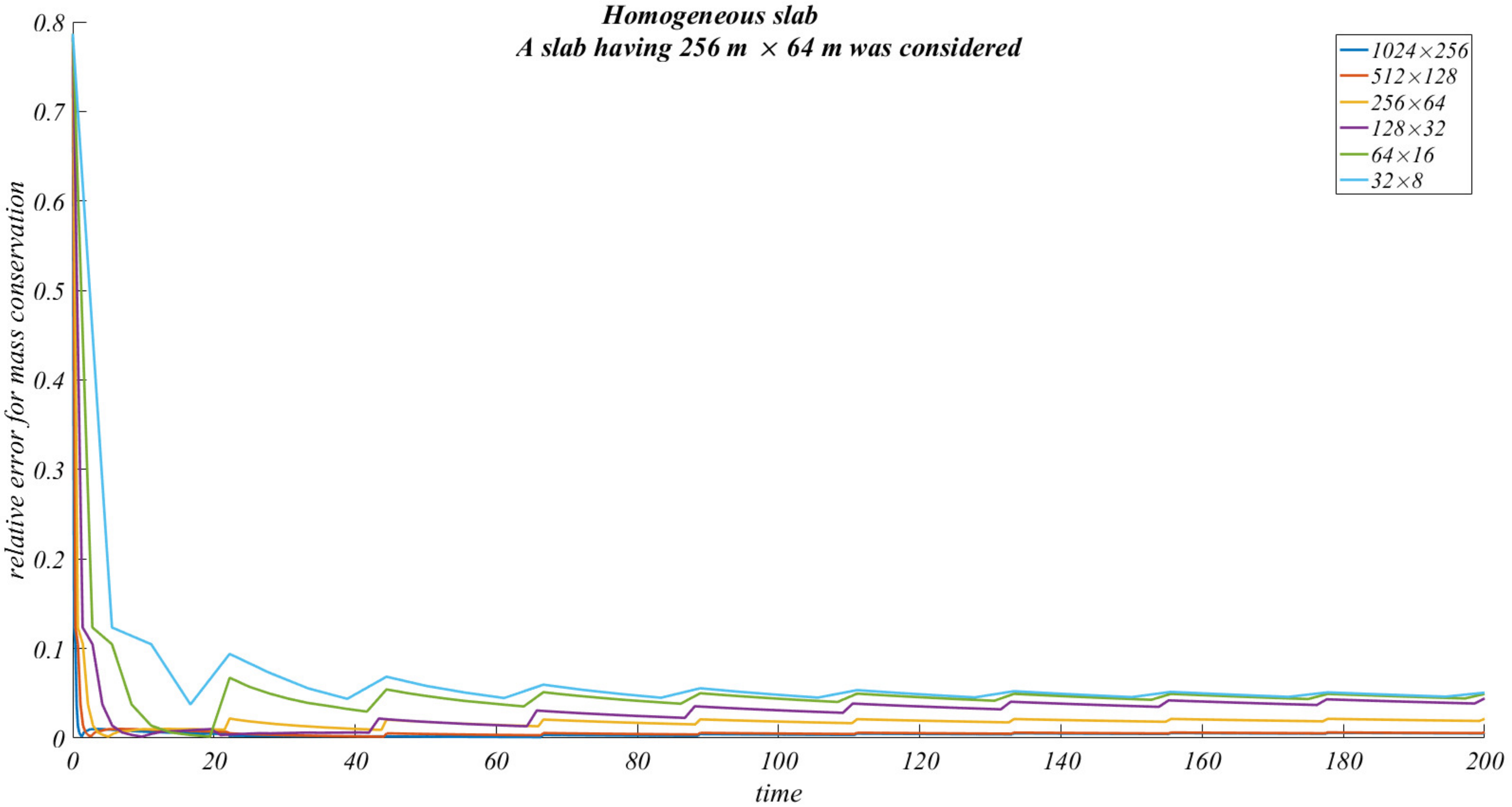}
\includegraphics[scale=0.18]{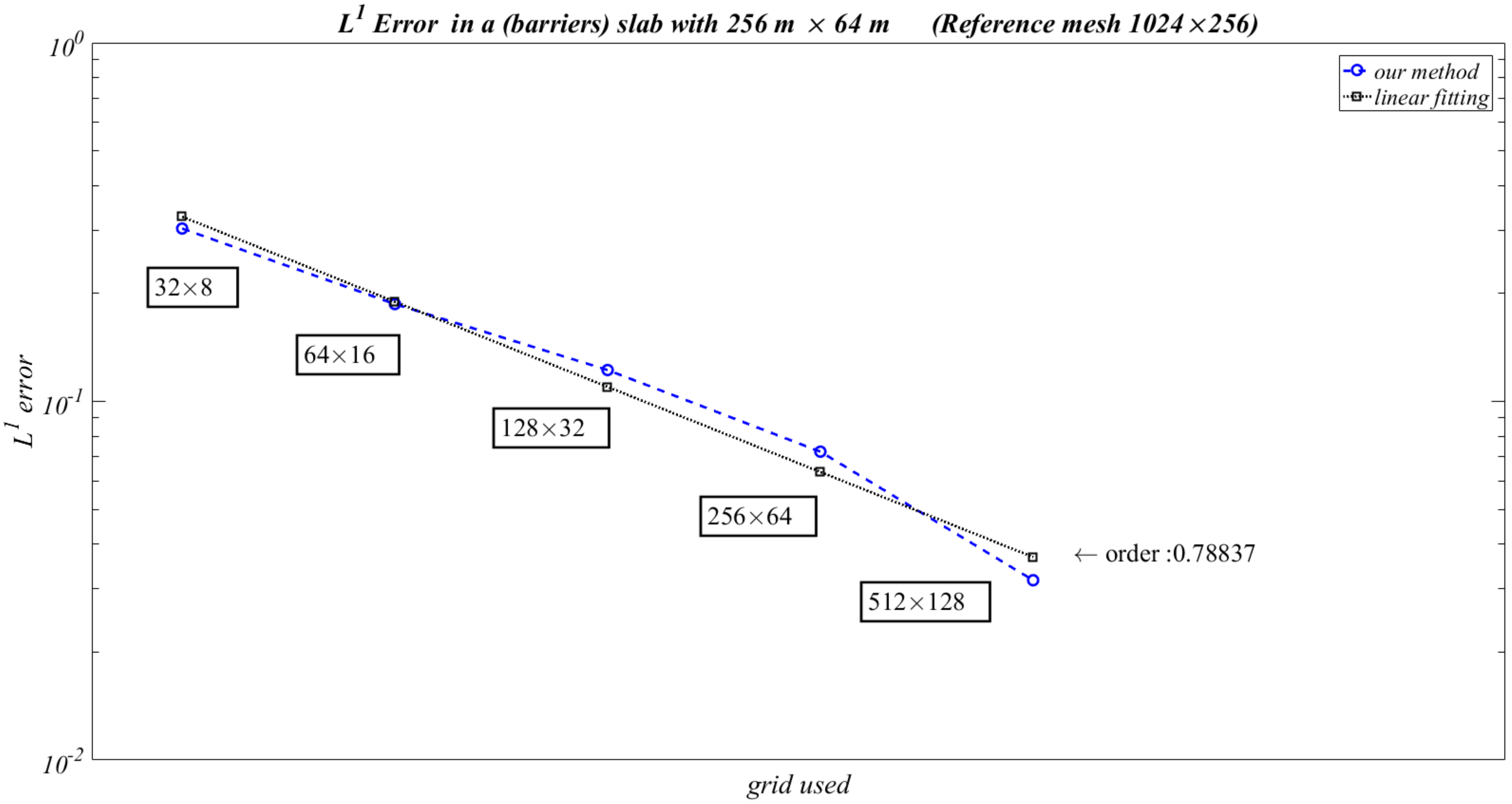}
\caption{A 2D homogeneous slab problem coupling test for the Darcy problem with
hyperbolic-transport having 256m $\times$ 64m: on the top a decreasing of
the relative error of mass under a mesh refinement study and on the bottom
we see evidence of numerical convergence of the full 
Darcy-hyperbolic-transport two-phase flow system.}
\label{Fig:Rel_Err_Mass_Test_1} 
\end{figure}

\begin{figure}[ht!]
\centering 
\includegraphics[scale=0.18]{MassRelativeErrorb-eps-converted-to.pdf}
\includegraphics[scale=0.18]{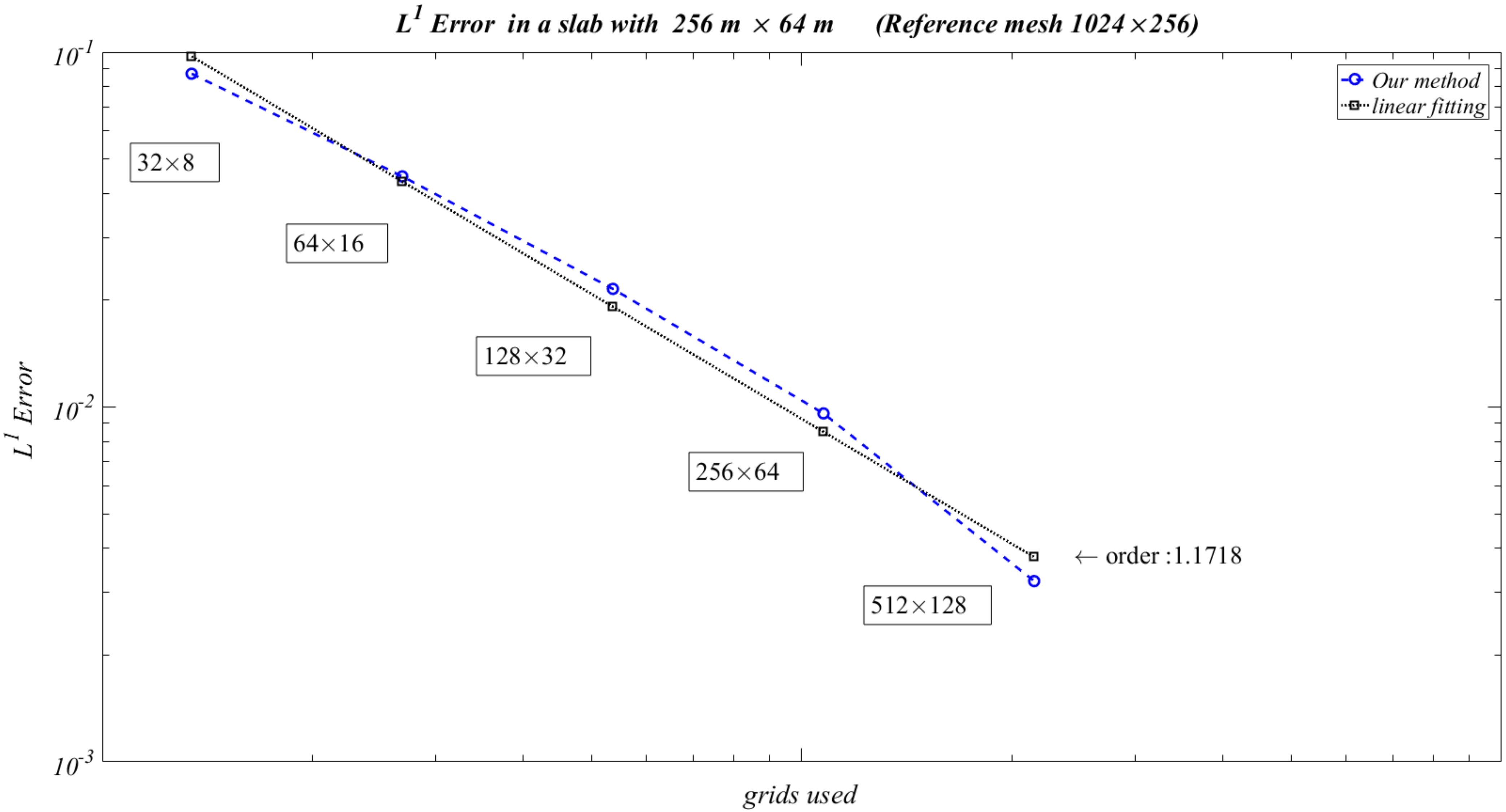}
\caption{A 2D barrier slab problem coupling test for the Darcy problem with
hyperbolic-transport having 256m $\times$ 64m: on the top a decreasing of
the relative error of mass under a mesh refinement study and on the bottom
we see evidence of numerical convergence of the full 
Darcy-hyperbolic-transport two-phase flow system.}
\label{Fig:Rel_Err_Mass_Test_2} 
\end{figure}

For both flow situations, we present the evolution 
of the waves front interaction of saturation as 
evolve in time; see frames in Figure \ref{fig:Test_1_Sat_Prof} 
(homogeneous medium) and \ref{fig:Test_2_Sat_Prof} 
(heterogeneous barrier high high-constrast permability). 
The evolution of velocity field correponding to the 
heterogeneous barrier flow situation is displayed in 
Figure \ref{fig:sfig2}. Indeed, we also present the 
computation of the relative mass errors computed 
with our multiscale coupling procedure through 
situations Test 1 and 2, depicted in Figure 
\ref{Fig:Rel_Err_Mass_Test_1} (homogeneous medium) 
and \ref{Fig:Rel_Err_Mass_Test_2} (heterogeneous 
barrier). We also present a numerical convergence 
study that corroborates our findings. Based on 
the reported results, we were able to show the 
promising methodology on the conservation 
properties in multiple scale coupling and 
simulation for Darcy flow with hyperbolic-transport 
in complex flows.

\newpage
\section{Concluding remarks and perspectives}\label{cpfinal}

In this paper, we are concerned with modeling, simulation 
and numerical analysis for approximate solutions in 
multiscale nonlinear PDE related to highly complex 
problems. A better comprehension of multiscale fluid 
flow in subsurface is very hard, challenging and 
undoubtedly still of current events and demand 
innovative multiscale approaches, since ingenious 
difficulties stem from the lack of regularity of 
solutions. We revisited a novel volumetric {\it 
locally conservative} and residual-based Lagrange 
multipliers saddle point reformulation of the standard 
high-order finite method , clarifying and simplifying 
the presentation of its conservative properties. A 
new robust and accurate dynamic forward tracking 
Lagrangian-Eulerian scheme for hyperbolic problems 
do deal with multiscale wave structures resulting 
from shock wave interactions is also introduced. We 
present numerical results with realistic high-contrast 
two-dimensional multiscale coefficients with coupling 
multiscale oil-water flow simulations along with 
convincing numerical tests of local and global local 
mass conservation. We expect to combine the novel 
approach into the framework of Generalized Multiscale 
Finite Element Methods as recently introduced in 
\cite{ADG19}; see also \cite{AAP19} with particular 
interest to the case of complex flow systems as 
discussed in \cite{EA14,EADC13,ABFL19} for real-file 
applications, but in which issues of existence, 
stability properties and uniqueness are not well 
understood in line of works 
\cite{BADM15,EA14,ABFL19,CAFM16,AGP19,GCJG19,FKMT17,HCCZC19}.

\section*{Acknowledgments}
E. Abreu thanks research grants as well as thanks to all the support given by the Brazilian funding agencies FAPESP 2019/20991-8 (São Paulo), CNPq 306385 /2019-8 (National) and PETROBRAS 2015/00398-0 and 2019/00538-7. 
Juan Galvis thanks partial support  from the European Union's Horizon 2020 research and innovation programme under
the Marie Sklodowska-Curie grant agreement No 777778 (MATHROCKS).
\bibliographystyle{plain}

\begin{thebibliography}{150}

\bibitem{JEA04}
J. E. Aarnes, On the use of a mixed multiscale finite element
method for greater flexibility and increased speed or improved
accuracy in reservoir simulation. Multiscale Model. 
Simul. 2(3) (2004) 421-439.

\bibitem{ADG19}
E. Abreu, C. Díaz and J. Galvis. 
A convergence analysis of Generalized Multiscale 
Finite Element Methods, Journal of Computational 
Physics, 396(1) (2019) 303-324.

\bibitem{AbGaDiSa2018}
E. Abreu, C. Díaz, J. Galvis, and M. Sarkis.
On high-order conservative finite element methods, 
Computers \& Mathematics with Applications, 75 (2018) 1852-1867.


{\color{black}
\bibitem{AELP20}
E. Abreu, A. Esp\'{i}rito Santo, W. Lambert and J. P\'{e}rez,
Convergence of a Lagrangian-Eulerian scheme via the weak 
asymptotic method, submitted.}

\bibitem{BADM15}
B. Andreianov and D. Mitrović, Entropy conditions for 
scalar conservation laws with discontinuous flux revisited, 
Ann. I. H. Poincaré - AN 32 (2015) 1307-1335.

\bibitem{ACP16}
E. Abreu, M. Colombeau and E. Y. Panov, 
Weak asymptotic methods for scalar equations and systems. 
Journal of Mathematical Analysis and Applications, 
444 (2016) 1203-1232.

\bibitem{ACP17}
E. Abreu, M. Colombeau and E. Y. Panov, 
Approximation of entropy solutions to degenerate 
nonlinear parabolic equations, Zeitschrift f\"ur 
angewandte Mathematik und Physik, 68 (2017) 133. 

\bibitem{EA14}
E. Abreu,
Numerical modelling of three-phase immiscible flow 
in heterogeneous porous media with gravitational 
effects, Mathematics and Computers in Simulation,
97 (2014) 234-259.

\bibitem{BACC14}
B. Andreianov, C. Canc\'es. A phase-by-phase upstream scheme 
that converges to the vanishing capillarity solution for 
countercurrent two-phase flow in two-rocks media. Comput. 
Geosci., 18(2) (2014) 211-226

\bibitem{EADC13}
E. Abreu and D. Conceição. 
Numerical modeling of degenerate equations in porous 
media flow, Journal of Scientific Computing, 55(3) 
(2013) 688-717.

\bibitem{EAHYP12}
E. Abreu, 
Numerical simulation of wave propagation in three-phase 
flows in porous media with spatially varying flux functions. 
In: 4th International Conference on Hyperbolic Problems: 
Theory, Numerics, Applications, 2014, Padova/Itália. The 
proceedings of HYP2012. American Institute of Mathematical 
Sciences (AIMS); Series on Applied Mathematics 8 (2014) 233-240.

\bibitem{ABFL19}
E. Abreu, A. Bustos, P. Ferraz and W. Lambert. 
A relaxation projection analytical-numerical approach 
in hysteretic two-phase flows in porous media, Journal 
of Scientific Computing 79 (2019) 1936-1980.

\bibitem{CAFM16}
E. Abreu, P. Castañeda, F. Furtado and D. Marchesin. 
On a universal structure for immiscible three-phase 
flow in virgin reservoirs. Comput. Geosci. 20(1) 
(2016) 171-185.

\bibitem{AGP18a}
E. Abreu, J. Galvis, M. Sarkis and C. Penedo.
On high-order conservative finite element methods. 
Computers \& Mathematics with Applications 75 (2018) 1852-1867.

\bibitem{AGP18b}
E. Abreu, J. Galvis and C. Penedo. 
On high-order approximation and stability with conservative 
properties. Special DD24 Int. Conf. on Domain 
Decomposition Methods. Domain Decomposition Methods in 
Science and Engineering XXIV in Svalbard, Norway (Lecture 
Notes in Computational Science and Engineering 125) 1st 
ed. 2018. 2019. xxii, 570 S. 19 SW-Abb., 110 Farbabb. $235$ mm.

\bibitem{AGP19}
E. Abreu, J. Galvis and C. Penedo. 
A convergence analysis of Generalized Multiscale Finite 
Element Methods. Journal of Computational Physics, 
396 (2019) 303-324.

\bibitem{EAJP19}
E. Abreu and J. Perez. 
A fast, robust, and simple Lagrangian-Eulerian solver 
for balance laws and applications, Computers and 
Mathematics with Applications, 77(9) (2019) 2310-2336.

\bibitem{APS17a}
E. Abreu, W. Lambert, J. Perez and  A. Santo. 
A new finite volume approach for transport models 
and related applications with balancing source 
terms, Math. Comput. Simulation 137 (2017) 2-28.

\bibitem{APS17b}
E. Abreu, J. Pérez, A. Santo. 
A conservative Lagrangian-Eulerian finite volume 
approximation method for balance law problems, 
Proc. Ser. Braz. Soc. Comput. Appl. Math. 5 
(1) (2017) {\footnotesize 010329-1/010329-7}.

\bibitem{APS18}
E. Abreu, J. Pérez, A. Santo. 
Lagrangian-eulerian approximation methods for balance 
laws and hyperbolic conservation laws, Rev. UIS Ing. 
17 (1) (2018) 191-200, 
http://dx.doi.org/10.18273/revuin.v17n1-2018018, 
Proceeding of the XI Congreso Colombiano de Métodos 
Numéricos 2017.

\bibitem{AMPRB19} 
E. Abreu, V. Matos, J. Perez and P. Rodríguez-Bermúdez. 
A class of Lagrangian-Eulerian shock-capturing schemes 
for first-order hyperbolic problems with forcing terms, submitted.

\bibitem{AAP19}
E. Abreu, J. Agudelo and J. Pérez. 
Convergence and error estimates for a new conservative 
Lagrangian-Eulerian method on triangular grids for 
hyperbolic conservation laws, in preparation.

\bibitem{CNMAC2016} 
E. Abreu, A. Santo and J. P\'{e}rez,
Solving hyperbolic conservation laws by using Lagrangian-Eulerian approach
XXXVI Congresso Nacional de Matemática Aplicada e Computacional (CNMAC),
05 a 09 de setembro de 2016, Gramado - RS.
	
\bibitem{AMVG05}
M. Adimurthi, G.D. Veerappa Gowda, Optimal entropy 
solutions for conservation laws with discontinuous 
flux-functions. J. Hyperbolic Differ. Equ. 2(4), 
783-837, pp. 2005.

\bibitem{GA92}
G. Allaire. 
Homogenization and two-scale convergence. 
SIAM J. Math. Anal. 23(6) (1992) 1482-1518.

\bibitem{ABJ06}
B. Amaziane, A. Bourgeat and M. Jurak. 
Effective macrodiffusion in solute transport through 
heterogeneous porous media. Multiscale Model. 
Simul. 5(1) (2006) 184-204.

\bibitem{AKR11}
B. Andreianov, K. H. Karlsen and N. H. Risebro. A Theory of 
$L^1$-Dissipative Solvers for Scalar Conservation Laws with 
Discontinuous Flux, Archive for Rational Mechanics and 
Analysis, 201(1), 2011, pp. 27-86.

\bibitem{BA15}
B. Andreianov. New approaches to describing admissibility 
of solutions of scalar conservation laws with discontinuous 
flux, ESAIM: ESAIM: Proceedings and Surveys, 50 2015, pp. 40-65.

\bibitem{BACC13}
B. Andreianov and  C. Cancès. 
Vanishing capillarity solutions of Buckley-Leverett 
equation with gravity in two-rocks' medium, 
Computational Geosciences 17 (2013) 551-572.

\bibitem{AHPV13}
R. Araya, C. Harder, D. Paredes, F. Valentin. 
Multiscale hybrid-mixed method, 
SIAM Journal on Numerical Analysis 51 (6) (2013) 3505-3531 (2013).

\bibitem{ADH90}
T. Arbogast, J. Douglas and U. Hornung. 
Derivation of the double porosity model of single phase 
flow via homogenization theory. SIAM J. Math. Anal. 
21(4) (1990) 823-836.

\bibitem{APWY07}
T. Arbogast, G. Pencheva, M. F. Wheeler and I. Yotov.
A multiscale mortar mixed finite element method. 
Multiscale Model. Simul. 6(1) (2007) 319-346.

\bibitem{ASFM10} 
A. Azevedo, A. Souza, F. Furtado, D. Marchesin, B. Plohr, 
The solution by the wave curve method of three-phase flow in 
virgin reservoirs, Transport in Porous Media 83 (2010) 99-125. 

\bibitem{ASFM14}  
A. Azevedo, A. Souza, F. Furtado, D. Marchesin. Uniqueness of the 
Riemann solution for three-phase flow in a porous medium. SIAM 
Journal on Applied Mathematics, 74-76 (2014) 1709-1741.

\bibitem{AMPZ02} 
A. Azevedo, D. Marchesin, B.J. Plohr, K. Zumbrun, Capillary 
instability in models for three-phase flow, Zeitschrift fur 
Angewandte Mathematikund Physik 53 (2002) 713-746. 

\bibitem{ABSH01}
A. Y. Beliaev S. M. Hassanizadeh. 
A theoretical model of hysteresis and dynamic effects in 
the capillary relation for two-phase flow in porous media. 
Transp. Porous Media 43 (2001) 487-510.

\bibitem{MR2168342}
M. Benzi, G. H. Golub, and J. Liesen.
Numerical solution of saddle point problems, 
Acta Numer., 14 (2005) 1-137.

\bibitem{BBK04a}
S. Berres, R. Burger, K.H. Karlsen, Central schemes and systems 
of conservation laws with discontinuous coefficients modeling 
gravity separation of polydisperse suspensions, Journal of 
Computational and Applied Mathematics (Proceedings of the 10th 
International Congress on Computational and Applied Mathematics)
164-165(1) (2004) 53-80.


{\color{black}
\bibitem{BLIS15} 
T. Buckmaster, C. De Lellis, P. Isett and L. Székelyhidi, Jr, Anomalous dissipation for 1/5-H\"older Euler flows, Annals of Mathematics 182 (2015), 127-172.
}

\bibitem{BKRT04}
R. Burger, K.H. Karlsen, N.H Risebro, J.D. Towers,
Well-posedness in $BV_t$ and convergence of a difference scheme 
for continuous sedimentation in ideal clarifier-thickener units, 
Numerische Mathematik 97(1) (2004) 25-65.



\bibitem{FBBP98}
F. Bouchut, B. Perthame. 
Kružkov's estimates for scalar conservation laws 
revisited, Trans. Amer. Math. Soc., 350(7) (1998) 2847-2870.

\bibitem{AB84}
A. Bourgeat. 
Homogenized behavior of two-phase flows in naturally 
fractured reservoirs with uniform fractures distribution. 
Comput. Methods Appl. Mech. Eng. 47(1-2) (1984) 205-216.

\bibitem{BGJMP17}
K. Brenner, M. Groza, L. Jeannin, R. Masson and J. Pellerin. 
Immiscible two-phase Darcy flow model accounting for 
vanishing and discontinuous capillary pressures: 
application to the flow in fractured porous media, 
Computational Geosciences 21(5-6) (2017) 1075-1094.

\bibitem{CC10B}
C. Canc\'es, Asymptotic behavior of two-phase flows in heterogeneous 
porous media for capillarity depending only on space. II. 
Nonclassical shocks to model oil-trapping. SIAM J. Math. Anal. 
42(2) (2010) 972-995.

\bibitem{MCLD10}
M. A. Cardoso and L. J. Durlofsky. 
Linearized reduced-order models for subsurface 
flow simulation, Journal of Computational Physics 
229 (2010) 681-700.

\bibitem{CNP18}
X. Cao, S. F. Nemadjieu and I. S. Pop. 
Convergence of an MPFA finite volume scheme for a two 
phase porous media flow model with dynamic capillarity. 
IMA J. Numer. Anal. 39 (2018) 512-544.

\bibitem{HCCP19}
H. Carrillo and C. Parés. 
Compact Approximate Taylor Methods for Systems of 
Conservation Laws, Journal of Scientific Computing 
80(3) (2019) 1832-1866.

\bibitem{CCH99}
C. Chainais-Hillairet. 
Finite volume schemes for a nonlinear hyperbolic equation. 
Convergence towards the entropy solution and error estimate. 
Mathematical Modelling and Numerical Analysis M2AN 33(1) 
(1999) 129-156.

\bibitem{ICBP08}
I. Christov and B. Popov. 
New non-oscillatory central schemes on unstructured triangulations 
for hyperbolic systems of conservation laws. 
Journal of Computational Physics, 227(11) (2008) 5736-5757.

\bibitem{GCJG19}
G-Q G.Chen and J. Glimm, 
Kolmogorov-type theory of compressible turbulence 
and inviscid limit of the Navier-Stokes equations 
in $\mathbf{R^3}$, Physica D: Nonlinear Phenomena 
400(15) (2019) 132138.

\bibitem{GGJG12}
G. Q. Chen and J. Glimm. 
Kolmogorov’s theory of turbulence and inviscid limit 
of the Navier-Stokes equations in R 3 . Comm. Math. 
Phys. 310 (1) (2012) 267-283.

\bibitem{CHM06}
Z. Chen, G. Huan and Y. Ma. 
Computational Methods for Multiphase Flows in Porous 
Media. Vol. 2. Philadelphia: SIAM (2006).

\bibitem{CDGW03}
Y. Chen, L. J. Durlofsky, M. Gerritsen, X.-H Wen, A coupled
local-global upscaling approach for simulating flow in highly
heterogeneous formations. Adv. Water Resour. 26(10) (2003) 1041-1060.

\bibitem{ECCL15}
E. Chiodaroli, C. De Lellis and O. Kreml. 
Global ill-posedness of the isentropic system of 
gas dynamics. Comm. Pure Appl. Math., 68(7) 
(2015) 1157-1190.

\bibitem{CLVW}
E.T. Chung, W. T. Leung, M. Vasilyeva, Y. Wang. 
Multiscale model reduction for transport and flow problems 
in perforated domains. J. Comput. Appl. Math. 330 (2018) 519-535.

\bibitem{TenhSPE}
M. Christie and M. Blunt. 
Tenth spe comparative solution project: A comparison 
of upscaling techniques, Society of Petroleum Engineers
SPE Reservoir Simulation Symposium, 11-14 February (2001) 
Houston, Texas, SPE-66599-MS.

\bibitem{CRMAJ80}
M. G. Crandall and A. Majda.
Monotone difference approximations for scalar conservation 
laws. Mathematics of Computation, 34(149) (1980) 1-21.

\bibitem{CCP15}
G. Crasta, V. De Cicco and G. De Philippis.
Kinetic Formulation and Uniqueness for Scalar Conservation 
Laws with Discontinuous Flux, Communications in Partial 
Differential Equations, 40(4) 2015, pp. 694-726.

\bibitem{SD96}
S. Diehl. A conservation law with point source and discontinuous 
flux function modelling continuous sedimentation. SIAM J. Appl. 
Math., 56(2) (1996) 388-419.

\bibitem{RJD85}
R. J. DiPerna. Measure valued solutions to conservation 
laws. Arch. Rational Mech. Anal., 88(3) (1985) 223-270.

\bibitem{RDAM87}
R. J. DiPerna and A. Majda. 
Oscillations and concentrations in weak solutions of 
the incompressible fluid equations. Comm. Math. Phys. 
108 (4) (1987) 667-689.

\bibitem{LJD91}
L. J. Durlofsky. 
Numerical calculation of equivalent grid block
permeability tensors for heterogeneous porous media. Water
Resour. Res. 27(5) (1991) 699-708.

\bibitem{LDHL16}
L. J. Durlofsky and H. Li.
Local-global upscaling for compositional subsurface flow 
simulation. Transp. Porous Media 111(3) (2016) 701-730.

\bibitem{EGH95a}
R. Eymard, T. Gallouët and R. Herbin. 
Existence and uniqueness of the entropy solution 
to a nonlinear hyperbolic equation. Chin. Ann. 
Math. B16 (1995) 1-14.

\bibitem{EGH95b}
R. Eymard, T. Gallouët and R. Herbin. 
Convergence of a finite volume scheme for a 
nonlinear hyperbolic equation, Proceedings 
of the Third colloquium on numerical analysis, 
edited by D. Bainov and V. Covachev. Elsevier (1995) 61-70.

\bibitem{FSYFY19}
W. Fan, H. Sun,  J. Yao, D. Fan and Y. Yang. 
An upscaled transport model for shale gas 
considering multiple mechanisms and heterogeneity 
based on homogenization theory. Journal of 
Petroleum Science and Engineering, 106392 (2019) 
doi:10.1016/j.petrol.2019.106392 

\bibitem{FKMT17} 
U. S. Fjordholm, R. Käppeli, S. Mishra, and E. Tadmor. 
Construction of approximate entropy measure-valued 
solutions for hyperbolic systems of conservation laws. 
Found. Comput. Math., 17(3) (2017) 763-827, 2017.

\bibitem{FLMW19}
U. S. Fjordholm, K. Lye, S. Mishra, F. Weber, 
Statistical solutions of hyperbolic systems of 
conservation laws: numerical approximation, 
https://arxiv.org/abs/1906.02536 (access in Oct 13, 2019).

\bibitem{HFVS03}   
H. Frid, V. Shelukhin. A Quasi-linear Parabolic System for 
Three-Phase Capillary Flow in Porous Media, SIAM J. Math. Anal., 
35(4) (2003) 1029-1041.

\bibitem {HFVS05}  
H. Frid, V. Shelukhin. Initial Boundary Value Problems for 
a Quasi-linear Parabolic System in Three-Phase Capillary Flow 
in Porous Media, SIAM J. Math. Anal., 36(5) (2005) 1407-1425.

\bibitem{HFEG19}
H. Florez and E. Gildin. 
Model-Order Reduction of Coupled Flow and Geomechanics 
in Ultra-Low Permeability ULP Reservoirs, SPE-193911-MS,  
SPE Reservoir Simulation Conference, 10-11 April, 
Galveston, Texas, USA 2019.

\bibitem{UFTu95}
U. Frisch. Turbulence. Cambridge University Press, 1995.

\bibitem{FFFP03}
F. Furtado and F. Pereira. 
Crossover from Nonlinearity Controlled to Heterogeneity 
Controlled Mixing in Two-Phase Porous Media Flows, 
Computational Geosciences 7(2) (2003) 115-13.


\bibitem{GMSFEM}
Efendiev, Yalchin, Juan Galvis, and Thomas Y. Hou. 
Generalized multiscale finite element methods (GMsFEM),
 Journal of Computational Physics 251 (2013): 116-135.

\bibitem{PRESHO2016376}
J. Galvis and M. Presho. 
A mass conservative generalized multiscale finite 
element method applied to two-phase flow in 
heterogeneous porous media, Journal of 
Computational and Applied Mathematics, 296 (2016) 376-388.

\bibitem{GCLW19}
W. Ge, Q. Chang, C. Li and J. Wang. 
Multiscale structures in particle-fluid systems: 
Characterization, modeling, and simulation, 
Chemical Engineering Science 198(28) (2019) 198-223.

\bibitem{GKNSVL19}
D. Groen, J. Knap, P. Neumann, D. Suleimenova, L. Veen and K. Leiter. 
Mastering the scales: a survey on the benefits of multiscale 
computing software.Phil. Trans. R. Soc. A (2019) 377: 20180147.

\bibitem{GKD06}
B. Gong, M. Karimi-Fard and L. J. Durlofsky. 
An upscaling procedure for constructing generalized 
dual-porosity/dual-permeability models from discrete 
fracture characterizations. SPE Annual Technical 
Conference and Exhibition. Society of Petroleum 
Engineers (2006).

\bibitem{GASP18}
R. T. Guiraldello, R. F. Ausas, F. S. Sousa, F. Pereira and G. C. Buscaglia. 
The multiscale Robin coupled method for flows in porous media. 
J. Comput. Phys. 355 (2018) 1-21.

\bibitem{HCCZC19}
A. G. Hoekstra , B. Chopard , D. Coster , S. P. Zwart  and P. V. Coveney. 
Multiscale computing for science and engineering in the era of 
exascale performance, Phil. Trans. R. Soc. A (2019) 377: 20180144.

\bibitem{HHNR95}
H. Holden, N. H. Risebro. A mathematical model of traffic 
flow on a network of unidirectional roads. SIAM J. Math. Anal., 
26(4) (1995) 999-1017.

\bibitem{THXW97}
T. Y. Hou and X.-H. Wu. 
A multiscale finite element method for elliptic problems 
in composite materials and porous media.
J. Comput. Phys. 134(1) (1997) 169-189.

\bibitem{JJLD17}
J. D. Jansen and L. J. Durlofsky. 
Use of reduced-order models in well control 
optimization, Optim Eng 18 (2017) 105-132.

\bibitem{JLT03}
P. Jenny, S.H. Lee and H. A. Tchelepi.
Multi-scale finite-volume method for elliptic 
problems in subsurface flow simulation. 
J. Comput. Phys. 187(1) (2003) 47-67.

\bibitem{GSJET98} 
G-S Jiang and E. Tadmor. 
Nonoscillatory central schemes for multidimensional hyperbolic 
conservation laws. 
SIAM Journal on Scientific Computing, 19(6):1892–1917, 1998.

\bibitem{JKO12}
V. V. Jikov, S. M. Kozlov and A. O. Oleinik. 
Homogenization of differential operators and integral 
functionals. Springer, Berlin Heidelberg (2012).

\bibitem{KEF99}
E. F. Kaasschieter, Solving the Buckley-Leverett equation 
with gravity in a heterogeneous porous medium. Comput. Geosci. 
3(1) (1999) 23-48.

\bibitem{KT04}
K. H. Karlsen, J. D. Towers, Convergence of the Lax-Friedrichs 
scheme and stability for conservation laws with a discontinuous 
space-time dependent flux. Chin. Ann. Math. Ser. B. 25(3) (2004) 
287-318.

\bibitem{KHKNHR02}
K. H. Karlsen, N. H. Risebro. Unconditionally 
stable methods for Hamilton-Jacobi equations. 
J. Comp. Phys. 180(2), 2002, pp. 710-735.

\bibitem{NNK76}
N. N. Kuznetsov. 
Accuracy of some approximate methods for computing 
the weak solutions of a first-order quasi-linear 
equation, USSR Comput. Math. and Math. Phys., 
16 (1976), pp. 105-119.


\bibitem{LWT08}
S. H. Lee, C. Wolfsteiner and H. A. Tchelepi. 
Multiscale finite volume formulation for multiphase 
flow in porous media: black oil formulation of 
compressible, three-phase flow with gravity.
Comput. Geosci. 12(3) (2008) 351-366.

\bibitem{CLLS09}
C. De Lellis and L. Szekelyhidi, Jr. 
The Euler equations as a differential inclusion. 
Ann. of Math. (2), 170(3) (2009) 1417-1436.

\bibitem{LVFEB02}
R. J. LeVeque.
{Finite volume methods for hyperbolic problems} 31.
Cambridge university press, 2002.
	
\bibitem{LGS08}
H. Lim, Y. Yu, J. Glimm, X. L. Li and D. H. Sharp. 
Chaos, transport and mesh convergence for fluid mixing. 
Act. Math. Appl. Sin., 24 (3) (2008) 355-368.

\bibitem{DMBP01} 
D. Marchesin, B. Plohr. Wave structure in wag recovery, 
Society of Petroleum Engineering Journal 71314 (2001) 209-219. 

\bibitem{MDD06}
A. Mikelić and V. Devigne and C. J. Van Duijn. 
Rigorous upscaling of the reactive flow through a pore, 
under dominant peclet and damkohler numbers. 
SIAM J. Math. Anal. 38(4) (2006) 1262-1287.

\bibitem{SMJJ10}
S. Mishra, J. Jaffr\'e. On the upstream mobility scheme 
for two-phase flow in porous e media. Comput. Geosci., 
14(1) (2010) 105-124.

\bibitem{OMKAL16a}
O. Moyner and K. A. Lie. 
A multiscale restriction-smoothed basis method for 
high contrast porous media represented on unstructured 
grids. J. Comput. Phys. 304 (2016) 46-71.

\bibitem{OMKAL16b}
O. Moyner and K. A. Lie. 
A multiscale restriction-smoothed basis method for 
compressible black-oil model. Society of Petroleum 
Engineers Journal 21(6) (2016) 2079-2096.

\bibitem{DNO99}
D. N. Ostrov. Viscosity solutions and convergence of monotone 
schemes for synthetic aperture radar shape-from-shading equations 
with discontinuous intensities. SIAM J. Appl. Math., 59(6) 
(1999) 2060-2085.

\bibitem{EP05}
E. Yu. Panov, 
Existence of strong traces for generalized solutions 
of multidimensional scalar conservation laws, 
J. Hyperbolic Differ. Equ. 2 (2005) 885-908.

\bibitem{EP07}
E. Yu. Panov, 
Existence of strong traces for quasi-solutions of 
multidimensional conservation laws, J. Hyperbolic 
Differ. Equ. 4 (2007) 729-770.

\bibitem{EP09}
E.Yu. Panov, 
On existence and uniqueness of entropy solutions to 
the Cauchy problem for a conservation law with 
discontinuous flux, J. Hyperbolic Differ. Equ. 3 
(2009) 525-548.

\bibitem{EP10}
E. Yu. Panov, 
Existence and strong pre-compactness properties for 
entropy solutions of a first-order quasilinear equation 
with discontinuous flux, Arch. Ration. Mech. Anal. 
195 (2010) 643-673.

\bibitem{RKNP18}
F. A. Radu, K. Kumar, J. M Nordbotten and I. S. Pop. 
A robust, mass conservative scheme for two-phase flow 
in porous media including Hölder continuous nonlinearities. 
IMA J. Numer. Anal. 38(2) (2018) 884-820.

\bibitem{RPK08}
F. A. Radu, I. S. Pop and P. Knabner.
Error estimates for a mixed finite element discretization 
of some degenerate parabolic equations. Numer. Math. 
109 (2008) 285-311.

\bibitem{RSS19}
A. M. Ruf, E. Sande and S. Solem. 
The Optimal Convergence Rate of Monotone Schemes
for Conservation Laws in the Wasserstein Distance, 
Journal of Scientific Computing 80 (2019) 1764-1776.

\bibitem{FB06}
F. Sabac. 
The Optimal Convergence Rate of Monotone Finite 
Difference Methods for Hyperbolic Conservation Laws, 
SIAM J. Numer. Anal., 34(6) (2006) 2306-2318. 

\bibitem{SV03}
N. Seguin and J. Vovelle.
Analysis and approximation of a scalar conservation 
law with a flux function with discontinuous coefficients, 
Math. Models Methods Appl. Sci. 13 (2) (2003) 221-257.

\bibitem{SLW19}
G. Singh, W. Leung and M. F. Wheeler. 
Multiscale methods for model order reduction of 
non-linear multiphase flow problems, 
Computational Geosciences 23 (2019) 305-323.

\bibitem{TT01}
T. Tang. 
Error Estimates of Approximate Solutions for Nonlinear 
Scalar Conservation Laws. In: Freistühler H., Warnecke 
G. (eds) Hyperbolic Problems: Theory, Numerics, 
Applications. ISNM International Series of Numerical 
Mathematics, 141. Birkhäuser, Basel (2001).

\bibitem{ATMV16}
A. Talonov and M. Vasilyeva. 
On numerical homogenization of shale gas 
transport. J. Comput. Appl. Math. 301 (2016) 44-52.

\bibitem{TTZT06}
T. Tang and Z-h Teng. 
On the Regularity of Approximate Solutions to 
Conservation Laws with Piecewise Smooth Solutions, 
SIAM J. Numer. Anal., 38(5) (2006) 1483-1495.

\bibitem{THT19}
Z-X. Tong,  Y-L. He and W-Q. Tao. 
A review of current progress in multiscale simulations 
for fluid flow and heat transfer problems: The frameworks, 
coupling techniques and future perspectives, International 
Journal of Heat and Mass Transfer, 137 (2019) 1263-1289.

\bibitem{STIA12}
S. Tveit, I. Aavatsmark. Errors in the upstream mobility 
scheme for countercurrent two-phase flow in heterogeneous 
porous media, Comput Geosci 16 (2012) 809-825.

\bibitem{VCEK19}
M. Vasilyeva, E. T. Chung, Y. Efendiev and J. Kim. 
Constrained energy minimization based upscaling for 
coupled flow and mechanics, 
Journal of Computational Physics 376(1) (2019) 660-674.

\bibitem{WEH02}
X.-H. Wu, Y. Efendiev, T. Y. Hou. 
Analysis of upscaling absolute permeability. 
Discrete and Continuous Dynamical Systems Series
B 2(2) (2002) 185-204.


\bibitem{GMsFEMrecent}
Wang, M., Cheung, S. W., Chung, E. T., Efendiev, Y., Leung, W. T.,  Wang, Y. (2019). Prediction of Discretization of GMsFEM Using Deep Learning. Mathematics, 7(5), 412.

\bibitem{EGW2009jcp}
Efendiev, Y., Galvis, J.,  Wu, X. H. (2009). Multiscale finite element and domain decomposition methods for high-contrast problems using local spectral basis functions.


\end{thebibliography}

\end{document}